\documentclass[11pt]{article}
\usepackage{graphicx}
\usepackage{amsmath}
\usepackage{amsthm}
\usepackage{amssymb}
\usepackage{latexsym}
\usepackage{mathrsfs}
\usepackage{textcase}
\usepackage[dvips]{epsfig}
\usepackage[small]{caption}
\usepackage[colorlinks=true,linkcolor=blue,citecolor=red]{hyperref}

\newtheoremstyle{dotless}{}{}{\itshape}{}{\bfseries}{}{ }{}
\theoremstyle{dotless}

\newenvironment{Proof}[1][Proof]%
 {\begin{trivlist} \item[]{\textbf{#1.} }}%
{\hspace*{\fill}$\rule{.4\baselineskip}{.4\baselineskip}$\end{trivlist}}

\newenvironment{Acknowledgment}
 {\begin{trivlist}\item[]\textbf{Acknowledgments }}{\end{trivlist}}

\numberwithin{equation}{section}
\newtheorem{lemma}{Lemma}[section]
\newtheorem{theorem}{Theorem}
\newtheorem{proposition}[lemma]{Proposition}
\newtheorem{definition}[lemma]{Definition}
\newtheorem{corollary}[lemma]{Corollary}
\newtheorem{hypothesis}[lemma]{Hypothesis}
\newtheorem{hypotheses}[lemma]{Hypotheses}
\newtheorem{remark}[lemma]{Remark}

\newcommand{\bt}{\begin{theorem}}
\newcommand{\et}{\end{theorem}}
\newcommand{\bl}{\begin{lemma}}
\newcommand{\el}{\end{lemma}}
\newcommand{\bp}{\begin{proposition}}
\newcommand{\ep}{\end{proposition}}
\newcommand{\bd}{\begin{definition}}
\newcommand{\ed}{\end{definition}}
\newcommand{\bc}{\begin{corollary}}
\newcommand{\ec}{\end{corollary}}
\newcommand{\br}{\begin{remark}}
\newcommand{\er}{\end{remark}}
\newcommand{\bh}{\begin{hypotheses}}
\newcommand{\eh}{\end{hypotheses}}
\newcommand{\be}{\begin{enumerate}}
\newcommand{\ee}{\end{enumerate}}
\newcommand{\beq}{\begin{equation}}
\newcommand{\eeq}{\end{equation}}
\newcommand{\beqs}{\begin{equation*}}
\newcommand{\eeqs}{\end{equation*}}
\newcommand{\bpf}{\begin{Proof}}
\newcommand{\epf}{\end{Proof}}
\newcommand{\bld}{\begin{aligned}}
\newcommand{\eld}{\end{aligned}}
\newcommand{\bhp}{\begin{hypothesis}}
\newcommand{\ehp}{\end{hypothesis}}
\newcommand{\bcs}{\begin{cases}}
\newcommand{\ecs}{\end{cases}}
\newcommand{\R}{\mathbb{R}}
\newcommand{\C}{\mathbb{C}}
\newcommand{\N}{\mathbb{N}}
\newcommand{\Z}{\mathbb{Z}}

\newcommand{\T}{\mathbb{T}}

\newcommand{\rmO}{\mathrm{O}}

\newcommand{\rmd}{\mathrm{d}}
\newcommand{\rme}{\mathrm{e}}
\newcommand{\rmi}{\mathrm{i}}
\newcommand{\id}{\mathrm{\,id}\,}

\renewcommand{\leq}{\leqslant}
\renewcommand{\geq}{\geqslant}

\newcommand{\rmnum}[1]{\romannumeral #1}
\newcommand{\Rmnum}[1]{\uppercase\expandafter{\romannumeral #1\relax}}

\def\re{\mathop{\mathrm{Re}}}

\def\span{\mathop{\mathrm{span}}}
\def\spec{\mathop{\mathrm{spec}}}

\def\rg{\mathop{\mathrm{Rg}}}

\makeatletter
\newsavebox{\@brx}
\newcommand{\llangle}[1][]{\savebox{\@brx}{\(\m@th{#1\langle}\)}%
  \mathopen{\copy\@brx\kern-0.5\wd\@brx\usebox{\@brx}}}
\newcommand{\rrangle}[1][]{\savebox{\@brx}{\(\m@th{#1\rangle}\)}%
  \mathclose{\copy\@brx\kern-0.5\wd\@brx\usebox{\@brx}}}
\makeatother

\makeatletter
\newcommand{\opnorm}{\@ifstar\@opnorms\@opnorm}
\newcommand{\@opnorms}[1]{%
  \left|\mkern-1.5mu\left|\mkern-1.5mu\left|
   #1
  \right|\mkern-1.5mu\right|\mkern-1.5mu\right|
}
\newcommand{\@opnorm}[2][]{%
  \mathopen{#1|\mkern-1.5mu#1|\mkern-1.5mu#1|}
  #2
  \mathclose{#1|\mkern-1.5mu#1|\mkern-1.5mu#1|}
}
\makeatother

\makeatletter\@addtoreset{figure}{section}\makeatother

\newfam\bifam
\font\tenbi=cmmib10 scaled \magstep1 \font\sevenbi=cmmib10 at 11pt
\font\fivebi=cmmib10 at 6pt \textfont\bifam = \tenbi
\scriptfont\bifam = \sevenbi \scriptscriptfont\bifam= \fivebi

\title{\textbf{Diffusive stability of Turing patterns via normal forms} }
\author{Arnd Scheel and  Qiliang Wu\\
School of Mathematics\\
University of Minnesota\\
206 Church St SE\\
Minneapolis, MN 55414}
\date{}

\begin{document}

\maketitle

\begin{abstract}
We investigate dynamics near Turing patterns in reaction-diffusion systems posed on the real line.
Linear analysis predicts diffusive decay of small perturbations. We construct a ``normal form'' coordinate system near
such Turing patterns which exhibits an approximate discrete conservation law. The key ingredients to the normal form is
a conjugation of the reaction-diffusion system on the real line to a lattice dynamical system. At each lattice site, we
decompose perturbations into neutral phase shifts and normal decaying components. As an application of our normal form construction,
we prove nonlinear stability of Turing patterns with respect to perturbations that are small in $L^1\cap L^\infty$, with sharp rates,
recovering and slightly improving on results in \cite{schneider_1996,johnsonzumbrun_2011}.

%
%
%
\end{abstract}

\vfill

\hrule
{\small
\begin{Acknowledgment}
This work was partially supported by the National Science
Foundation through grant NSF-DMS-0806614.
\end{Acknowledgment}}

\section{Introduction}

Turing predicted that the simple interplay of reaction and diffusion can lead to stable, spatially periodic patterns \cite{turing_1952}.
His ideas proved quite influential in the general area of pattern formation,
where one seeks to understand the formation and dynamics of self-organized spatio-temporal structures.
One can easily envision simple reaction-diffusion systems with two species that exhibit diffusion-driven instabilities of
spatially homogeneous equilibria. Typical examples are activator-inhibitor systems such as the Gray-Scott or the Gierer-Meinhard equation;
see for instance \cite{murray1,  murray2}.
Perturbations of the homogeneous unstable equilibrium grow exponentially at an initial stage,
with fastest growth for distinct spatial wavenumbers. This wavenumber is roughly independent of boundary conditions
in large enough domains. As a final result, one often finds a spatially periodic pattern, up to narrow,
exponentially localized boundary layers. In order to understand such nonlinear spatially periodic patterns and the process of wavenumber selection,
one is therefore naturally led to considering reaction-diffusion systems on idealized unbouned domains.

To fix ideas, consider
\begin{equation*}
\mathbf{u}_t=D\triangle\mathbf{u} + \mathbf{f}(\mathbf{u}),
\end{equation*}
for $\mathbf{u}(t,x)\in\R^n$, with $x\in\R^N$, with smooth reaction-kinetics $\mathbf{f}$ and positive diagonal diffusion
matrix $D=\mathrm{diag}\,(d_j)>0$. Here, and in the following, the term ``smooth'' refers to functions
with sufficiently many derivatives.
In many circumstances, one can show that there exist families of spatially periodic striped solutions,
\begin{equation*}
\mathbf{u}(t,x)=\mathbf{u}_\star(k x_1;k),\qquad \mathbf{u}_\star(\xi;k)=\mathbf{u}_\star(\xi+2\pi;k),
\end{equation*}
parameterized by the spatial wavenumber $k>0$. In fact, such families occur for an open class of reaction-diffusion systems,
including but not limited to systems of activator-inhibitor type mentioned above.

As a first predictor on the stability of such solutions with respect to perturbations, one analyzes the linearization,
\begin{equation}\label{e:glrd}
\mathbf{v}_t=D\triangle\mathbf{v} + \mathbf{f}'(\mathbf{u}_\star(kx;k))\mathbf{v}.
\end{equation}
It turns out that, again for open classes of reaction-diffusion systems including the above examples,
solutions to this linear equation are bounded for bounded initial data,
for an open subset of patterns $\mathbf{u}_\star(\cdot;k)$ in the family. We refer to such patterns as \emph{linearly stable Turing patterns}.
We will discuss detailed assumptions that guarantee such linear stability later in this section.

The presence of a family of patterns, parameterized by the wavenumber, and, even more obviously, by translations of the pattern in $x$,
implies that solutions to (\ref{e:glrd}) with general initial conditions will not decay.
More explicitly, $\mathbf{v}(t,x)=\partial_x \mathbf{u}_\star(kx;k)$ and
$\mathbf{v}(t,x)=\frac{\rmd}{\rmd k} \mathbf{u}_\star(kx;k)$  are constant in time and solve (\ref{e:glrd}).

In fact, one can show that under typical assumptions, initial conditions $\mathbf{v}(t=0,x)\in L^1(\R^N,\R^n)$ will give rise to diffusive decay,
$\sup_x |\mathbf{v}(t,x)|\leq C t^{-N/2}$. Such algebraic decay is in general not strong enough to ensure nonlinear decay in dimensions $N\leq 3$.
The simplest example is the nonlinear heat equation
\beqs
u_t=\triangle u + u^2,
\eeqs
which exhibits blowup of arbitrarily small, smooth, positive initial data at finite time in dimensions $N\leq 3$ \cite{hervel_1993,denglevine_2000}.
In the seminal paper \cite{schneider_1996}, Schneider recognized that diffusive decay near Turing patterns is not altered
by the presence of nonlinear terms due to cancellations in a Bloch-wave expansion. He studied the most difficult case, $N=1$,
where diffusion is weak and nonlinearity potentially most dangerous, in the specific example of the Swift-Hohenberg equation.
His proof has later been generalized, simplified, and adapted; see \cite{uecker_1999,johnson_2009, johnsonzumbrun_2010, johnsonzumbrunpascal_2011,
johnsonzumbrun_2011sj, johnsonzumbrun_2011, gallayscheel_2011,sssu_2012}.
Our focus here is, again, on the one-dimensional case, in a general reaction-diffusion setting.
Our goal is to find coordinates that show explicitly  why nonlinear terms do not alter linear decay near Turing patterns.
Going back to the scalar heat equation,  the interaction of nonlinear terms with diffusion can be categorized as relevant,
critical, or irrelevant; \cite{bricmontkupiainen_1992,bricmontkupiainen_1994}. Explicitly,
in the heat equation $u_t=u_{xx}+f(u,u_x,u_{xx})$,
\begin{itemize}
\item[(\rmnum{1})] Nonlinear terms such as $f(u,u_x,u_{xx})=u u_{xx},u_x^2, u^p$, where $p>3$ are irrelevant;
\item[(\rmnum{2})] Nonlinear terms such as $f(u,u_x,u_{xx})=u u_x, u^3$ are critical;
\item[(\rmnum{3})] Nonlinear terms such as $f(u)=u^2$ are relevant.
\end{itemize}
Without pretending to fully explain this phenomenon, notice that, for $L^1$-initial data, assuming Gaussian decay,
we find $u_{xx}\sim t^{-3/2}$ in $L^\infty$. Irrelevant nonlinear terms decay with rate $t^{-\alpha}$, $\alpha>3/2$,
critical terms have $\alpha=3/2$, and relevant terms have $\alpha<3/2$.

Perturbations $\mathbf{v}$ of Turing patterns solve a system
\begin{equation*}
\mathbf{v}_t=\partial_{xx}\mathbf{v} +\mathbf{f}'(\mathbf{u}_\star(x))\mathbf{v}+\mathbf{g}(x,\mathbf{v}),
\end{equation*}
where $\mathbf{g}(x,\mathbf{v})=\rmO(|\mathbf{v}|^2)$. Note that from here on, we fix the wavenumber $k=1$, without loss of generality,
and write $\mathbf{u}_\star(x):=\mathbf{u}_\star(x;1)$. In particular, the nonlinearity $\mathbf{g}$ has potentially dangerous quadratic terms.
Roughly speaking, our goal is to find coordinates in which the nonlinearity involves at least two ``derivatives'',
which according to the numerology for the scalar heat equation would be sufficient to guarantee nonlinear decay.
The reason to hope for derivatives is the presence of a conservation law associated with the translation symmetry,
which in turn generates the neutral decay in the linearization.

To be precise, we now consider reaction diffusion systems
\beq
\label{e:21}
\mathbf{u}_t=D\partial_{xx} \mathbf{u}+\mathbf{f}(\mathbf{u}),
\eeq
where $\mathbf{u},\mathbf{f}\in\mathbb{R}^n, x\in \mathbb{R}, t\in(0,+\infty), D\in \mathbb{R}^{n\times n}$
is a diagonal matrix with strictly positive diagonal entries and $\mathbf{f}$ is smooth.
Firstly, we assume the existence of a Turing pattern of the system.
\begin{hypothesis}[\textbf{existence}]
\label{h:21}
The system of ordinary differential equations $D\partial_{xx} \mathbf{u}+\mathbf{f}(\mathbf{u})=0$
possesses a smooth periodic even solution $\mathbf{u}_{\star}$.
\end{hypothesis}
Without loss of generality, we assume that the period is $2\pi$. Our aim is to study nonlinear stability of
this temporal equilibrium under general small non-periodic perturbations.
To this end, we introduce an initial condition
\beq
\label{e:22}
\mathbf{u}(0,x)=\mathbf{u}_\star(x)+\mathbf{v}^0(x).
\eeq
Then assuming that $\mathbf{u}(t,x)=\mathbf{u}_{\star}(x)+\mathbf{v}(t,x)$ is a solution to (\ref{e:21})
with the given initial condition (\ref{e:22}), we have
\begin{equation}
\label{e:23}
\bcs
\mathbf{v}_t=A\mathbf{v}+\mathbf{g}(x,\mathbf{v}),\\
\mathbf{v}(0)=\mathbf{v}^0,
\ecs
\end{equation}
where
\begin{equation}
\label{e:A}
 \begin{matrix}
 A:&X^1& \longrightarrow & X\\
 &\mathbf{v}&\longmapsto&D\partial_{xx} \mathbf{v}+\mathbf{f}^\prime(\mathbf{u}_{\star})\mathbf{v}.
\end{matrix}
\end{equation}
Here we define
\beq
\label{e:XX1}
X=(L^1(\R))^n\cap (L^\infty(\R))^n,\quad X^1=(W^{2,1}(\R))^n\cap (W^{2,\infty}(\R))^n,
\eeq
with norms
\beqs
\|\cdot\|_X=\|\cdot\|_{L^1}+\|\cdot\|_{L^\infty},\quad \|\cdot\|_{X^1}=\|\cdot\|_{W^{2,1}}+\|\cdot\|_{W^{2,\infty}}.
\eeqs
Note that from now on, we suppress $n$ and $\R$ if there is no ambiguity.
Moreover,
$\mathbf{g}:\mathbb{T}_{2\pi}\times \mathbb{R}^n\rightarrow \mathbb{R}^n $ is smooth,
$\mathbf{g}(x,\mathbf{v})=\mathbf{f}(\mathbf{u}_\star+\mathbf{v}(x))-\mathbf{f}(\mathbf{u}_\star)-
\mathbf{f}^\prime(\mathbf{u}_\star)\mathbf{v}$, so that
$\mathbf{g}(x,0)\equiv 0$ and $\partial_{\mathbf{v}}\mathbf{g}(x,0)\equiv 0$.

According to Bloch wave decomposition, let us introduce the family of Bloch operators, for $\sigma\in[-\frac{1}{2},\frac{1}{2}]$,
\begin{equation}
\label{e:B}
\begin{matrix}
 B(\sigma):
&(H^2(\mathbb{T}_{2\pi}))^n &\longrightarrow & (L^2(\mathbb{T}_{2\pi}))^n\\
& \mathbf{v} &\longmapsto & D(\partial_x+\rmi \sigma)^2 \mathbf{v}+\mathbf{f}^\prime(\mathbf{u}_{\star})\mathbf{v}.
\end{matrix}
\end{equation}
 For further reading on Bloch wave decomposition and Bloch operators,
 we refer to Section \ref{ss:62} and \cite{reedsimon}. Note that one obtains $B(\sigma)$ formally by applying $A$ to functions of the form
 $\mathbf{u}=\rme^{\rmi \sigma x}\mathbf{v}$.
\begin{hypothesis}[\textbf{spectral stability}]
\label{h:22}
The family of Bloch wave operators $B(\sigma)$ has the following properties.
\be
 \item [(\rmnum{1})]$\spec(B(\sigma))\bigcap\{Re\lambda\geq0\}=\emptyset$, for $\sigma\neq 0$;
 \item [(\rmnum{2})]$\spec(B(0))\bigcap\{Re\lambda\geq0\}=\{0\}$ and $0$ is simple with $\span\{\mathbf{u}_{\star}^\prime\}$ as its eigenspace;
 \item [(\rmnum{3})]Near $\sigma=0$, the only eigenvalue $\lambda$ is a smooth function of $\sigma$ and the expression of $\lambda(\sigma)$ reads:
 $\lambda(\sigma)=-d\sigma^2+\rmO(|\sigma|^3)$, where $d>0$ is a constant.
\ee
\end{hypothesis}
\br
The expansion in (\rmnum{3}) is a consequence of the simplicity of
$\lambda=0$ at $\sigma=0$ and the evenness of $\mathbf{u}_\star$.
In fact, we have an ``explicit'' expression for $d$; see Section \ref{ss:65}.
\er

Given the above hypotheses, we can state our main result.
\bt[\textbf{nonlinear stability}]
\label{t:1}
Assume Hpotheses \ref{h:21} and \ref{h:22} hold. There are $C,\sigma>0$ so that, for any
$\|\mathbf{v}^0\|_{X}<\sigma$, where $X=(L^1(\R))^n\cap (L^\infty(\R))^n$, the solution $\mathbf{v}(t)$ to the system \eqref{e:23}
exists for time $t\in[0,\infty)$ and satisfies the estimate
\beq
\|\mathbf{v}(t)\|_{(L^\infty(\R))^n}\leq C\frac{\|\mathbf{v}^0\|_{X}}{(1+t)^{\frac{1}{2}}}.
\eeq
\et

The rest of the paper contains three main contributions. First, we construct normal form coordinates,
where the neutral mode is represented by a discrete phase $\theta_j$, which decays according to a linear
discrete diffusion equation $\dot{\theta}_j=d(\theta_{j+1}-2\theta_{j}+\theta_{j-1})$.
The idea is to capture the leading order dynamics of perturbations using an ansatz of the
type $\mathbf{u}(t,x)=\mathbf{u}_\star(x-\theta_j)+\mathbf{w}_j(t,x)$ on intervals $x\in[2\pi (j-1/2),2\pi(j+1/2)]$,
where $\mathbf{w}_j(t,x)$ lies in a linear strong stable fiber. The coordinate change mimics the much
simpler coordinate change in \cite{gallayscheel_2011}, where strong stable fibers of a \emph{temporally periodic},
but spatially homogeneous solution were straightened out.

Our second main contribution are decay estimates for the linearization in these coordinates.
In particular, we show that the $\mathbf{w}_j$ indeed decay with higher algebraic rate than the $\theta_j$.

Our third main contribution is the computation of nonlinear terms in the new coordinate systems.
Leading nonlinear terms turn out to involve \emph{discrete derivatives}, associated with the discrete translational
symmetry near the periodic pattern. Similarly to the scalar case, these discrete derivatives render the nonlinearity irrelevant.
From a different view point, dependence on derivatives, only, indicates the presence of a conservation law:
An equation $u_t=u_{xx}+f(u_x)$ can be rewritten as $v_t=v_{xx}+(f(v))_x$, for $v=u_x$, and the gain in decay is now clear
from an integration by parts in the variation of constant formula. An analogous observation applies to the
$\theta-\mathbf{W}$ system, where discrete derivatives in the nonlinearity reflect a discrete conservation law.

Together, these observations quite readily imply a nonlinear stability result--Theorem \ref{t:1} as shown above.

%
The remainder of this paper is organized as follows.
In Section \ref{s:3}, we construct the normal form.
Section \ref{s:4} contains linear estimates in Fourier-Bloch space.
Section \ref{s:5} converts those decay estimates into $L^p-L^q$ decay estimates in physical space.
Section \ref{s:6} contains the proof of the nonlinear stability result.
We relegate a detailed description of the nonlinearity, and the spectral properties and the analytic semigroup results of the linear
operator to the appendix.

\paragraph{Notation}Throughout we will use the following notation.
\begin{itemize}
 \item $(\cdot,\cdot)$ is the standard inner product on $\R^n$ given by
       $$
        (\mathbf{u},\mathbf{v})=\sum_{j=1}^nu_jv_j, \text{ for any }\mathbf{u}=\{u_j\}_{j=1}^n,\mathbf{v}=\{v_j\}_{j=1}^n\in\R^n.
        $$
 \item $\langle\cdot,\cdot\rangle$ is the standard inner product on the Hilbert space $(L^2(-\pi,\pi))^n$ given by
        $$
        \langle\mathbf{u},\mathbf{v}\rangle=\int_{-\pi}^{\pi}(\mathbf{u}(x),\mathbf{v}(x))\rmd x,
        \text{ for any } \mathbf{u},\mathbf{v}\in(L^2(-\pi,\pi))^n.
        $$
 \item $\llangle \cdot,\cdot\rrangle$ is the standard inner products on $(\ell^2)^n$,
       or the $(\ell^p)^n$--$(\ell^q)^n$ pairing, given by
       $$
        \llangle\underline{\mathbf{u}},\underline{\mathbf{v}}\rrangle=
        \sum_{j\in\Z}(\mathbf{u}_j,\mathbf{v}_j),
        \text{ for any }\mathbf{u}=\{u_j\}_{j\in\Z},\mathbf{v}=\{v_j\}_{j\in\Z}.
        $$
\end{itemize}
We denote the Euclidean norm in Euclidean spaces as $|\cdot|$, the norm in a general Banach space $\mathscr{X}$
as $\|\cdot\|_{\mathscr{X}}$, and the norm of a linear operator from a Banach space $\mathscr{X}$ to $\mathscr{Y}$ as
$\opnorm{\cdot}_{\mathscr{X}\rightarrow\mathscr{Y}}$. For the case $\mathscr{Y}=\mathscr{X}$, the last norm notation simply becomes
$\opnorm{\cdot}_{\mathscr{X}}$.

\section{Normal form}\label{s:3}
As we pointed out in the introduction, the linear system for the perturbation
\beq
\label{e:perlin}
\mathbf{v}_t=\partial_{xx}\mathbf{v} +\mathbf{f}'(\mathbf{u}_\star(x))\mathbf{v},
\eeq
is expected to exhibit diffusive decay for the linear part. This weak decay is not obviously strong enough to conclude nonlinear
decay because of quadratic and cubic terms in the nonlinearity. Our approach here converts the system \eqref{e:perlin} into
an infinite-dimensional lattice dynamical system for $\underline{\mathbf{V}}=\{\mathbf{V}_j\}_{j\in\Z}$, where $\mathbf{V}_j=
(\theta_j,\mathbf{W}_j)\in\R\times L^p(-\pi,\pi)$ for all $j\in\Z$. Here, the scalar component $\theta_j$ of $\mathbf{V}_j$ measures
\emph{local shifts} of the primary periodic pattern, and the infinite-dimensional component, $\mathbf{W}_j$, represents
\emph{local complements}. In such a representation, one expects diffusive decay of $\theta_j$ and faster decay of $\mathbf{W}_j$.
We will make this precise in Section \ref{s:5}. In fact, the linear asymptotics of $\theta_j$ are equivalent to the discrete diffusion
\beqs
\dot{\theta}_j=d(\theta_{j+1}-2\theta_j+\theta_{j-1}).
\eeqs
The key idea is that in this lattice system, nonlinear terms in the $\theta$-equations involve discrete derivatives, $\theta_{j+1}-\theta_j$,
rather than $\theta_j$ alone. Roughly speaking, we expect $\theta$-dependence to disappear when $\theta_j=\theta_{j+1}$
for all $j\in\Z$ due to shift invariance of the original system.
Just like in the continuous scalar heat equation, these derivatives decay faster, so that terms like
$(\theta_{j+1}-\theta_j)^2$ are now irrelevant, that is, they do not alter linear diffusive decay.

In summary, we will find a system, where the linear part exhibits diffusive decay, and where nonlinearities are \emph{explicitly}
irrelevant. In this sense, our transformation has eliminated lower-order terms in the system, that turn out not to contribute to leading
order dynamics. The term normal form alludes to this elimination of lower-order terms by comparing with normal form theory in ODEs, where
coordinate changes are used to simplify equations and systems at least locally, mostly through removing lower-order terms in the Taylor jet
of the equation or system.

The remainder of this section is organized as follows. We discuss local well-posedness and ``chopping-up'', the first key step in the transformation
to a lattice system in Section \ref{ss:31}. The ultimate transformation towards a quasilinear lattice dynamical system is constructed in Section
\ref{ss:32}. Key steps involve separation of the neutral phase $\theta_j$ and a smoothing procedure at the chopping boundaries.

\subsection{Well-posedness: spatially extended system and lattice system}\label{ss:31}
%
We first show local-in-time well-posedness of the system \eqref{e:23} on the space $X=L^1\cap L^\infty$.
\bl
\label{l:31}
The initial value problem of the semi-linear parabolic system \eqref{e:23}
is locally well-posed in $X$. To be precise, the following assertions hold:
\begin{enumerate}
 \item (\textbf{existence and uniqueness}) For any given $\mathbf{v}^0\in X$, there exists some $T>0$,
 depending only on $\|\mathbf{v}^0\|_X$, such that the system \eqref{e:23}
 admits a unique mild solution
 $$\mathbf{v}\in C^0([0,T],(L^1(\R))^n)\cap C^0((0,T],(L^\infty(\R))^n).$$
 Here a mild solution solves the integral-equation variant of \eqref{e:23}.
 \item (\textbf{regularity}) The solution $\mathbf{v}(t,x)$ to \eqref{e:23} is smooth for $t\in(0,T]$. Moreover,
there exists $C>0$ such that, for all $t\in (0,T]$,
\beqs
\|\mathbf{v}(t)\|_{H^2}\leq C t^{-1}\|\mathbf{v}^0\|_{X}.
\eeqs
\end{enumerate}
\el
\bpf
The existence and uniqueness follow directly from \cite{henry} and \cite{lunardi}.
To show that $\|\mathbf{v}(t)\|_{H^2}\leq C t^{-1}\|\mathbf{v}^0\|_{X}$,
we first note that for any $T_0\in(0, T)$, there exists $C(T_0)>0$ such that
\beq
\label{e:part1}
\|\mathbf{v}(t)\|_{H^2}\leq C(T_0)\|\mathbf{v}^0\|_{L^2},
\text{ for all }t\in(T_0, T).
\eeq
Moreover, by \cite[Thm. 7.1.5]{lunardi}, for every $\mathbf{v}^0\in L^2$, there are $T_1>0$ and $C(T_1)>0$ such that
\beqs
\|\mathbf{v}(t/2)\|_{H^1}\leq C(T_1)(t/2)^{-1/2}\|\mathbf{v}^0\|_{L^2},\quad
\|\mathbf{v}(t)\|_{H^2}\leq C(T_1)(t/2)^{-1/2}\|\mathbf{v}(t/2)\|_{H^1}, \text{ for all }t\in(0, T_1),
\eeqs
which implies that
\beq
\label{e:part2}
\|\mathbf{v}(t)\|_{H^2}\leq \frac{C(T_1)}{2}t^{-1}\|\mathbf{v}^0\|_{L^2},\text{ for all }t\in(0, T_1).
\eeq
Combining \eqref{e:part1} and \eqref{e:part2}, we conclude our proof.
\epf
\br
By Lemma \ref{l:31}, we can assume without loss of generality that, in the proof of Theorem \ref{t:1},
the initial perturbation is small in $X\cap H^2$.
\er

Now suppose that $\mathbf{v}(t,x)$ is a solution to (\ref{e:23}), close to $0$. In particular,
$\mathbf{v}(t,x)$ is close to $0$ on all intervals $[2\pi(j-1/2), 2\pi(j+1/2)]$, $j\in \mathbb{Z}$.
Then instead of solving (\ref{e:23}), we claim that it is equivalent to solve the infinite-dimensional system, for all $j\in \mathbb{Z}$,
\beq
\label{e:313}
\begin{cases}
\partial_{t}\mathbf{v}_j=D\partial_{xx}\mathbf{v}_j+\mathbf{f}^\prime(\mathbf{u}_\star)\mathbf{v}_j
+\mathbf{g}(x,\mathbf{v}_j)\\
\mathbf{v}_j(t,\pi)=\mathbf{v}_{j+1}(t,-\pi)\\
\partial_{x}\mathbf{v}_j(t,\pi)=\partial_{x}\mathbf{v}_{j+1}(t, -\pi).
\end{cases}
\eeq
In order to justify the well-posedness of (\ref{e:313}), we first introduce the chopped space
\beq
\label{e:Xchop}
X_{\mathrm{ch}}=\ell^1(\Z,(L^1(-\pi,\pi))^n)\cap \ell^\infty(\Z,(L^\infty(-\pi,\pi))^n),
\eeq
with the norm defined as
\beqs
\|\underline{\mathbf{w}}\|_{X_{\mathrm{ch}}}=\sum_{j\in \mathbb{Z}}\|\mathbf{w}_j\|_{L^1}+\sup_{j\in \mathbb{Z}}\|\mathbf{w}_j\|_{L^\infty},
\text{ for any } \underline{\mathbf{w}}=\{\mathbf{w}_j\}_{j\in \mathbb{Z}}\in X_{\mathrm{ch}}.
\eeqs
We then consider the chopping map
\beq
\label{e:Tchop}
\begin{matrix}
\mathscr{T}_{\mathrm{ch}}: &    X_{\mathrm{ch}}    &\longrightarrow & X \\
\quad  & \underline{\mathbf{v}}&\longmapsto  &\mathscr{T}_{\mathrm{ch}}(\underline{\mathbf{v}}),
\end{matrix}
\eeq
where $X$ is defined in \eqref{e:XX1} and
$\mathscr{T}_{\mathrm{ch}}(\underline{\mathbf{v}})(2\pi j+x)=\mathbf{v}_j(x), \text{for all }x\in[-\pi,\pi] \text{ and }j\in\mathbb{Z}$.
It is not hard to see that $\mathscr{T}_{\mathrm{ch}}$ is an isomorphism and thus we have the diagram
\beqs
\begin{matrix}
X^1   & \overset{A}{\longrightarrow}         &X        \\
\mathscr{T}_{\mathrm{ch}}\uparrow    &                                        &\mathscr{T}_{\mathrm{ch}}\uparrow\\
X_{\mathrm{ch}}^1       & \overset{A_{\mathrm{ch}}}{\longrightarrow}               &X_{\mathrm{ch}},
\end{matrix}
\eeqs
where $X_{\mathrm{ch}}^1:=\mathscr{T}_{\mathrm{ch}}^{-1}(X^1) $ and
\beq
\label{e:A1}
\begin{matrix}
 A_{\mathrm{ch}}: & X_{\mathrm{ch}}^1 & \longrightarrow & X_{\mathrm{ch}}\\
 &\underline{\mathbf{v}}&\longmapsto&\mathscr{T}_{\mathrm{ch}}^{-1}A\mathscr{T}_{\mathrm{ch}}\underline{\mathbf{v}}.
\end{matrix}
\eeq
More specifically, $(A_{\mathrm{ch}}\underline{\mathbf{v}})_j=D\partial_{xx}\mathbf{v}_j+\mathbf{f}^\prime(\mathbf{u}_\star)\mathbf{v}_j$.
To describe $X_{\mathrm{ch}}^1$,
we define
\begin{eqnarray*}
 &\widetilde{\mathscr{D}}(A_{\mathrm{ch}},X_{\mathrm{ch}}):=\ell^1(W^{2,1}(-\pi,\pi))\cap \ell^\infty(W^{2,\infty}(-\pi,\pi)),\\
 &\mathscr{D}(A_{\mathrm{ch}},X_{\mathrm{ch}}):=\{\underline{\mathbf{v}}\in \widetilde{\mathscr{D}}(A_{\mathrm{ch}},X_{\mathrm{ch}})\mid
\mathbf{v}_j^{(k)}(t,\pi)=\mathbf{v}_{j+1}^{(k)}(t,-\pi), t\geq 0, j\in\Z, k=0,1
\}.
\end{eqnarray*}
\bl
We have $X_{\mathrm{ch}}^1=\mathscr{D}(A_{\mathrm{ch}},X_{\mathrm{ch}})$.
\el
\bpf
From the definition, we find $X_{\mathrm{ch}}^1\subseteq\mathscr{D}(A_{\mathrm{ch}},X_{\mathrm{ch}})$.
We only need to show that for any given $\underline{\mathbf{v}}\in\mathscr{D}(A_{\mathrm{ch}},X_{\mathrm{ch}})$,
we have $\mathbf{v}=\mathscr{T}_{\mathrm{ch}}(\underline{\mathbf{v}})\in X^1$. In fact, for arbitrary $\mathbf{w}\in C_c^\infty$, we obtain
\begin{equation*}
  \langle\mathbf{v},\mathbf{w}^\prime\rangle_{L^2(\R)}=\sum_{j\in\Z}\langle\mathbf{v}_j(x),\mathbf{w}^\prime(2\pi j+x)\rangle
  =-\sum_{j\in\Z}\langle\mathbf{v}^\prime_j(x),\mathbf{w}(2\pi j+x)\rangle
 =-\langle\mathscr{T}_{\mathrm{ch}}(\{\mathbf{v}_j^\prime\}_{j\in\Z}),\mathbf{w}\rangle_{L^2(\R)},
\end{equation*}
which shows that $\mathbf{v}^\prime=\mathscr{T}_{\mathrm{ch}}(\{\mathbf{v}_j^\prime\}_{j\in\Z})\in X$. Similarly, we have
$\mathbf{v}^{\prime\prime}\overset{a.e.}{=}\mathscr{T}_{\mathrm{ch}}(\{\mathbf{v}_j^{\prime\prime}\}_{j\in\Z})\in X$.
\epf

In all, we conclude that our initial value problem for a spatially extended system (\ref{e:23}) is equivalent to an
initial value problem for a lattice system as follows.
\begin{equation}
\label{e:313m}
\begin{cases}
\partial_{t}\underline{\mathbf{v}}=A_{\mathrm{ch}}\underline{\mathbf{v}}+\mathbf{G}(\underline{\mathbf{v}}),\quad x\in(-\pi,\pi), t>0,\\
\mathbf{v}_j^{(k)}(t,\pi)=\mathbf{v}_{j+1}^{(k)}(t,-\pi), \quad k=0,1, j\in\Z, t\geq 0,\\
\mathbf{v}_j(0,x)=\mathbf{v}^0(2\pi j+x),\quad x\in[-\pi,\pi],j\in\Z,
\end{cases}
\end{equation}
where $\mathbf{G}(\underline{\mathbf{v}})=\{\mathbf{g}(x,\mathbf{v}_j)\}_{j\in\Z}$.

\br For any solution $\underline{\mathbf{v}}$ to (\ref{e:313m}), we have all higher matching boundary conditions, that is,
$\partial_{x}^m\mathbf{v}_j(t,\pi)=\partial_{x}^m\mathbf{v}_{j+1}(t,-\pi)$, for all $t>0$, $j\in\mathbb{Z}$
and $m\in\mathbb{Z}^+$.
\er

\subsection{Lattice system: phase decomposition and boundary-condition matching}\label{ss:32}
We start with sketching the construction of the normal form step by step without rigorous justification.
We first decompose each $2\pi$-long piece $\mathbf{v}_j(x)=\mathbf{v}(2\pi j+x)$
into a linearly neutral phase and a stable phase and then match the boundary conditions for the stable phase. This two-step
smooth phase decomposition procedure will be summarized and justified rigorously in a lemma at the end of this section.

We now decompose each $\mathbf{v}_j$ according to
\beqs
\bcs
\mathbf{v}_j(x)=\mathbf{w}_j(x)+\mathbf{u}_\star (x-\theta_j)-\mathbf{u}_\star (x) \\
\langle\mathbf{w}_j(x), \mathbf{u}_{\mathrm{ad}}(x-\theta_j)\rangle=0,\\
\ecs
\eeqs
where $\mathbf{u}_{\mathrm{ad}}$
is an element in the kernel of the adjoint operator of $B(0)$ with $\langle\mathbf{u}_{\star}^\prime, \mathbf{u}_{\mathrm{ad}}\rangle=1$.
Substituting this expression into \eqref{e:23},
we can therefore formally derive a system for $\theta_j$ and $\mathbf{w}_j$, which takes the explicit form
\beq
\label{e:310}
\bcs
\bld
\dot{\theta}_j=&\frac{1}{-1+\langle\mathbf{w}_j(x),\mathbf{u}_{\mathrm{ad}}^\prime(x-\theta_j)\rangle}
[-(\mathbf{w}_j(\pi)-\mathbf{w}_j(-\pi), D\mathbf{u}
_{\mathrm{ad}}^\prime(\pi-\theta_j))\\
&+(\partial_x\mathbf{w}_j(\pi)-\partial_x\mathbf{w}_j(-\pi), D\mathbf{u}_{\mathrm{ad}}(\pi-
\theta_j))+\langle\tilde{\mathbf{g}}(\theta_j,\mathbf{w}_j),\mathbf{u}_{\mathrm{ad}}(x-\theta_j)\rangle]\\
\dot{\mathbf{w}}_j =&D\partial_{xx}\mathbf{w}_j+\mathbf{u}_{\star, \theta}(x-\theta_j)\dot{\theta}_j+
\mathbf{f}(\mathbf{w}_j+\mathbf{u}_\star(x-\theta_j))-\mathbf{f}(\mathbf{u}_\star(x-\theta_j)),
\eld
\ecs
\eeq
with boundary conditions
\beqs
\bcs
\bld
&\partial_x^{m}\mathbf{w}_j(\pi)-\partial_x^{m}\mathbf{w}_{j+1}(-\pi)=
\mathbf{u}_\star^{(m)}(\pi-\theta_{j+1})-\mathbf{u}_\star^{(m)}(\pi-\theta_j),
\text{ for } m=0,1\\
&\langle\mathbf{w}_j(x),\mathbf{u}_{\mathrm{ad}}(x-\theta_j)\rangle=0,\\
\eld
\ecs
\eeqs
where
\beqs
\tilde{\mathbf{g}}(\theta_j,\mathbf{w}_j)=\mathbf{f}(\mathbf{w}_j+\mathbf{u}_\star(x-\theta_j))-
\mathbf{f}(\mathbf{u}_\star(x-\theta_j))-\mathbf{f}^\prime(\mathbf{u}_\star(x-\theta_j))\mathbf{w}_j.
\eeqs
\br
In the second equation of (\ref{e:310}), $\dot{\theta}_j$ represents the right hand side of the first equation.
\er

Note that $\mathbf{w}_j$ is in a codimension-one subspace depending on $\theta_j$. More formally, we mapped every $\mathbf{v}_j$
into a vector bundle. Also, the boundary conditions are now nonlinear.
These facts generate technical difficulties so that we find it easier to work with a further modified system,
where, for all $j\in\mathbb{Z}$, we substitute
\beq
\label{e:PD}
\mathbf{w}_j(x)=\mathbf{W}_j(x)+\mathbf{H}(x,\theta_{j-1},\theta_j,\theta_{j+1},\mathbf{W}_j).
\eeq
For simplicity, we denote $\mathbf{H}_j(x)=\mathbf{H}(x,\theta_{j-1},\theta_j,\theta_{j+1},\mathbf{W}_j)$.
In the new coordinates $\underline{\mathbf{V}}=(\underline{\theta},\underline{\mathbf{W}})$, where
$\underline{\theta}=\{\theta_j\}_{j\in \mathbb{Z}}$ and $\underline{\mathbf{W}}=\{\mathbf{W}_j\}_{j\in \mathbb{Z}}$,
we will have again ``homogeneous matching boundary conditions" and all $\mathbf{W}_j$'s are in a fixed codimension-1 subspace,
that is, for all $j\in\mathbb{Z}$,
\begin{eqnarray}
&\partial_x^{m}\mathbf{W}_j(\pi)-\partial_x^{m}\mathbf{W}_{j+1}(-\pi)=0,\text{ for } m=0,1,	 \label{e:33a}\\
&\langle \mathbf{W}_j(x), \mathbf{u}_{\mathrm{ad}}(x)\rangle=0.\label{e:33b}
\end{eqnarray}

We now construct $\underline{\mathbf{H}}=\{\mathbf{H}_j(x)\}_{{j\in \mathbb{Z}}}$ explicitly in the form
\beq
\label{e:HD}
\mathbf{H}_j=\mathbf{H}_j^1+\mathbf{H}_j^2,
\eeq
where $\mathbf{H}_j^1$ accomplishes ``homogeneous matching boundary conditions" \eqref{e:33a} and $\mathbf{H}_j^2$ corrects
so that every $\mathbf{W}_j$ is perpendicular to $\mathbf{u}_{\mathrm{ad}}$ \eqref{e:33b}. First, we construct $\mathbf{H}_j^1$.
In order to accomplish \eqref{e:33a}, one readily verifies that we need
\beqs
\partial_x^m\mathbf{H}_j^1(\pi)-\partial_x^m\mathbf{H}_{j+1}^1(-\pi)=
\mathbf{u}_\star^{(m)}(\pi-\theta_{j+1})-\mathbf{u}_\star^{(m)}(\pi-\theta_j),\text{ for } m=0,1,
\eeqs
which can be achieved by choosing
\beqs
\bcs
\bld
&\mathbf{H}_j^1(x)=\frac{1}{2}(\mathbf{u}_\star(x-\theta_{j+1})-\mathbf{u}_\star(x-\theta_j)), \text{for } x\sim\pi,\\
&\mathbf{H}_j^1(x)=\frac{1}{2}(\mathbf{u}_\star(x-\theta_{j-1})-\mathbf{u}_\star(x-\theta_j)), \text{for } x\sim-\pi.
\eld
\ecs
\eeqs
In light of this observation, we let
\beq
\label{e:H1}
\mathbf{H}_j^1(x)=
\frac{1}{2}\phi(x)(\mathbf{u}_\star(x-\theta_{j+1})-\mathbf{u}_\star(x-\theta_{j-1}))+
\frac{1}{4}(\mathbf{u}_\star(x-\theta_{j+1})+\mathbf{u}_\star(x-\theta_{j-1})-2\mathbf{u}_\star(x-\theta_j)),
\eeq
where $\phi$ is a smooth odd, increasing function on $[-\pi,\pi]$ such that
\beqs
\phi(x)=
\bcs
\bld
&\frac{1}{2}, \quad \quad\text{for }x>\frac{\pi}{2},\\
&-\frac{1}{2},\quad \text{for }x<-\frac{\pi}{2}.
\eld
\ecs
\eeqs
To be specific, we can choose
\beqs
\phi(x)=[\eta\ast \chi_{[0,\infty)}](x)\cdot\chi_{[-\pi,\pi]}(x)-\frac{1}{2},
\eeqs
where $\chi_{J}$ is the characteristic function of the interval $J$ and $\eta$ is a smooth nonnegative even mollifier such that
\[
\int_{\mathbb{R}}\eta(x)\rmd x=1, \text{and }|\eta(x)|=0, \text{ for all } |x|>\frac{\pi}{2}.
\]
In order to keep $\mathbf{H}_j$ identical with $\mathbf{H}_j^1$ near $\pm\pi$, $\mathbf{H}_j^2$ has to be $0$ near $\pm\pi$.
We first note that there exists an odd function $\psi\in (C_c^\infty(-\pi,\pi))^n$ such that $\langle\psi,\mathbf{u}_{\mathrm{ad}}\rangle=1$
since $(C_c^\infty(-\pi,\pi))^n$ is dense in $(L^2(-\pi,\pi))^n$ and $\langle\mathbf{u}_\star^\prime,\mathbf{u}_{\mathrm{ad}}\rangle=1$.
We then define
\beq
\label{e:H2}
\mathbf{H}_j^2=c_j\psi(x-\theta_j),
\eeq
where
\beq
\label{e:cj}
c_j=-\langle\mathbf{H}_j^1,\mathbf{u}_{\mathrm{ad}}(x-\theta_j)\rangle-
\langle\mathbf{W}_j,\mathbf{u}_{\mathrm{ad}}(x-\theta_j)-\mathbf{u}_{\mathrm{ad}}(x)\rangle.
\eeq
Noting  that $\theta_j$ and $\mathbf{W}_j$ are small,
this concludes the construction of $\mathbf{H}_j$.

Defining $X_{\mathrm{ch}}^\perp=\{\underline{\mathbf{v}}\in X_{\mathrm{ch}}\mid \langle\mathbf{v}_j,
\mathbf{u}_{\mathrm{ad}}\rangle=0, \text{ for all } j\in \mathbb{Z}\}$,
where $X_{\mathrm{ch}}$ is defined in \eqref{e:Xchop}, we summarize the
``smooth phase decomposition'' procedure, denoted as $\mathscr{T}_{\mathrm{phd}}$, in the following lemma.
\bl
\label{l:33}
The ``smooth phase decomposition'' operator $\mathscr{T}_{\mathrm{phd}}$, as constructed above, is a smooth local diffeomorphism.
More precisely, there are two neighborhoods of zero $\mathscr{U}\in X_{\mathrm{ch}}$, $\mathscr{V}\in \ell^1\times X_{\mathrm{ch}}^{\perp}$
such that the nonlinear transformation
\beqs
\begin{matrix}
\mathscr{T}_{\mathrm{phd}}:&\mathscr{V}&\longrightarrow&\mathscr{U}\\
&\underline{\mathbf{V}}=(\underline{\theta}, \underline{\mathbf{W}}(x))&\longmapsto&\{\mathbf{W}_j(x)+\mathbf{H}_j(x)+
\mathbf{u}_\star(x-\theta_j)-\mathbf{u}_\star(x)\}_{j\in \mathbb{Z}}
\end{matrix}
\eeqs
is invertible with $\mathscr{T}_{\mathrm{phd}}$ and $\mathscr{T}_{\mathrm{phd}}^{-1}$ smooth.
Its derivative at the origin is
\beq
\label{e:T2}
\begin{matrix}
\mathscr{L}_{\mathrm{phd}}:=\mathscr{T}_{\mathrm{phd}}^\prime(0): &\ell^1\times X_{\mathrm{ch}}^{\perp}&\longrightarrow&X_{\mathrm{ch}}\\
&\underline{\mathbf{V}}=(\underline{\theta}, \underline{\mathbf{W}})&\longmapsto&\underline{\mathbf{W}}+\underline{\mathbf{E}}*\underline{\theta},
\end{matrix}
\eeq
where $\underline{\mathbf{E}}$ is defined in \eqref{e:E} below.
\el
\bpf
We claim that
\begin{itemize}
\item[(\rmnum{1})]$\mathscr{T}_{\mathrm{phd}}(0)=0$;
\item[(\rmnum{2})]$\mathscr{T}_{\mathrm{phd}}$ is $C^\infty$;
\item[(\rmnum{3})]$\mathscr{T}_{\mathrm{phd}}^\prime(0)$, denoted as $\mathscr{L}_{\mathrm{phd}}$,
 is an invertible bounded linear operator.
\end{itemize}
Property (\rmnum{1}) is straightforward. As for (\rmnum{2}), $\mathscr{T}_{\mathrm{phd}}$ is smooth with respect to
$\underline{\mathbf{W}}$ due to the fact that $\mathscr{T}_{\mathrm{phd}}$ is linear in $\underline{\mathbf{W}}$ for fixed $\underline{\theta}$.
On the other hand, the smoothness of $\mathscr{T}_{\mathrm{phd}}$ with respect to $\underline{\theta}$
can be readily reduced to the smoothness of the mapping
$\underline{\theta}\mapsto\{\mathbf{u}_*(x-\theta_j)-\mathbf{u}_*(x)\}_{j\in\Z}$. A direct calculation shows that,
for given $m\in\Z^+$,
the $m$th-derivative mapping at $\underline{\theta}$ is
$\underline{\eta}\mapsto\{\frac{1}{m!}\mathbf{u}_*^{(m)}(\theta_j-x)\eta_j\}_{j\in\Z}$.
We now only have to show that (\rmnum{3}) is true.
In fact, the linear part of $\{\mathbf{H}_j+\mathbf{u}_\star(x-\theta_{j})-\mathbf{u}_\star(x)\}_{j\in\Z}$ with respect to
$(\underline{\theta}, \underline{\mathbf{W}})$ around $(0,0)$ is
$\underline{\mathbf{E}}*\underline{\theta}=\{\sum_{k\in\Z}\mathbf{E}_{j-k}\theta_k\}_{j\in\Z}$,
where $\underline{\mathbf{E}}=\{\mathbf{E}_j\}_{j\in\Z}$ with
\beq
\label{e:E}
\mathbf{E}_j=\begin{cases}
              \frac{1}{4}\psi(x)-(\frac{1}{4}+\frac{1}{2}\phi(x))\mathbf{u}_{\star}^\prime(x), & j=-1,\\
              -\frac{1}{2}(\psi(x)+\mathbf{u}_{\star}^\prime(x)), & j=0,\\
              \frac{1}{4}\psi(x)-(\frac{1}{4}-\frac{1}{2}\phi(x))\mathbf{u}_{\star}^\prime(x), & j=1,\\
              0,& others.
             \end{cases}
\eeq
Then we have the linear phase decomposition operator
\beqs
\begin{matrix}
\mathscr{L}_{\mathrm{phd}}: &\ell^1\times X_{\mathrm{ch}}^{\perp}&\longrightarrow&X_{\mathrm{ch}}\\
&\underline{\mathbf{V}}=(\underline{\theta}, \underline{\mathbf{W}})&\longmapsto&\underline{\mathbf{W}}+\underline{\mathbf{E}}*\underline{\theta}.
\end{matrix}
\eeqs
Moreover, through direct calculation, it is not hard to obtain the bounded inverse of $\mathscr{L}_{\mathrm{phd}}$
\beqs
\begin{matrix}
\mathscr{L}_{\mathrm{phd}}^{-1}: &X_{\mathrm{ch}}&\longrightarrow&\ell^1\times X_{\mathrm{ch}}^{\perp}\\
&\underline{\mathbf{v}}&\longmapsto&(F\underline{\mathbf{v}}, \underline{\mathbf{v}}-\mathbf{E}* F\underline{\mathbf{v}}),
\end{matrix}
\eeqs
where
\beq
\label{e:F}
\begin{matrix}
 F:&X_{\mathrm{ch}}&\longrightarrow&\ell^1\\
 &\underline{\mathbf{v}}=\{\mathbf{v}_j\}_{j\in\mathbb{Z}}&\longmapsto&\{-\langle\mathbf{v}_j,
 \mathbf{u}_{\mathrm{ad}}\rangle\}_{j\in\mathbb{Z}}.
\end{matrix}
\eeq
By (\rmnum{1}), (\rmnum{2}) and the inverse function theorem, the conclusion of the lemma follows.
\epf
\br
The above lemma still holds when replacing $X_{\mathrm{ch}}$ with $\mathscr{T}_{\mathrm{ch}}(X\cap H^2)$ and the proof is similar.
\er

In the new coordinates, the system contains lengthy expressions. We therefore introduce some simplifying notation first.
\begin{equation}
\label{e:deltaGamma}
\begin{matrix}
&\delta_+: &  \mathbb{C}^{\mathbb{Z}} &\longrightarrow & \mathbb{C}^{\mathbb{Z}}\\
&       & \underline{x}=\{x_j\}_{j\in\mathbb{Z}}&\longmapsto & \{x_{j+1}-x_j\}_{j\in\mathbb{Z}}.\\
 &&&&\\
 &\delta_-: & \mathbb{C}^{\mathbb{Z}} &\longrightarrow & \mathbb{C}^{\mathbb{Z}}\\
 &       & \underline{x}&\longmapsto & \{x_{j}-x_{j-1}\}_{j\in\mathbb{Z}}.\\
&&&&\\
&\Gamma:&(C([-\pi,\pi],\R^n))^{\mathbb{Z}}&\longrightarrow & \mathbb{R}^{\mathbb{Z}}\\
&       &\underline{\mathbf{v}}&\longmapsto &\{(\mathbf{v}_j(-\pi), D\mathbf{u}_{\mathrm{ad}}^\prime(\pi))\}_{j\in\mathbb{Z}}.\\
\end{matrix}
\end{equation}
Now, sorting out the linear terms, our lattice system is
\beq
\label{e:34}
\begin{pmatrix}
\dot{\underline{\theta}}\\ \dot{\underline{\mathbf{W}}}
\end{pmatrix}
=A_{\mathrm{nf}}\begin{pmatrix}
\underline{\theta}\\ \underline{\mathbf{W}}
\end{pmatrix}+\begin{pmatrix}
\mathbf{N}^\theta(\underline{\theta},\underline{\mathbf{W}})\\ \mathbf{N}^\mathbf{w}(\underline{\theta},\underline{\mathbf{W}})
\end{pmatrix}
\eeq
with boundary-matching and phase-decomposition conditions (\ref{e:33a}), (\ref{e:33b}),
where $\mathbf{N}^{\theta/\mathbf{w}}$ represent the nonlinear terms of the system and
\beq
\label{e:A2}
A_{\mathrm{nf}}=\mathscr{L}_{\mathrm{phd}}^{-1} A_{\mathrm{ch}} \mathscr{L}_{\mathrm{phd}}
=\begin{pmatrix}
F\\ \id-\mathbf{E}*F
\end{pmatrix}
A_{\mathrm{ch}}
\begin{pmatrix}
\mathbf{E}*&\id
\end{pmatrix}
=\begin{pmatrix}
0&\delta_+\Gamma \\
A_{\mathrm{ch}}\mathbf{E}*&A_{\mathrm{ch}}-\mathbf{E}*\delta_+\Gamma
\end{pmatrix},
\eeq
where $A_{\mathrm{ch}}$ is the linear operator acting on the chopped variables; see \eqref{e:A1}.
\br
\begin{itemize}
\item[(\rmnum{1})]Due to the fact that $\mathscr{T}_{\mathrm{ch}\slash \mathrm{nf}}$ are isomorphisms,
$A_{\mathrm{nf}}=\mathscr{L}_{\mathrm{phd}}^{-1}\mathscr{T}_{\mathrm{ch}}^{-1}
A \mathscr{T}_{\mathrm{ch}}\mathscr{L}_{\mathrm{phd}}$ shares many properties with $A$.
For example, $A_{\mathrm{nf}}$ is sectorial in $\ell^1\times X_{\mathrm{ch}}^{\perp}$ since $A$ is sectorial in $X$.
Here we use the definition of a sectorial operator from \cite{lunardi} which does not require the operator to have a dense domain.
\item[(\rmnum{2})]We relegate the detailed estimates of the nonlinear terms to Lemma \ref{l:61}
   in the appendix since expressions are lengthy. We have, roughly,
    \beqs
    \bcs
    &|\mathbf{N}^\theta|\sim |(\delta_+\underline{\theta})^2|+|\underline{\theta}^3||\delta_+\underline{\theta}|+
    (|\underline{\theta}|+|\underline{\mathbf{W}}|)(|\underline{\mathbf{W}}|+|\delta_+\partial_{xx}\underline{\mathbf{W}}|)
    +|\underline{\mathbf{W}}^2|\\
    &|\mathbf{N}^\mathbf{w}|\sim |\underline{\theta}||\delta_+\underline{\theta}|+
    (|\underline{\theta}|+|\underline{\mathbf{W}}|)(|\underline{\mathbf{W}}|+|\delta_+\partial_{xx}\underline{\mathbf{W}}|)
    +|\underline{\mathbf{W}}^2|.
    \ecs
    \eeqs
\item[(\rmnum{3})]Since the branch of continuous spectrum connected to $\lambda=0$ may intersect the branches of continuous
spectrum in $\re \lambda<0$, it is in general not clear how to globally separate neutral from stable modes even linearly.
Phase decompositions have been achieved globally in the case of weak pulse interaction, that is,
in the regime where $\mathbf{u}_\star(x)$ is close to a homoclinic orbit in the ordinary differential system
$D\mathbf{u}_{xx}+\mathbf{f}(\mathbf{u})=0$ ; see \cite{sslsp_2001} for a linear analysis and \cite{zemi_2009} for a nonlinear reduction.
\end{itemize}
\er

\section{Linear Fourier-Bloch estimates  }
\label{s:4}
In Section \ref{s:4} and Section \ref{s:5}, we derive linear diffusive decay in our linear normal form
\beqs
\dot{\underline{\mathbf{V}}}=A_{\mathrm{nf}}\underline{\mathbf{V}}.
\eeqs
To illustrate the idea, we again use the linear heat equation
\beqs
u_t(t,x)=\triangle u(t,x).
\eeqs
In order to obtain the diffusive decay on $\rme^{\triangle t}$, we apply the Fourier transform and obtain
the ``diagonalized'' equation
\beqs
\widehat{u}_t(t,k)=-k^2\widehat{u}(t,k).
\eeqs
Then we have that $|\widehat{u}(t,k)|=\rme^{-k^2t}|\widehat{u}(0,k)|$, for all $t>0$ and $k\in\R$, which, combined with
Young's inequality, will give us diffusive decay for the scalar heat equation.

In light of this procedure, we exploit Fourier transforms and the
Bloch wave decomposition of $A$ to construct an isomorphism diagram, from which we
obtain a direct integral representation of $A_{\mathrm{nf}}$, that is, $\widehat{A}_{\mathrm{nf}}=\int_{-1/2}^{1/2}
\widehat{A}_{\mathrm{nf}}(\sigma)\rmd \sigma$.
Unlike the explicit expression of $\rme^{-k^2t}$, the estimates on $\rme^{\widehat{A}_{\mathrm{nf}}(\sigma)t}$
are more intricate and their derivation will occupy most of this section.

To show the conjugacy between the linear normal form and its counterpart in a Fourier-Bloch space, we
build a commutative isomorphism diagram involving the underlying spaces for these two operators,
the linear operator $A$ and its Bloch wave decomposition.
To this end, we recall the definitions of the linearized operator $A$ in \eqref{e:A},
the chopping operator $\mathscr{T}_{\mathrm{ch}}$ in \eqref{e:Tchop},
and the linear phase decomposition operator $\mathscr{L}_{\mathrm{phd}}$ in \eqref{e:T2}
from above. We now consider these operators on $L^2\slash\ell^2$-based spaces, that is, with new notation,
\beq
\label{e:opL2}
\begin{matrix}
\widetilde{A}:&(H^2(\R))^n&\longrightarrow&(L^2(\R))^n,\\
&&&\\
\widetilde{\mathscr{T}}_{\mathrm{ch}}:&\ell^2(\Z,(L^2(\T_{2\pi}))^n)&\longrightarrow&(L^2(\R))^n,\\
&&&\\
\widetilde{\mathscr{L}}_{\mathrm{phd}}:&\ell^2\times\ell^2_\perp(\Z,(L^2(\T_{2\pi}))^n)
&\longrightarrow&\ell^2(\Z,(L^2(\T_{2\pi}))^n),
\end{matrix}
\eeq
where $\mathbb{T}_\alpha=\mathbb{R}/{\alpha \mathbb{Z}}$ is the one-dimensional torus of length $\alpha$ and
\beqs
\ell^2_\perp(\Z,(L^2(\T_{2\pi}))^n)=\{\underline{\mathbf{w}}\in \ell^2(\Z,(L^2(\T_{2\pi}))^n)
\mid \langle\mathbf{w}_j,\mathbf{u}_{\mathrm{ad}}\rangle=0, \text{ for all } j\in \mathbb{Z}\}.
\eeqs
We write $\widehat{\mathbf{u}}=\int_{\mathbb{R}}\mathbf{u}(x)\rme^{-\rmi k x}\rmd x$ and introduce several Fourier transform variants as follows:
\beq
\label{e:FT}
\begin{matrix}
\mathscr{F}:& \ell^2 &\longrightarrow &L^2(\mathbb{T}_1)\\
&\underline{\theta}=\{\theta_j\}_{j\in\Z}&\longmapsto & \sum_{j\in\mathbb{Z}}\theta_j\rme^{-\rmi 2\pi j \sigma},\\
&&&\\
\mathscr{F}_n:& (L^2(\mathbb{T}_{2\pi}))^n&\longrightarrow &(\ell^2)^n\\
&\mathbf{u}(x)&\longmapsto & \underline{\mathbf{u}}=\{\int_{-\pi}^{\pi}\mathbf{u}(x)\rme^{-\rmi\ell x}\rmd x\}_{\ell\in\mathbb{Z}},\\
&&&\\
\mathscr{F}_{\mathrm{ch}}:& \ell^2(\mathbb{Z}, (L^2(\mathbb{T}_{2\pi}))^n)&\longrightarrow &L^2(\mathbb{T}_1, (\ell^2)^n)\\
&\underline{\mathbf{u}}(x)=\{\mathbf{u}_j(x)\}_{j\in\Z}&\longmapsto & \underline{\widehat{\mathbf{u}}}
(\sigma)=\{\sum_{j\in\Z}\int_{\T_{2\pi}}\mathbf{u}_j(x)\rme^{-\rmi(\sigma+
\ell)(2\pi j+x)}\}_{\ell\in\mathbb{Z}},\\
&&&\\
\mathscr{F}_{\mathrm{nf}}:& \ell^2\times\ell^2_\perp(\Z,(L^2(\T_{2\pi}))^n)&\longrightarrow &L^2(\T_1)\times L^2_\perp(\T_1,(\ell^2)^n)\\
&(\underline{\theta},\underline{\mathbf{u}})^T&\longmapsto & (\mathscr{F}(\underline{\theta}),\mathscr{F}_{\mathrm{ch}}(\underline{\mathbf{u}}))^T,
\end{matrix}
\eeq
where
\beqs
L^2_\perp(\T_1,(\ell^2)^n)=\{\underline{\mathbf{w}}\in L^2(\T_1,(\ell^2)^n)
\mid \llangle\underline{\mathbf{w}}(\sigma),\mathscr{F}_n(\rme^{-\rmi \sigma x}\mathbf{u}_{\mathrm{ad}})
\rrangle=0, \text{ for all } \sigma\in \T_1\},
\eeqs
We then have a commutative diagram of isomorphisms as follows,
\beq
\label{diag:31}
\begin{matrix}
(L^2(\R))^n&\overset{\widetilde{\mathscr{T}}_{\mathrm{ch}}}{\longleftarrow}&
\ell^2(\Z,(L^2(\T_{2\pi}))^n)&\overset{\widetilde{\mathscr{L}}_{\mathrm{phd}}}{\longleftarrow}
&\ell^2\times\ell^2_\perp(\Z,(L^2(\T_{2\pi}))^n)\\
\downarrow\mathscr{B}^{-1}&&\downarrow\mathscr{F}_{\mathrm{ch}}&&\downarrow\mathscr{F}_{\mathrm{nf}}\\
L^2(\T_1,(L^2(\T_{2\pi}))^n)&\overset{\widehat{\mathscr{T}}_{\mathrm{ch}}}{\longleftarrow}&L^2(\T_1,(\ell^2)^n)
&\overset{\widehat{\mathscr{L}}_{\mathrm{phd}}}{\longleftarrow}&L^2(\T_1)\times L^2_\perp(\T_1,(\ell^2)^n),
\end{matrix}
\eeq
where $\mathscr{B}^{-1}$ is the inverse of the direct integral defined in \eqref{e:DI}, Section \ref{ss:62} and
\beqs
\begin{matrix}
\widehat{\mathscr{T}}_{\mathrm{ch}}:& L^2(\T_1,(\ell^2)^n)&\longrightarrow &L^2(\T_1,(L^2(\mathbb{T}_{2\pi}))^n)\\
&\underline{\mathbf{u}}(\sigma)=\{\mathbf{u}_j(\sigma)\}_{j\in\Z}&\longmapsto & \mathbf{u}(\sigma)=
(2\pi)^{\frac{1}{2}}\mathscr{F}_n^{-1}\underline{\mathbf{u}}(\sigma),\\
&&&\\
\widehat{\mathscr{L}}_{\mathrm{phd}}:&L^2(\T_1)\times L^2_\perp(\T_1,(\ell^2)^n)&\longrightarrow&L^2(\T_1,(\ell^2)^n)\\
&(\theta(\sigma),\underline{\mathbf{w}}(\sigma))&\longmapsto&
\theta(\sigma)\underline{\widehat{\mathbf{E}}}(\sigma)+\underline{\mathbf{w}}(\sigma).
\end{matrix}
\eeqs
Here we have
\beq
\label{e:Esigma}
\underline{\widehat{\mathbf{E}}}(\sigma)=\mathscr{F}_{\mathrm{ch}}(\underline{\mathbf{E}}),
\eeq
with $\underline{\mathbf{E}}$ defined in \eqref{e:E}.
The inverse of
$\widehat{\mathscr{L}}_{\mathrm{phd}}$, which will be used later, has the expression
\beqs
\begin{matrix}
\widehat{\mathscr{L}}_{\mathrm{phd}}^{-1}:&L^2(\T_1,(\ell^2)^n)&\longrightarrow&L^2(\T_1)\times L^2_\perp(\T_1,(\ell^2)^n)\\
&\underline{\mathbf{w}}(\sigma)&\longmapsto&
(\widehat{F}(\underline{\mathbf{w}}(\sigma)),\underline{\mathbf{w}}(\sigma)-
\underline{\widehat{\mathbf{E}}}(\sigma)\widehat{F}(\underline{\mathbf{w}}(\sigma))),
\end{matrix}
\eeqs
where
\beqs
\begin{matrix}
\widehat{F}:&L^2(\T_1,(\ell^2)^n)&\longrightarrow&L^2(\T_1)\\
&\underline{\mathbf{w}}(\sigma)&\longmapsto&
-\llangle\underline{\mathbf{w}}(\sigma),\mathscr{F}_n(\rme^{-\rmi\sigma x}\mathbf{u}_{\mathrm{ad}})\rrangle.
\end{matrix}
\eeqs

We now use tildes for operators in physical space and hats for their conjugates in Fourier space.
The index ``$\mathrm{ch}$'' refers to the chopped operators, the index ``$\mathrm{phd}$'' refers to
the smooth phase decomposition operators, and  the index ``$\mathrm{nf}$''
refers to the normal form operators.
We then define
\begin{eqnarray}
\label{e:39}
\widetilde{A}_{\mathrm{ch}}:=\widetilde{\mathscr{T}}_{\mathrm{ch}}^{-1}
\widetilde{A}\widetilde{\mathscr{T}}_{\mathrm{ch}},&
\widehat{A}_{\mathrm{ch}}:=\widehat{\mathscr{T}}_{\mathrm{ch}}^{-1}\widehat{A}\widehat{\mathscr{T}}_{\mathrm{ch}},\\
\label{e:315}
\widetilde{A}_{\mathrm{nf}}:=\widetilde{\mathscr{L}}_{\mathrm{phd}}^{-1}\widetilde{\mathscr{T}}_{\mathrm{ch}}^{-1}
\widetilde{A}\widetilde{\mathscr{T}}_{\mathrm{ch}}\widetilde{\mathscr{L}}_{\mathrm{phd}},&
\widehat{A}_{\mathrm{nf}}:=\widehat{\mathscr{L}}_{\mathrm{phd}}^{-1}
\widehat{\mathscr{T}}_{\mathrm{ch}}^{-1}\widehat{A}\widehat{\mathscr{T}}_{\mathrm{ch}}\widehat{\mathscr{L}}_{\mathrm{phd}},
\end{eqnarray}
where, according to the Bloch wave decomposition from Theorem \ref{t:4} in the appendix, we have
\beq
\label{e:316}
 \mathscr{B}^{-1}\widetilde{A}\mathscr{B}=\widehat{A}=\int_{-\frac{1}{2}}^{\frac{1}{2}}B(\sigma)\rmd\sigma,
\eeq
with $B(\sigma)$ defined in \eqref{e:B}.
Therefore, by the commutative diagram \eqref{diag:31} and the equivalence relations in (\ref{e:315}), (\ref{e:316}),
we find the conjugacy
\beq
\label{e:lpf}
\widetilde{A}_{\mathrm{nf}}=\mathscr{F}_{\mathrm{nf}}^{-1}\widehat{A}_{\mathrm{nf}}\mathscr{F}_{\mathrm{nf}}.
\eeq

Just as we pointed out at the beginning of this section, based on this conjugacy, in order to obtain estimates on
$\rme^{\widetilde{A}_{\mathrm{nf}}t}$, we only need to derive estimates on $\rme^{\widehat{A}_{\mathrm{nf}}t}$. To this end,
we first derive an explicit direct integral expression of $\widehat{A}_{\mathrm{nf}}$.
From the equivalence relations in (\ref{e:39}),(\ref{e:316}), it is straightforward to see that
\beq
\label{e:Asigma}
\widehat{A}_{\mathrm{ch}}=\int_{-\frac{1}{2}}^{\frac{1}{2}}\widehat{A}_{\mathrm{ch}}(\sigma)\rmd\sigma,
\text{ with }\widehat{A}_{\mathrm{ch}}(\sigma):=\mathscr{F}_nB(\sigma)\mathscr{F}_n^{-1},
\text{ for all }\sigma\in[-1/2,1/2].
\eeq
Moreover, for any given $(\theta(\sigma),\underline{\mathbf{w}}(\sigma))\in L^2(\T_1)\times L^2_\perp(\T_1,(\ell^2)^n)$ and
fixed $\sigma\in[-\frac{1}{2},\frac{1}{2}]$, by definition, we have,
\beq
\label{e:A2sigma}
\bld
\left(\widehat{A}_{\mathrm{nf}}\begin{pmatrix}\theta\\ \underline{\mathbf{w}}\end{pmatrix}\right)(\sigma)
&=\left(\widehat{\mathscr{L}}_{\mathrm{phd}}^{-1}\widehat{\mathscr{T}}_{\mathrm{ch}}^{-1}
\widehat{A}\widehat{\mathscr{T}}_{\mathrm{ch}}\widehat{\mathscr{L}}_{\mathrm{phd}}
\begin{pmatrix}\theta\\ \underline{\mathbf{w}}\end{pmatrix}\right)(\sigma)\\
&=
\begin{pmatrix} \widehat{F}(\sigma) \\ \id -  \underline{\widehat{\mathbf{E}}}(\sigma)\widehat{F}(\sigma) \end{pmatrix}
\widehat{A}_{\mathrm{ch}}(\sigma)
\begin{pmatrix} \underline{\widehat{\mathbf{E}}}(\sigma) & \id \end{pmatrix}
\begin{pmatrix}
\theta(\sigma)\\
\underline{\mathbf{w}}(\sigma)\\
\end{pmatrix}\\
&=
\begin{pmatrix}
0&R(\sigma)\\
\widehat{A}_{\mathrm{ch}}(\sigma)\underline{\widehat{\mathbf{E}}}(\sigma)&\widehat{A}_
{\mathrm{ch}}(\sigma)-\underline{\widehat{\mathbf{E}}}(\sigma)R(\sigma)\\
\end{pmatrix}
\begin{pmatrix}
\theta(\sigma)\\
\underline{\mathbf{w}}(\sigma)\\
\end{pmatrix}\\
&=:
\widehat{A}_{\mathrm{nf}}(\sigma)
\begin{pmatrix}
\theta(\sigma)\\
\underline{\mathbf{w}}(\sigma)\\
\end{pmatrix},
\eld
\eeq
where
\beq
\label{e:Rsigma}
\begin{matrix}
\widehat{F}(\sigma):&(\ell^2)^n&\longrightarrow & \C\\
&\underline{\mathbf{w}}&\longmapsto&
-\llangle\underline{\mathbf{w}},\mathscr{F}_n(\rme^{-\rmi\sigma x}\mathbf{u}_{\mathrm{ad}})\rrangle,\\
&&&\\
R(\sigma):&(\ell^1)^n&\longrightarrow&\C\\
&\underline{\mathbf{w}}&\longmapsto&\rmi\frac{\sin{\pi\sigma}}{\pi}
(\sum_\ell(-1)^\ell\mathbf{w}_\ell,D\mathbf{u}_{\mathrm{ad}}^\prime(\pi)).
\end{matrix}
\eeq
We now conclude that
\beq
\label{e:foulin}
\widehat{A}_{\mathrm{nf}}=\int_{-\frac{1}{2}}^{\frac{1}{2}}\widehat{A}_{\mathrm{nf}}(\sigma)\rmd\sigma
=\int_{-\frac{1}{2}}^{\frac{1}{2}}
\mathscr{L}_{\mathrm{phd}}(\sigma)^{-1}\widehat{A}_{\mathrm{ch}}(\sigma)\mathscr{L}_{\mathrm{phd}}(\sigma)\rmd\sigma,
\eeq
where
\beqs
\begin{matrix}
\mathscr{L}_{\mathrm{phd}}(\sigma):&\C\times(\ell^2)^n(\sigma)&\longrightarrow&(\ell^2)^n\\
&(\theta,\underline{\mathbf{w}})&\longmapsto&\theta\widehat{\underline{\mathbf{E}}}(\sigma)+\underline{\mathbf{w}}.
\end{matrix}
\eeqs
Here $(\ell^2)^n(\sigma)=\{\underline{\mathbf{w}}\in(\ell^2)^n\mid
\llangle\underline{\mathbf{w}},
\mathscr{F}_n(\rme^{-\rmi\sigma x}\mathbf{u}_{\mathrm{ad}})\rrangle=0\}$.
We also recall that $\widehat{A}_{\mathrm{ch}}(\sigma)$ is defined in \eqref{e:Asigma}
and $\underline{\widehat{\mathbf{E}}}(\sigma)$ defined in \eqref{e:Esigma}.
\br
We note that for any $\underline{\mathbf{u}}\in L^2(\T_1,(\ell^2)^n)$ and
$\underline{\mathbf{v}}\in \mathscr{D}(\widehat{A}_{\mathrm{ch}})$,
$$(\widehat{F}\underline{\mathbf{u}})(\sigma)=\widehat{F}(\sigma)\underline{\mathbf{u}}(\sigma),\quad
(\mathscr{F}\delta_+\Gamma\mathscr{F}_{\mathrm{ch}}^{-1}\underline{\mathbf{v}})(\sigma)=R(\sigma)\underline{\mathbf{v}}(\sigma),
\text{ for a.e. }\sigma\in[-\frac{1}{2},\frac{1}{2}].$$
In addition, for any $(\theta,\underline{\mathbf{w}})\in L^2(\T_1)
\times L^2_\perp(\T_1,(\ell^2)^n)$,
$$\left(\widehat{\mathscr{L}}_{\mathrm{phd}}\begin{pmatrix}\theta\\ \underline{\mathbf{w}}\end{pmatrix}\right)(\sigma)
=\mathscr{L}_{\mathrm{phd}}(\sigma)\begin{pmatrix}\theta(\sigma)\\ \underline{\mathbf{w}}(\sigma)\end{pmatrix},
\text{ for a.e. }\sigma\in[-\frac{1}{2},\frac{1}{2}].$$
\er

We now consider the family of linear systems,
\beq
\label{e:3pl}
\begin{pmatrix}
\dot{\theta} \\ \dot{\underline{\mathbf{w}}}
\end{pmatrix}=
\widehat{A}_{\mathrm{nf}}(\sigma)\begin{pmatrix}
\theta \\ \underline{\mathbf{w}}
\end{pmatrix},
\text{ for all }\sigma\in[-\frac{1}{2},\frac{1}{2}].
\eeq
While we obtained these operators based on $L^2\slash\ell^2$ spaces, we can also consider them on $L^q\slash\ell^q$-based spaces.
To be more precise, we first define a family of projections
\beq
\label{e:Pqtilde}
\begin{matrix}
\widetilde{P}_q(\sigma): &Y_q & \longrightarrow & Y_q \\
&\underline{\mathbf{w}} &\longmapsto &\underline{\mathbf{w}}-\frac{1}{2\pi}\llangle\underline{\mathbf{w}},
\mathscr{F}_n(\rme^{-\rmi\sigma x}\mathbf{u}_{\mathrm{ad}})\rrangle\mathscr{F}_n(\rme^{-\rmi\sigma x}\mathbf{u}_{\star}^\prime),
\end{matrix}
\eeq
where
\beq
\label{e:space}
Y_q=
\bcs
(\ell^q)^n, &\text{for }1\leq q<\infty,\\
(\ell_0^\infty)^n,&\text{for } q=\infty.
\ecs
\eeq
Here we have $\ell_0^\infty=\{x\in \ell^\infty\mid\lim_{|n|\rightarrow \infty}|x_n|=0\}$ with the supremum norm.
For any $q\in[1,\infty]$, the projection $P_q(\sigma)$ is well-defined. In fact,
$\mathscr{F}_n(\rme^{-\rmi\sigma x}\mathbf{u}_{\mathrm{ad}}),  \mathscr{F}_n(\rme^{-\rmi\sigma x}\mathbf{u}_\star^\prime)\in Y_1$
since $\mathbf{u}_{\mathrm{ad}}(\pm\pi)=\mathbf{u}_\star^\prime(\pm\pi)=0$.
We now denote $\widetilde{Y}_{q,\mathrm{s}}(\sigma)=\rg{\widetilde{P}_q(\sigma)}$, and,
in the following lemma, define $\widehat{A}_{\mathrm{nf}}(\sigma)$ on $L^q\slash\ell^q$-based space.
\bl
\label{l:T2A2}
For $q\in [1,\infty]$ and $\sigma\in[-\frac{1}{2},\frac{1}{2}]$,
\beqs
\begin{matrix}
\mathscr{L}_{\mathrm{phd}}(\sigma):&\C\times \widetilde{Y}_{q,\mathrm{s}}(\sigma)&\longrightarrow&Y_q\\
&(\theta,\underline{\mathbf{w}})&\longmapsto&\theta\widehat{\underline{\mathbf{E}}}(\sigma)+\underline{\mathbf{w}},
\end{matrix}
\eeqs
is uniformly bounded and invertible with its inverse
\beqs
\begin{matrix}
\mathscr{L}_{\mathrm{phd}}(\sigma)^{-1}:&Y_q&\longrightarrow&\C\times \widetilde{Y}_{q,\mathrm{s}}(\sigma)\\
&\underline{\mathbf{v}}&\longmapsto&(\widehat{F}(\sigma)\underline{\mathbf{v}},
\underline{\mathbf{v}}-\widehat{\underline{E}}(\sigma)\widehat{F}(\sigma)\underline{\mathbf{v}}).
\end{matrix}
\eeqs
Moreover,
\beqs
\begin{array}{rcl}
\widehat{A}_{\mathrm{nf}}(\sigma): \C\times(\widetilde{Y}_{q,\mathrm{s}}(\sigma)\cap \mathscr{D}_q(\widehat{A}_{\mathrm{ch}}(\sigma)) )&
\to &\C\times\widetilde{Y}_{q,\mathrm{s}}(\sigma)
\end{array}
\eeqs
is well-defined and sectorial. Here $\mathscr{D}_q(\widehat{A}_{\mathrm{ch}}
(\sigma))=\{\underline{\mathbf{w}}\in Y_q\mid \{(1+m^2)\mathbf{w}_m\}_{m\in\Z}\in Y_q\}$
is the domain of $\widehat{A}_{\mathrm{ch}}(\sigma)$ in $Y_q$.
\el
\bpf
The assertions for $\mathscr{L}_{\mathrm{phd}}(\sigma)$ are straightforward. In order to show that
$\widehat{A}_{\mathrm{nf}}(\sigma)$ is well-defined, we recall the definition
of $\widehat{A}_{\mathrm{nf}}(\sigma)$ in \eqref{e:A2sigma}, which indicates that we only need to show
\beqs
\widehat{A}_{\mathrm{ch}}(\sigma)\underline{\widehat{\mathbf{E}}}(\sigma)\in\widetilde{Y}_{q,\mathrm{s}}(\sigma),\quad
\rg{(\widehat{A}_{\mathrm{ch}}(\sigma)-\underline{\widehat{\mathbf{E}}}(\sigma)R(\sigma))}\subseteq \widetilde{Y}_{q,\mathrm{s}}(\sigma).
\eeqs
We claim that $\widehat{A}_{\mathrm{ch}}(\sigma)\underline{\widehat{\mathbf{E}}}
(\sigma)\in\rg{(\widehat{A}_{\mathrm{ch}}(\sigma)-\underline{\widehat{\mathbf{E}}}(\sigma)R(\sigma))}$.
In fact, recall the definition of $R(\sigma)$ in \eqref{e:Rsigma} and
define $\mathbf{E}(\sigma, x):= (\sum_j \mathbf{E}_j(x)\rme^{-\rmi2\pi j\sigma})\rme^{-\rmi\sigma x}\in (C^\infty)^n(\mathscr{T}_{2\pi})$,
we have
\beqs
R(\sigma)\underline{\widehat{\mathbf{E}}}(\sigma)
=2\rmi\sin{\pi\sigma}(\mathbf{E}(\sigma, \pi), D\mathbf{u}_{\mathrm{ad}}^\prime(\pi))=0,
\eeqs
which means that
$\widehat{A}_{\mathrm{ch}}(\sigma)\underline{\widehat{\mathbf{E}}}(\sigma)
=(\widehat{A}_{\mathrm{ch}}(\sigma)-\underline{\widehat{\mathbf{E}}}(\sigma)R(\sigma))
(\underline{\widehat{\mathbf{E}}}(\sigma))\in\rg{(\widehat{A}_{\mathrm{ch}}(\sigma)-
\underline{\widehat{\mathbf{E}}}(\sigma)R(\sigma))}$.

We now only have to show $\rg{(\widehat{A}_{\mathrm{ch}}(\sigma)-\underline{\widehat{\mathbf{E}}}
(\sigma)R(\sigma))}\subseteq \widetilde{Y}_{q,\mathrm{s}}(\sigma)$.
Actually, for any $\underline{\mathbf{w}}\in \mathscr{D}_q(\widehat{A}_{\mathrm{ch}}(\sigma))$ with finitely many nonzero components, we have
\beqs
\bld
\llangle \widehat{A}_{\mathrm{ch}}(\sigma)\underline{\mathbf{w}}-\underline{\widehat{\mathbf{E}}}
(\sigma)R(\sigma)\underline{\mathbf{w}},\mathscr{F}_n(\rme^{-\rmi\sigma x}\mathbf{u}_{\mathrm{ad}})\rrangle
=&\llangle \widehat{A}_{\mathrm{ch}}(\sigma)\underline{\mathbf{w}},\mathscr{F}_n(\rme^{-\rmi\sigma x}\mathbf{u}_{\mathrm{ad}})\rrangle+2\pi
R(\sigma)\underline{\mathbf{w}}\\
=&2\pi\langle A(\rme^{\rmi\sigma x}\mathscr{F}_n^{-1}\underline{\mathbf{w}}), \mathbf{u}_{\mathrm{ad}}\rangle+
2\pi(\rme^{\rmi\sigma x}\mathscr{F}_n^{-1}\underline{\mathbf{w}}, D\mathbf{u}_{\mathrm{ad}}^\prime(x))|^{\pi}_{-\pi}\\
=&2\pi\langle \rme^{\rmi\sigma x}\mathscr{F}_n^{-1}\underline{\mathbf{w}}, B^*(0)\mathbf{u}_{\mathrm{ad}}\rangle=0,
\eld
\eeqs
and $\{\underline{\mathbf{w}}\in D_q(\widehat{A}_{\mathrm{ch}}(\sigma))|\underline{\mathbf{w}}\text{ has finite many nonzero elements}\}$ is
dense in $ \mathscr{D}_q(\widehat{A}_{\mathrm{ch}}(\sigma))$ under the graph norm of
$\widehat{A}_{\mathrm{ch}}(\sigma)$. Therefore, $\widehat{A}_{\mathrm{nf}}(\sigma)$ is well-defined.

Next, $\widehat{A}_{\mathrm{nf}}(\sigma)$ is sectorial, due to the facts that $\widehat{A}_
{\mathrm{nf}}(\sigma)=\mathscr{L}_{\mathrm{phd}}(\sigma)^{-1}\widehat{A}_{\mathrm{ch}}(\sigma)\mathscr{L}_{\mathrm{phd}}(\sigma)$
and $\widehat{A}_{\mathrm{ch}}(\sigma)$ is sectorial (for details, see Section \ref{ss:64} in the appendix).
\epf

Now we are ready to obtain the estimates for the time evolution of system \eqref{e:3pl}, for any given $\sigma\in[-\frac{1}{2},\frac{1}{2}]$.
Our discussion is split into the case $\sigma$ close to $0$ and the case $\sigma$ away from $0$.

For the case $\sigma\sim 0$, the derivation of the estimate relies on a diagonalized normal form, that is,
a complete separation of the netural and stable phase.
First, we notice that $\spec(\widehat{A}_{\mathrm{nf}}(\sigma))=\spec(\widehat{A}_{\mathrm{ch}}(\sigma))$
is independent of the choice of $q\in[1,\infty]$ and $\sigma\in[-\frac{1}{2},
\frac{1}{2}]$, which we will prove in Proposition \ref{p:601}. Moreover,  for $\sigma$ sufficiently small,
there is a unique continuation of the eigenvalue $0$, denoted as $\lambda(\sigma)$. The set
$\Lambda_1:=\{\lambda(\sigma)\}$ is a spectral set; see Section \ref{ss:65}, \ref{ss:66} for detailed treatment. Hence, let
\beq
\label{e:Pq}
\begin{matrix}
P_q(\sigma): &Y_q&\longrightarrow &Y_q\\
&\underline{\mathbf{w}}&\longmapsto&\underline{\mathbf{w}}-\frac{1}{2\pi}\llangle\underline{\mathbf{w}},
\mathscr{F}_n(\mathbf{e}^*(\sigma))\rrangle\mathscr{F}_n(\mathbf{e}(\sigma))
\end{matrix}
\eeq
be the spectral projection associated with $\Lambda_2:=\spec(\widehat{A}_{\mathrm{ch}}(\sigma))\backslash \{\lambda(\sigma)\}$.
Here $\mathbf{e}(\sigma)$ (respectively, $\mathbf{e}^*(\sigma)$) is the eigenvector of the Bloch wave operator
$B(\sigma)$ (respectively, the adjoint operator $B^*(\sigma)$) according to $\lambda(\sigma)$ with
\beq
\label{e:e}
\mathbf{e}(0)=\mathbf{u}^\prime_\star,\quad \mathbf{e}^*(0)=\mathbf{u}_{\mathrm{ad}},\quad \langle\mathbf{e}(\sigma),\mathbf{e}^*(\sigma)\rangle=1.
\eeq
We refer to Section \ref{ss:65} in the appendix for more details on $\mathbf{e}(\sigma)$ and $\mathbf{e}^*(\sigma)$.
We now denote
\beq
\label{e:cs}
\bld
Y_{q,\mathrm{c}}(\sigma)=\span\{\mathbf{e}(\sigma)\}, \quad Y_{q,\mathrm{s}}(\sigma)=\rg{P_q(\sigma)}, \\
\widehat{A}_{\mathrm{ch}}(\sigma)|_{Y_{q,\mathrm{c}}(\sigma)}=\widehat{A}_{\mathrm{c}}(\sigma),
\quad \widehat{A}_{\mathrm{ch}}(\sigma)|_{Y_{q,\mathrm{s}}}=\widehat{A}_{\mathrm{s}}(\sigma).
\eld
\eeq
We then introduce the following diagonalized operator
\beq
\label{e:Ascr}
\widehat{A}_{\mathrm{dg}}(\sigma)=
\begin{pmatrix}
\lambda(\sigma)&0\\
0&\widehat{A}_{\mathrm{s}}(\sigma)\\
\end{pmatrix}.
\eeq
It is not hard to conclude that
for $\sigma$ sufficiently small,
\beqs
\widehat{A}_{\mathrm{dg}}(\sigma): \C\times (Y_{q,\mathrm{s}}(\sigma)\cap \mathscr{D}^q(\widehat{A}_{\mathrm{ch}}(\sigma)))
\to \C\times Y_{q,\mathrm{s}}(\sigma)
\eeqs
is a well-defined operator.

The key step here is to find an invertible bounded linear transformation
\beq
\label{e:T3}
\begin{array}{rcl}
\widehat{\mathscr{T}}_{\mathrm{dg}}(\sigma)=
\begin{pmatrix}
\widehat{T}_{00}(\sigma)&\widehat{T}_{01}(\sigma)\\
\widehat{T}_{10}(\sigma)&\widehat{T}_{11}(\sigma)\\
\end{pmatrix}: \C\times\widetilde{Y}_{q,\mathrm{s}}(\sigma) \to  \C\times Y_{q,\mathrm{s}}(\sigma)
\end{array}
\eeq
such that
$\widehat{\mathscr{T}}_{\mathrm{dg}}(\sigma)\widehat{A}_{\mathrm{nf}}(\sigma)=
\widehat{A}_{\mathrm{dg}}(\sigma)\widehat{\mathscr{T}}_{\mathrm{dg}}(\sigma)$.
We note that the choice of $\widehat{\mathscr{T}}_{\mathrm{dg}}(\sigma)$ is not unique since there are nontrivial invertible
operators that commute with $\widehat{A}_{\mathrm{dg}}(\sigma)$.
\bl
\label{l:nft}
For $\sigma$ sufficiently small (that is, $|\sigma|\leq \gamma_0$) and $q\in[1,\infty]$,
\beq
\label{e:T3con}
\widehat{\mathscr{T}}_{\mathrm{dg}}(\sigma)=
\begin{pmatrix}
\mu(\delta)& S(\sigma)|_{\widetilde{Y}_{q,\mathrm{s}}(\sigma)}\\
P_q(\sigma)\underline{\widehat{\mathbf{E}}}(\sigma)&P_q(\sigma)|_{\widetilde{Y}_{q,\mathrm{s}}(\sigma)}\\
\end{pmatrix}
\eeq
satisfies the relation $\widehat{\mathscr{T}}_{\mathrm{dg}}(\sigma)\widehat{A}_
{\mathrm{nf}}(\sigma)=\widehat{A}_{\mathrm{dg}}(\sigma)\widehat{\mathscr{T}}_{\mathrm{dg}}(\sigma)$. Here
 we have that $\mu(\sigma)=-\frac{1}{2\pi}\llangle\widehat{\underline{\mathbf{E}}}(\sigma),
 \mathscr{F}_n(\mathbf{e}^*(\sigma))\rrangle$ and
\beqs
\begin{matrix}
 S(\sigma):&Y_q&\longrightarrow&\C\\
 &\underline{\mathbf{w}}&\longmapsto&-\frac{1}{2\pi}
 \llangle\underline{\mathbf{w}},\mathscr{F}_n(\mathbf{e}^*(\sigma))\rrangle.
\end{matrix}
\eeqs
Moreover, we have
\beq
\label{e:T3inv}
\widehat{\mathscr{T}}_{\mathrm{dg}}^{-1}=(\widehat{T}_{00}-\widehat{T}_{01}\widehat{T}_{11}^{-1}\widehat{T}_{10})^{-1}
\begin{pmatrix}
1&-\widehat{T}_{01}\widehat{T}_{11}^{-1}\\
-\widehat{T}_{11}^{-1}\widehat{T}_{10}&(\widehat{T}_{00}
-\widehat{T}_{01}\widehat{T}_{11}^{-1}\widehat{T}_{10})\widehat{T}_{11}^{-1}+\widehat{T}_{11}^{-1}
\widehat{T}_{10}\widehat{T}_{01}\widehat{T}_{11}^{-1}\\
\end{pmatrix},
\eeq
in which we suppress $\sigma$-dependence for simplicity.
\el
\bpf
 We recall from \eqref{e:A2sigma} that
 \beqs
 \widehat{A}_{\mathrm{nf}}(\sigma)=\begin{pmatrix} \widehat{F}(\sigma) \\
 \id-\underline{\widehat{\mathbf{E}}}(\sigma)\widehat{F}(\sigma) \end{pmatrix}
\widehat{A}_{\mathrm{ch}}(\sigma)
\begin{pmatrix} \underline{\widehat{\mathbf{E}}}(\sigma) & \id \end{pmatrix}.
 \eeqs
Therefore, in order to find a $\widehat{\mathscr{T}}_{\mathrm{dg}}$ as required,
we only need to find an invertible bounded linear operator
\beqs
\widehat{\mathscr{T}}_{\mathrm{int}}(\sigma)=
\begin{pmatrix}
 \widehat{\mathscr{T}}_1(\sigma) \\
 \widehat{\mathscr{T}}_2(\sigma)
\end{pmatrix}: Y_q\longrightarrow \C\times Y_{q,\mathrm{s}}(\sigma)
\eeqs
such that
\beqs
\widehat{\mathscr{T}}_{\mathrm{int}}(\sigma)\widehat{A}_{\mathrm{ch}}(\sigma)
=\widehat{A}_{\mathrm{dg}}(\sigma)\widehat{\mathscr{T}}_{\mathrm{int}}(\sigma),
\eeqs
which is equivalent to
\beqs
\begin{cases}
 \widehat{\mathscr{T}}_1(\sigma)(\lambda(\sigma)-\widehat{A}_{\mathrm{ch}}(\sigma))=0,\\
 \widehat{\mathscr{T}}_2(\sigma)\widehat{A}_{\mathrm{ch}}(\sigma)-\widehat{A}_{\mathrm{ch}}
 (\sigma)\widehat{\mathscr{T}}_2(\sigma)=0.
\end{cases}
\eeqs
While the choice of $\widehat{\mathscr{T}}_{1\slash2}(\sigma)$ satisfying the above equation is apparently not unique,
we choose that $\widehat{\mathscr{T}}_1(\sigma)=S(\sigma)$ and $\widehat{\mathscr{T}}_2(\sigma)=P_q(\sigma)$.
As a result, we have
\beqs
\widehat{\mathscr{T}}_{\mathrm{dg}}(\sigma)=
\begin{pmatrix}
 \widehat{\mathscr{T}}_1(\sigma) \\
 \widehat{\mathscr{T}}_2(\sigma)
\end{pmatrix}
\begin{pmatrix} \underline{\widehat{\mathbf{E}}}(\sigma) & \id|_{\widetilde{Y}_{q,\mathrm{s}}(\sigma)} \end{pmatrix}
=
 \begin{pmatrix}
 \mu(\delta)& S(\sigma)|_{\widetilde{Y}_{q,\mathrm{s}}(\sigma)}\\
 P_q(\sigma)\underline{\widehat{\mathbf{E}}}(\sigma)&P_q(\sigma)|_{\widetilde{Y}_{q,\mathrm{s}}(\sigma)}\\
 \end{pmatrix}.
\eeqs
To show that $\widehat{\mathscr{T}}_{\mathrm{dg}}(\sigma)^{-1}$ in (\ref{e:T3inv}) is correct, we only need to verify that
\beqs
\widehat{\mathscr{T}}_{\mathrm{dg}}(\sigma)^{-1}\widehat{\mathscr{T}}_{\mathrm{dg}}(\sigma)=
\begin{pmatrix}
 \id & 0\\ 0 & \id|_{\widetilde{Y}_{q,\mathrm{s}}(\sigma)}
\end{pmatrix},\quad
\widehat{\mathscr{T}}_{\mathrm{dg}}(\sigma)\widehat{\mathscr{T}}_{\mathrm{dg}}(\sigma)^{-1}=
\begin{pmatrix}
 \id & 0\\ 0 & \id|_{Y_{q,\mathrm{s}}(\sigma)}
\end{pmatrix},
\eeqs
which is clearly true.
\epf

Based on this lemma, we now derive the estimate for $\rme^{\widehat{A}_{\mathrm{nf}}(\sigma)t}$
when $\sigma$ is close to zero. We first introduce new notation
$M(t, \sigma):=\rme^{\widehat{A}_{\mathrm{nf}}t}$ and $\mathscr{M}(t):=\rme^{\widetilde{A}_{\mathrm{nf}}t}$ with
\beq
\label{e:M}
\begin{matrix}
 M(t, \sigma)=\begin{pmatrix}M_{00}(t, \sigma)&M_{01}(t, \sigma)\\M_{10}(t, \sigma)&M_{11}(t, \sigma)\end{pmatrix}
 :& \C\times \widetilde{Y}_{q,\mathrm{s}}(\sigma)&\longrightarrow&\C\times \widetilde{Y}_{q,\mathrm{s}}(\sigma),\\
 \mathscr{M}(t)=\begin{pmatrix}\mathscr{M}_{00}(t)&\mathscr{M}_{01}(t)\\\mathscr{M}_{10}(t)&\mathscr{M}_{11}(t)\end{pmatrix}:&
 \ell^2\times\ell^2_\perp(\Z,(L^2(\T_{2\pi}))^n)&\longrightarrow&\ell^2\times\ell^2_\perp(\Z,(L^2(\T_{2\pi}))^n).\\
\end{matrix}
\eeq
To make sense of the derivatives and Taylor expansions with respect to $\sigma$ of entries in $\mathscr{T}_{\mathrm{dg}}(\sigma)$,
we extend $\widehat{T}_{01}(\sigma)$ and $\widehat{T}_{11}(\sigma)$ continuously as operators on $Y_q$, that is,
\beq
\label{e:new}
\widehat{T}_{01}(\sigma)=S(\sigma)\widetilde{P}_q(\sigma),\quad \widehat{T}_{11}(\sigma)=P_q(\sigma)\widetilde{P}_q(\sigma).
\eeq
The same argument applies to operators $\widehat{\mathscr{T}}_{\mathrm{dg}}^{-1}(\sigma)$ and $M(t,\sigma)$.
\bl
\label{l:pen0}
For $\sigma$ sufficiently small (that is, $|\sigma|\leq \gamma_0$) and $q\in[1,\infty]$, there exist positive constants $C(q)$
and $\widetilde{d}$ such that, for all $t\geq 0$,
\beq
\begin{pmatrix}
 |M_{00}(t,\sigma)|&\opnorm{M_{01}(t,\sigma)}_{Y_q\rightarrow\C}\\
\opnorm{M_{10}(t,\sigma)}_{\C\rightarrow Y_q} &\opnorm{M_{11}(t,\sigma)}_{Y_q}
\end{pmatrix}
\leq C(q)
\begin{pmatrix} 1 & \frac{1}{\sqrt{1+t}}  \\ \frac{1}{\sqrt{1+t}}  & \frac{1}{1+t} \end{pmatrix}\rme^{-\widetilde{d}\sigma^2t}.
\eeq
Moreover, we have a higher regularity result for $M_{11}(t,\sigma)$, that is,
for any given $\sigma\in[-\gamma_0,\gamma_0]$, $q\in[1,\infty]$ and $\alpha>0$,
there exists $C(q,\alpha)>0$ such that
\beqs
\opnorm{M_{11}(t,\sigma)}_{Y_q\rightarrow Y_q^\alpha}\leq C(q,\alpha)[(1+t^{-\alpha})\rme^{-\frac{\gamma_1}{2}t}+
\frac{1}{1+t}\rme^{-\frac{d}{2}\sigma^2 t}].
\eeqs
\el
\bpf
The idea is to evaluate $M(t,\sigma)=\rme^{\widehat{A}_{\mathrm{nf}}(\sigma)t}$ based on
$\rme^{\widehat{A}_{\mathrm{nf}}(\sigma)t}=\widehat{\mathscr{T}}_{\mathrm{dg}}(\sigma)^
{-1}\rme^{\widehat{A}_{\mathrm{dg}}(\sigma)t}\widehat{\mathscr{T}}_{\mathrm{dg}}(\sigma)$.
We first state the following estimate from Proposition \ref{p:605}: For all $q\in[1,\infty]$ and $\sigma\in[-\gamma_0,\gamma_0]$,
there exists a constant $C(q)>0$ such that
\beqs
|\rme^{\lambda(\sigma)t}|\leq C(q)\rme^{-\frac{d}{2}\sigma^2t},\quad
\opnorm{\rme^{\widehat{A}_{\mathrm{s}}(\sigma)t}}_q \leq C(q)\rme^{-\frac{\gamma_1}{2}t}.
\eeqs
To obtain estimates on $\widehat{\mathscr{T}}_{\mathrm{dg}}$ and its inverse,
we start by computing the Taylor expansions of entries in $\widehat{\mathscr{T}}_{\mathrm{dg}}(\sigma)$.
A straightforward calculation using \eqref{e:e}, \eqref{e:Esigma} and \eqref{e:E} shows that
\begin{equation*}
\bld
 \mathbf{e}(\sigma)&=\mathbf{u}_\star^\prime+\rmi\sigma\mathbf{e}_1+\mathcal{O}(\sigma^2), &\quad\quad
 \mathbf{e}^*(\sigma)&=\mathbf{u}_{\mathrm{ad}}+\rmi\sigma\mathbf{e}^*_1+\mathcal{O}(\sigma^2),\\
 \rme^{-\rmi \sigma x}&=1-\rmi \sigma x+\mathcal{O}(\sigma^2),&\quad\quad
 \underline{\widehat{\mathbf{E}}}(\sigma)&=\mathscr{F}_n
 (-\mathbf{u}_{\star}^\prime-2\pi\rmi\sigma\phi\mathbf{u}_{\star}^\prime+\rmi \sigma x\mathbf{u}_{\star}^\prime)+\mathcal{O}(\sigma^2),
\eld
\end{equation*}
where
$\mathbf{e}_1$(respectively, $\mathbf{e}_1^*$) is even and nonzero due to the fact that $B(0)\mathbf{e}_1=-2D\mathbf{u}_\star^{\prime\prime}$
(respectively, $B^*(0)\mathbf{e}_1^*=-2D\mathbf{u}_{\mathrm{ad}}^{\prime}$). Then,
plugging these expansions into $\widehat{\mathscr{T}}_{\mathrm{dg}}$, and
using \eqref{e:new}, we obtain
\beq
\label{e:T3exp}
\widehat{\mathscr{T}}_{\mathrm{dg}}(\sigma)=
\begin{pmatrix}
 \mu(\sigma)&S(\sigma)\widetilde{P}_q(\sigma)\\
 P_q(\sigma)\underline{\widehat{\mathbf{E}}}(\sigma)&P_q(\sigma)\widetilde{P}_q(\sigma)
\end{pmatrix}=
\begin{pmatrix}
 1&-2\pi\rmi \sigma\Psi\\
 -2\pi\rmi \sigma\mathscr{F}_n(\Phi) & P_q(0)
\end{pmatrix}+
\begin{pmatrix}
\mathcal{O}(\sigma^2)&
\mathcal{O}(\sigma^2)\\
\mathcal{O}(\sigma^2)&
\mathcal{O}(\sigma)
\end{pmatrix},
\eeq
where $\Phi(x)=-\frac{x}{2\pi}\mathbf{u}_{\star}^\prime+\phi\mathbf{u}_{\star}^\prime-\frac{\mathbf{e}_1}{2\pi}$ and
\beqs
\begin{matrix}
 \Psi: &Y_q&\longrightarrow &\C\\
 &\underline{\mathbf{w}}&\longmapsto&\frac{1}{4\pi^2}
 \llangle\underline{\mathbf{w}},\mathscr{F}_n(x\mathbf{u}_{\mathrm{ad}}+\mathbf{e}_1^*)\rrangle.
\end{matrix}
\eeqs
Therefore, for any $\sigma\in[-\gamma_0,\gamma_0]$, $q\in[1,\infty]$, there exist positive constants
$C(q)$ and $\widetilde{d}$ such that, for all $(\theta,\underline{\mathbf{w}})\in \C\times Y_q$, we have the following estimate.
\beqs
\bld
\begin{pmatrix}
 |M_{00}(t,\sigma)\theta|&|M_{01}(t,\sigma)\underline{\mathbf{w}}|\\
\|M_{10}(t,\sigma)\theta\|_{Y_q} &\|M_{11}(t,\sigma)\underline{\mathbf{w}}\|_{Y_q}
\end{pmatrix}
&\leq C(q)
\begin{pmatrix} 1 & |\sigma| \\|\sigma| & 1 \end{pmatrix}
\begin{pmatrix} \rme^{-\frac{d}{2}\sigma^2 t} & 0 \\ 0 & \rme^{-\frac{\gamma_1}{2}t} \end{pmatrix}
\begin{pmatrix} 1 & |\sigma| \\ |\sigma| & 1 \end{pmatrix}\begin{pmatrix}|\theta|\\ \|\underline{\mathbf{w}}\|_{Y_q}\end{pmatrix}
\\
&\leq C(q)\left[\begin{pmatrix} 1 & |\sigma| \\|\sigma| & |\sigma|^2 \end{pmatrix}\rme^{-\frac{d}{2}\sigma^2 t}+
\begin{pmatrix} |\sigma|^2 & |\sigma| \\|\sigma| & 1 \end{pmatrix}\rme^{-\frac{\gamma_1}{2}t}\right]
\begin{pmatrix}|\theta|\\ \|\underline{\mathbf{w}}\|_{Y_q}\end{pmatrix}\\
&\leq C(q)\rme^{-\tilde{d}\sigma^2t}\begin{pmatrix} 1 & \frac{1}{\sqrt{1+t}} \\ \frac{1}{\sqrt{1+t}}  & \frac{1}{1+t} \end{pmatrix}
\begin{pmatrix}|\theta|\\ \|\underline{\mathbf{w}}\|_{Y_q}\end{pmatrix}.
\eld
\eeqs
Using \eqref{e:T3}, \eqref{e:T3inv} and \eqref{e:M}, we now expand $M_{11}(t,\sigma)$ and obtain
\beqs
\bld
\|M_{11}(t,\sigma)\underline{\mathbf{W}}\|_{Y_q^\alpha}\leq &
C(\|\widehat{T}_{11}(\sigma)^{-1}\widehat{T}_{10}(\sigma)\|_{Y_q^\alpha}
|\rme^{\lambda(\sigma)t}\widehat{T}_{01}(\sigma)\underline{\mathbf{W}}|+
\|\widehat{T}_{11}(\sigma)^{-1}\rme^{\widehat{A}_{\mathrm{s}}(\sigma)t}\widehat{T}_{11}
(\sigma)\underline{\mathbf{W}}\|_{Y_q^\alpha}+\\
&\|\widehat{T}_{11}^{-1}(\sigma)\widehat{T}_{10}(\sigma)\|_{Y_q^\alpha}
|\widehat{T}_{01}(\sigma)\widehat{T}_{11}(\sigma)^{-1}
\rme^{\widehat{A}_{\mathrm{s}}(\sigma)t}\widehat{T}_{11}(\sigma)\underline{\mathbf{W}}|)\\
\leq& C(q,\alpha)[|\sigma|^2\rme^{-\frac{d}{2}\sigma^2t}+(t^{-\alpha}+1)\rme^{-\frac{\gamma_1}{2}t}
+|\sigma|^2\rme^{-\frac{\gamma_1}{2}t}]\|\underline{\mathbf{W}}\|_{Y_q}\\
\leq & C(q,\alpha)[(t^{-\alpha}+1)\rme^{-\frac{\gamma_1}{2}t}+
\frac{1}{1+t}\rme^{-\frac{d}{2}\sigma^2t}]\|\underline{\mathbf{W}}\|_{Y_q},
\eld
\eeqs
where in the second inequality we used \eqref{e:T3exp}, and Proposition \ref{p:605}.
\epf
\br
We point out that, in the above lemma, the estimate for $M_{10}(\sigma)$ can not be improved, since
 $\mathscr{F}_n(\Phi)\neq0$. In fact, due to the fact that
$\Phi(x)\in (\mathbf{C}^{\infty}(\mathbb{T}_{2\pi}))^n$ and $\phi\mathbf{u}_{\star}^\prime$ is a nonzero even function, we have
\beqs
B(0)\Phi=\frac{1}{2\pi}[-2D\mathbf{u}_\star^{\prime\prime}-B(0)\mathbf{e}_1+2\pi B(0)(\phi\mathbf{u}_{\star}^\prime)]
=B(0)(\phi\mathbf{u}_{\star}^\prime)\neq 0.
\eeqs
On the other hand, the estimate for $M_{01}(\sigma)$ can be improved given suitable additional assumptions.
For example, if we assume that $\mathbf{u}_{\mathrm{ad}}^\prime(\pm \pi)=0$,
then $x\mathbf{u}_{\mathrm{ad}}+\mathbf{e}_1^*$ is zero, which leads to a
better estimate.
\er

For the case $\sigma$ away from $0$, we have the following result.
\bl
\label{l:pea0}
For $\sigma$ away from zero (i.e., for $\gamma_0\leq|\sigma|\leq \frac{1}{2}$) and $q\in[1,\infty]$,
there exist constants $C(q), \gamma_2>0$ such that
\beq
\label{e:36}
\opnorm{M(t,\sigma)}_{\C\times Y_q}\leq C(q)\rme^{-\gamma_2t}
\eeq
Moreover, we also have a higher regularity estimate for $M_{11}(t,\sigma)$, that is,
for any given $\gamma_0\leq |\sigma| \leq \frac{1}{2}$, $q\in[1,\infty]$ and $\alpha>0$,
there exists $C(q,\alpha)>0$ such that
\beqs
\opnorm{M_{11}(t,\sigma)}_{Y_q\rightarrow Y_q^\alpha}\leq C(q,\alpha)(1+t^{-\alpha})\rme^{-\gamma_2t}.
\eeqs
\el
\bpf
Recall that $\rme^{\widehat{A}_{\mathrm{nf}}(\sigma)t}=\mathscr{L}_{\mathrm{phd}}
(\sigma)^{-1}\rme^{\widehat{A}_{\mathrm{ch}}(\sigma)t}\mathscr{L}_{\mathrm{phd}}(\sigma)$.
The inequality (\ref{e:36}) is true due to the uniform boundedness of $\mathscr{T}(\sigma)$ in Lemma \ref{l:T2A2}
and the fact that $\opnorm{\rme^{\widehat{A}_{\mathrm{ch}}(\sigma)t}}_{Y_q}\leq C(q)\rme^{-\gamma_2t}$,
for $\sigma$ away from 0, in Proposition \ref{p:605}.
Moreover, by the expressions of $\mathscr{L}_{\mathrm{phd}}(\sigma)$ and its inverse in Lemma \ref{l:T2A2}, we have
$M_{11}(t,\sigma)=(\id-\widehat{\underline{E}}(\sigma)\widehat{F}(\sigma))\rme^{\widehat{A}_{\mathrm{ch}}(\sigma)t}$. Applying Proposition \ref{p:605},
we conclude that
$$\opnorm{M_{11}(t,\sigma)}_{Y_q\rightarrow Y_q^\alpha}=
\opnorm{(\id-\widehat{\underline{E}}(\sigma)\widehat{F}(\sigma))\rme^{\widehat{A}_{\mathrm{ch}}(\sigma)t}}_{Y_q\rightarrow Y_q^\alpha}
\leq C(q,\alpha)(1+t^{-\alpha})\rme^{-\gamma_2t}.$$
\epf

Lemma \ref{l:pen0} and \ref{l:pea0} give the following proposition.
\bp[\bf Fourier-Bloch estimates]
\label{pro:31}
For any $\sigma\in[-\frac{1}{2},\frac{1}{2}]$, $q\in[1,\infty]$, there exist constants
$C(q)$, $c>0$ such that $\widehat{A}_{\mathrm{nf}}(\sigma)$ is sectorial and
\beq
\begin{pmatrix}
 |M_{00}(t,\sigma)|&\opnorm{M_{01}(t,\sigma)}_{Y_q\rightarrow\C}\\
\opnorm{M_{10}(t,\sigma)}_{\C\rightarrow Y_q} &\opnorm{M_{11}(t,\sigma)}_{Y_q}
\end{pmatrix}\leq C(q)
\begin{pmatrix} 1 & \frac{1}{\sqrt{t+1}} \\ \frac{1}{\sqrt{t+1}} & \frac{1}{t+1} \end{pmatrix}\rme^{-c\sigma^2t},
\text{ for all }t\geq 0.
\eeq
Moreover, we have a higher regularity estimate on $M_{11}(t,\sigma)$, that is,
for any $\sigma\in[-\frac{1}{2},\frac{1}{2}]$, $q\in[1,\infty]$ and $\alpha>0$, there exist constants
$C(q,\alpha)$, $\gamma>0$ such that
\beq
\opnorm{M_{11}(t,\sigma)}_{Y_q\rightarrow Y_q^\alpha}\leq C(q,\alpha)\left((1+t^{-\alpha})\rme^{-\gamma t}+
\frac{1}{1+t}\rme^{-\frac{d}{2}\sigma^2t}\right), \text{ for all }t>0.
\eeq
\ep
We also need the Fourier-Bloch estimates for the derivative $\partial_\sigma M(t, \sigma)$ in the following lemma.
\bp [\bf Fourier-Bloch estimates for derivatives]
\label{pro:32}
For any $\sigma\in[-\frac{1}{2},\frac{1}{2}]$, $q\in[1,\infty]$
and $\beta\in(\frac{1}{2}(1-\frac{1}{q}),1)$, there exist positive constants $C(q,\beta)$ and $\widetilde{c}$
such that, for all $t\geq 0$,
\beq
\label{e:de}
\begin{pmatrix}
 |(\partial_\sigma M)_{00}(t,\sigma)|&\opnorm{(\partial_\sigma M)_{01}(t,\sigma)}_{Y_q\rightarrow \C}\\
\opnorm{(\partial_\sigma M)_{10}(t,\sigma)}_{\C\rightarrow Y_q} &\opnorm{(\partial_\sigma M)_{11}(t,\sigma)}_{Y_q}
\end{pmatrix}
\leq C(q,\beta)\begin{pmatrix}1 & \frac{1}{\sqrt{1+t}}\\ \frac{1}{\sqrt{1+t}} & \frac{1}{1+t}\end{pmatrix}
(t^{\frac{1}{2}}+t^{1-\beta})\rme^{-\widetilde{c}\sigma^2t}.
\eeq
Moreover, we have a higher regularity estimate on $(\partial_\sigma M)_{11}(t,\sigma)$, that is,
for $\sigma \in[-\frac{1}{2},\frac{1}{2}]$, $q\in[1,\infty]$, $\beta\in(\frac{1}{2}(1-\frac{1}{q}),1)$
and $\alpha\in(0,1)$, there exist $C(q,\alpha,\beta)>0$ and $\widetilde{\gamma}>0$ such that
\beqs
\bld
\opnorm{(\partial_\sigma M)_{11}(t,\sigma)}_{Y_q\rightarrow Y_q^\alpha}\leq &
C(q,\alpha,\beta)\left(\frac{t^{\frac{1}{2}}+t^{1-\beta}}{1+t}\rme^{-\frac{d}{2}\sigma^2t}+
(t^{\frac{1}{2}-\alpha}+t^{1-\beta})\rme^{-\widetilde{\gamma}t}\right),\text{ for all }t>0.
\eld
\eeqs
\ep
\bpf
On the one hand, we take the partial derivative of the following system with respect to $\sigma$
\beqs
\begin{pmatrix} \theta(t,\sigma)  \\ \underline{\widehat{\mathbf{W}}}(t,\sigma) \end{pmatrix}=
M(t, \sigma)\begin{pmatrix} \theta(0,\sigma)  \\ \underline{\widehat{\mathbf{W}}}(0,\sigma) \end{pmatrix},
\eeqs
 and obtain
\beqs
\begin{pmatrix} \partial_\sigma\theta(t,\sigma)  \\ \partial_\sigma\underline{\widehat{\mathbf{W}}}(t,\sigma) \end{pmatrix}=
M(t, \sigma)\begin{pmatrix} \partial_\sigma\theta(0,\sigma)
\\ \partial_\sigma\underline{\widehat{\mathbf{W}}}(0,\sigma) \end{pmatrix}
+(\partial_\sigma M(t, \sigma))\begin{pmatrix} \theta(0,\sigma)  \\ \underline{\widehat{\mathbf{W}}}(0,\sigma) \end{pmatrix}.
\eeqs
On the other hand, we have that
\beqs
\begin{pmatrix} \dot{\underline{\theta}}(t,\sigma)  \\ \dot{\widehat{\underline{\mathbf{W}}}}(t,\sigma) \end{pmatrix}=
\widehat{A}_{\mathrm{nf}}(\sigma)\begin{pmatrix}\underline{\theta}(t,\sigma)  \\ \widehat{\underline{\mathbf{W}}}(t,\sigma) \end{pmatrix}.
\eeqs
Taking the partial derivative with respect to $\sigma$, the equation becomes
\beqs
\begin{pmatrix} \dot{(\partial_\sigma\underline{\theta})}(t,\sigma)
\\ \dot{(\partial_\sigma\widehat{\underline{\mathbf{W}}})}(t,\sigma) \end{pmatrix}=
\widehat{A}_{\mathrm{nf}}(\sigma)\begin{pmatrix}\partial_\sigma\underline{\theta}(t,\sigma)
\\ \partial_\sigma\widehat{\underline{\mathbf{W}}}(t,\sigma) \end{pmatrix}+
\widehat{A}_{\mathrm{nf}}^\prime(\sigma)\begin{pmatrix}\underline{\theta}(t,\sigma)  \\ \widehat{\underline{\mathbf{W}}}(t,\sigma) \end{pmatrix},
\eeqs
for which the variation of constant formula gives
\beqs
\begin{pmatrix} \partial_\sigma\theta(t,\sigma)  \\ \partial_\sigma\underline{\widehat{\mathbf{W}}}(t,\sigma) \end{pmatrix}=
M(t, \sigma)\begin{pmatrix} \partial_\sigma\theta(0,\sigma)
\\ \partial_\sigma\underline{\widehat{\mathbf{W}}}(0,\sigma) \end{pmatrix}+
\int_0^tM(t-s, \sigma) \widehat{A}_{\mathrm{nf}}^\prime(\sigma)M(s, \sigma)\begin{pmatrix}
\theta(0,\sigma)  \\ \underline{\widehat{\mathbf{W}}}(0,\sigma) \end{pmatrix}\rmd s.
\eeqs
Therefore, one has
\beq
\label{e:pwd}
\bld
&\partial_\sigma M(t, \sigma)=\int_0^tM(t-s, \sigma) \widehat{A}_{\mathrm{nf}}^\prime(\sigma)M(s, \sigma)\rmd s\\
=&\int_0^t\mathscr{L}_{\mathrm{phd}}(\sigma)^{-1}\rme^{\widehat{A}_{\mathrm{ch}}
(\sigma)(t-s)}\mathscr{L}_{\mathrm{phd}}(\sigma)(\mathscr{L}_{\mathrm{phd}}(\sigma)^{-1}
\widehat{A}_{\mathrm{ch}}(\sigma)\mathscr{L}_{\mathrm{phd}}(\sigma))^{\prime}\mathscr{L}_
{\mathrm{phd}}(\sigma)^{-1}\rme^{\widehat{A}_{\mathrm{ch}}(\sigma)s}\mathscr{L}_{\mathrm{phd}}(\sigma)\rmd s\\
=&\int_0^t\mathscr{L}_{\mathrm{phd}}(\sigma)^{-1}\rme^{\widehat{A}_{\mathrm{ch}}(\sigma)(t-s)}
\mathscr{N}(\sigma)\rme^{\widehat{A}_{\mathrm{ch}}(\sigma)s}\mathscr{L}_{\mathrm{phd}}(\sigma)\rmd s,
\eld
\eeq
where
\beq
\label{e:N}
\mathscr{N}(\sigma)=\widehat{A}_{\mathrm{ch}}^\prime(\sigma)+\widehat{A}_{\mathrm{ch}}
(\sigma)\widehat{\underline{E}}^\prime(\sigma)
\widehat{F}(\sigma)-\widehat{\underline{E}}^\prime(\sigma)\widehat{F}(\sigma)\widehat{A}_{\mathrm{ch}}(\sigma).
\eeq
We recall that $\widehat{\underline{E}}(\sigma)$ is defined in \eqref{e:Esigma}, $\widehat{F}(\sigma)$
in \eqref{e:Rsigma} and $\widehat{A}_{\mathrm{ch}}(\sigma)$ in \eqref{e:Asigma}.

For $|\sigma|\leq \gamma_0$, by Lemma \ref{l:nft} and the above equation \eqref{e:pwd},
we have
\beq
\label{e:Mde0}
\partial_\sigma M(t,\sigma)=
\int_0^t\widehat{\mathscr{T}}_{\mathrm{dg}}(\sigma)^{-1}\widetilde{\mathscr{N}}
(\sigma,t,s)\widehat{\mathscr{T}}_{\mathrm{dg}}(\sigma)\rmd s,
\eeq
where
\beq
\label{e:Ntilde}
\bld
\widetilde{\mathscr{N}}(\sigma,t,s)&=\rme^{\widehat{A}_{\mathrm{dg}}(\sigma)(t-s)}
\begin{pmatrix}S(\sigma)\\ P_q(\sigma)\end{pmatrix}
\mathscr{N}(\sigma)\begin{pmatrix}-\mathscr{F}_n\mathbf{e}(\sigma),\id \end{pmatrix}
\rme^{\widehat{A}_{\mathrm{dg}}(\sigma)s}\\
&=\begin{pmatrix}
-\rme^{\lambda(\sigma)t}S(\sigma)\mathscr{N}(\sigma)\mathscr{F}_n\mathbf{e}(\sigma)&
\rme^{\lambda(\sigma)(t-s)}S(\sigma)\mathscr{N}(\sigma)\rme^{\widehat{A}_{\mathrm{s}}(\sigma)s}\\
-\rme^{\lambda(\sigma)s}\rme^{\widehat{A}_{\mathrm{s}}(\sigma)(t-s)}P_q(\sigma)\mathscr{N}(\sigma)\mathscr{F}_n\mathbf{e}(\sigma)&
\rme^{\widehat{A}_{\mathrm{s}}(\sigma)(t-s)}P_q(\sigma)\mathscr{N}(\sigma)\rme^{\widehat{A}_{\mathrm{s}}(\sigma)s}
\end{pmatrix}\\
&=:\begin{pmatrix}\widetilde{N}_{00}(\sigma,t,s)&\widetilde{N}_{01}(\sigma,t,s)
\\\widetilde{N}_{10}(\sigma,t,s)&\widetilde{N}_{11}(\sigma,t,s)\end{pmatrix}.
\eld
\eeq
We now evaluate the entries of $\widetilde{\mathscr{N}}$ with expansions combining \eqref{e:N}
and \eqref{e:Ntilde}.
First, recall the definitions of $\widehat{A}_{\mathrm{dg}}(\sigma)$ in \eqref{e:Asigma}, $S(\sigma)$ in Lemma \ref{l:nft},
$P_q(\sigma)$ in \eqref{e:Pq}, $\mathscr{F}_n$ in \eqref{e:FT}, and $\mathbf{e}(\sigma)$ in \eqref{e:e}.

For $\widetilde{N}_{00}$, note that
it is smooth with respect to $\sigma$ and
\beqs
\widetilde{N}_{00}(0,t,s)=-S(0)\left(\widehat{A}_{\mathrm{ch}}^\prime(0)+\widehat{A}_
{\mathrm{ch}}(0)\widehat{\underline{E}}^\prime(0)
\widehat{F}(0)-\widehat{\underline{E}}^\prime(0)\widehat{F}(0)\widehat{A}_{\mathrm{ch}}(0)\right)\mathscr{F}_n\mathbf{e}(0).
\eeqs
We claim that $\widetilde{N}_{00}(0,t,s)=0$.
In fact, since $\widehat{A}_{\mathrm{ch}}^\prime(0)\mathscr{F}_n(\mathbf{e}(0))$ and $\widehat{A}_{\mathrm{ch}}(0)\widehat{\underline{E}}^\prime(0)$
are orthogonal to $\mathscr{F}_n(\mathbf{u}_{\mathrm{ad}})$ in $Y_2$,
$S(0)\widehat{A}_{\mathrm{ch}}^\prime(0)\mathscr{F}_n(\mathbf{e}(0))=0$ and $S(0)\widehat{A}_{\mathrm{ch}}(0)\widehat{\underline{E}}^\prime(0)=0$.
Moreover, $\widehat{F}(\sigma)\widehat{A}_{\mathrm{ch}}(\sigma)=R(\sigma)$,
which is defined in \eqref{e:Rsigma} with $R(0)=0$. Therefore, there exists a positive constant
$C$ such that
\beqs
|\widetilde{N}_{00}|\leq C|\sigma|\rme^{-\frac{d}{2}\sigma^2t}.
\eeqs
For $\widetilde{N}_{10}$, due to Proposition \ref{p:605}
and the fact that $\widetilde{N}_{10}$ is smooth in $\sigma$ with $\widetilde{N}_{10}(0,t,s)\neq 0$,
there exists a positive constant $C$ such that
\beqs
\opnorm{\widetilde{N}_{10}}_{\C\rightarrow Y_{q,\mathrm{s}}(\sigma)}
\leq C\rme^{-\frac{\gamma_1}{2}(t-s)}\rme^{-\frac{d}{2}\sigma^2s}
\leq C\rme^{-\frac{d}{2}\sigma^2t}\rme^{-\frac{\gamma_1}{4}(t-s)}.
\eeqs
For $\widetilde{N}_{01}$,
we have, for any $q\in[1,\infty]$ and $\beta>\frac{1}{2}(1-\frac{1}{q})$,
\beqs
\bld
\opnorm{\widetilde{N}_{01}}_{Y_{q,\mathrm{s}}(\sigma)\rightarrow \C}\leq&
C|\rme^{\lambda(\sigma)(t-s)}|\left(\opnorm{\widehat{A}_{\mathrm{ch}}^\prime(\sigma)\rme^{\widehat{A}_{\mathrm{s}}(\sigma)s}}
_{Y_{q,\mathrm{s}}(\sigma)\rightarrow Y_q}+\right.\\
&\left.|S(\sigma)\widehat{A}_{\mathrm{ch}}(\sigma)\widehat{\underline{E}}^\prime(\sigma)
|\opnorm{\rme^{\widehat{A}_{\mathrm{s}}(\sigma)s}}_{Y_{q,\mathrm{s}}(\sigma)}+\opnorm{\widehat{F}(\sigma)
\widehat{A}_{\mathrm{ch}}(\sigma)\rme^{\widehat{A}_{\mathrm{s}}(\sigma)s}}_{Y_{q,\mathrm{s}}(\sigma)\rightarrow \C}\right)\\
\leq &C\rme^{-\frac{d}{2}\sigma^2(t-s)}\left(\opnorm{\rme^{\widehat{A}_{\mathrm{s}}(\sigma)s}}
_{Y_{q,\mathrm{s}}(\sigma)\rightarrow Y_q^{\frac{1}{2}}}+|\sigma|\opnorm{\rme^{\widehat{A}_{\mathrm{s}}(\sigma)s}}_{Y_{q,\mathrm{s}}(\sigma)}+
|\sigma|\opnorm{\rme^{\widehat{A}_{\mathrm{s}}(\sigma)s}}_{Y_{q,\mathrm{s}}(\sigma)\rightarrow Y_1}\right)\\
\leq &C(\beta)\rme^{-\frac{d}{2}\sigma^2(t-s)}
\left(\opnorm{\rme^{\widehat{A}_{\mathrm{s}}(\sigma)s}}_{Y_{q,\mathrm{s}}(\sigma)\rightarrow Y_q^{\frac{1}{2}}}
+\opnorm{\rme^{\widehat{A}_{\mathrm{s}}(\sigma)s}}_{Y_{q,\mathrm{s}}(\sigma)}+
|\sigma|\opnorm{\rme^{\widehat{A}_{\mathrm{s}}(\sigma)s}}_{Y_{q,\mathrm{s}}(\sigma)\rightarrow Y_q^{\beta}}\right),
\eld
\eeqs
where the last inequality results from the fact that, for any $q\in[1,\infty]$ and $\beta>\frac{1}{2}(1-\frac{1}{q})$, we
have a continuous imbedding
\beqs Y_q^\beta \hookrightarrow Y_1. \eeqs
Now, using Proposition \ref{p:605} and \ref{p:604}, we can further conclude that
\beqs
\bld
\opnorm{\widetilde{N}_{01}}_{Y_{q,\mathrm{s}}(\sigma)\rightarrow \C}
\leq &C(q,\beta)(s^{-\frac{1}{2}}+1+|\sigma|s^{-\beta})\rme^{-\frac{d}{2}\sigma^2(t-s)}\rme^{-\frac{\gamma_1}{2}s}\\
\leq&
C(q,\beta)\rme^{-\frac{d}{2}\sigma^2t}(s^{-\frac{1}{2}}+|\sigma|+|\sigma|s^{-\beta})\rme^{-\frac{\gamma_1}{4}s}.
\eld
\eeqs
For $\widetilde{N}_{11}$, we have, for any $q\in[1,\infty]$ and $\beta>\frac{1}{2}(1-\frac{1}{q})$,
\beqs
\bld
\opnorm{\widetilde{N}_{11}}_{Y_{q,\mathrm{s}}(\sigma)}\leq&
C\opnorm{\rme^{\widehat{A}_{\mathrm{s}}(\sigma)(t-s)}}_{Y_{q,\mathrm{s}}(\sigma)}
\left(\opnorm{\widehat{A}_{\mathrm{ch}}^\prime(\sigma)
\rme^{\widehat{A}_{\mathrm{s}}(\sigma)s}}_{Y_{q,\mathrm{s}}(\sigma)\rightarrow Y_q}\right.\\
&\left.+\opnorm{\rme^{\widehat{A}_{\mathrm{s}}(\sigma)s}}_{Y_{q,\mathrm{s}}(\sigma)}+
\opnorm{\widehat{F}(\sigma)\widehat{A}_{\mathrm{ch}}(\sigma)
\rme^{\widehat{A}_{\mathrm{s}}(\sigma)s}}_{Y_{q,\mathrm{s}}(\sigma)\rightarrow \C}\right)\\
\leq &C(q,\beta)\left(s^{-\frac{1}{2}}+1+|\sigma|s^{-\beta}\right)\rme^{-\frac{\gamma_1}{2}t}.
\eld
\eeqs
Therefore, combining \eqref{e:Mde0}, \eqref{e:T3exp}, and the above estimates for entries,
we conclude that, for $|\sigma|\in[-\gamma_0,\gamma_0]$, $q\in[1,\infty]$ and $\beta\in(\frac{1}{2}(1-\frac{1}{q}),1)$,
there exist positive constants $C(q,\beta)$ and $c_1\leq \frac{d}{2}$ such that
\beq
\label{e:de1}
\bld
&\begin{pmatrix}
 |(\partial_\sigma M)_{00}(t,\sigma)|&\opnorm{(\partial_\sigma M)_{01}(t,\sigma)}_{Y_q\rightarrow \C}\\
\opnorm{(\partial_\sigma M)_{10}(t,\sigma)}_{\C\rightarrow Y_q} &\opnorm{(\partial_\sigma M)_{11}(t,\sigma)}_{Y_q}
\end{pmatrix}\\
\leq&
\begin{pmatrix}1 & |\sigma|\\ |\sigma| & 1\end{pmatrix}
\int_0^t\begin{pmatrix}
|\widetilde{N}_{00}(\sigma,t,s)|&\opnorm{\widetilde{N}_{01}(\sigma,t,s)}_{Y_{q,\mathrm{s}}(\sigma)\rightarrow \C}\\
\opnorm{\widetilde{N}_{10}(\sigma,t,s)}_{\C\rightarrow Y_{q,\mathrm{s}}(\sigma)}&\opnorm{\widetilde{N}_{11}(\sigma,t,s)}_{Y_{q,\mathrm{s}}(\sigma)}
\end{pmatrix}\rmd s
\begin{pmatrix}1 & |\sigma|\\ |\sigma| & 1\end{pmatrix}\\
\stackrel{*}{\leq}&C(q,\beta)\begin{pmatrix}1 & |\sigma|\\ |\sigma| & 1\end{pmatrix}
\begin{pmatrix}|\sigma| t & \frac{t^{\frac{1}{2}}}{\sqrt{1+t}}+|\sigma|\frac{t^{1-\beta}}{(1+t)^{1-\beta}}\\
\frac{t^{\frac{1}{2}}}{\sqrt{1+t}} & \frac{t^{\frac{1}{2}}+t^{1-\beta}}{1+t}\end{pmatrix}
\begin{pmatrix}1 & |\sigma|\\ |\sigma| & 1\end{pmatrix}\rme^{-\frac{d}{2}\sigma^2t}\\
\leq &C(q,\beta)\begin{pmatrix}1 & \frac{1}{\sqrt{1+t}}\\ \frac{1}{\sqrt{1+t}} & \frac{1}{1+t}\end{pmatrix}
(t^{\frac{1}{2}}+t^{1-\beta})\rme^{-c_1\sigma^2t}.
\eld
\eeq
Here the inequality $(*)$ relies on the fact that for any $\beta\in(0,1)$, there exists a positive constant $C(\beta)$
such that
$$
\int_0^t \rme^{-\frac{\gamma_1}{4}s}\rmd s\leq C(\beta)\frac{t^{\beta}}{(1+t)^{\beta}},\quad
\int_0^ts^{-\beta}\rme^{-\frac{\gamma_1}{4}s}\rmd s\leq C(\beta)\frac{t^{1-\beta}}{(1+t)^{1-\beta}}.
$$

On the other hand, for $\gamma_0\leq|\sigma|\leq \frac{1}{2}$, $q\in[1,\infty]$ and $\beta\in(\frac{1}{2}(1-\frac{1}{q}),1)$,
by the expression \eqref{e:pwd} and Proposition \ref{p:605},
there exist positive constants $C(q,\beta)$ and $c_2$ such that
\beq
\label{e:de2}
\bld
\opnorm{\partial_{\sigma}M(t,\sigma)}_{\C\times Y_q}\leq&
C(q)
\int_0^t\opnorm{\rme^{\widehat{A}_{\mathrm{ch}}(\sigma)(t-s)}\mathscr{N}(\sigma)\rme^{\widehat{A}_{\mathrm{ch}}(\sigma)s}}_{q}\rmd s\\
\leq&C(q,\beta)
\rme^{-\gamma_2t}\int_0^t(s^{-\frac{1}{2}}+1+|\sigma|s^{-\beta})\rmd s\\
\leq &
C(q,\beta)
(t^{\frac{1}{2}}+t+|\sigma|t^{1-\beta})\rme^{-\gamma_2t}.
\eld
\eeq
By \eqref{e:de1} and \eqref{e:de2},we now conclude that, for any $\sigma\in[-\frac{1}{2},\frac{1}{2}]$, $q\in[1,\infty]$
and $\beta\in(\frac{1}{2}(1-\frac{1}{q}),1)$, there exists positive constant $C(q,\beta)$ and $\widetilde{c}$
such that
\beqs
\begin{pmatrix}
 |(\partial_\sigma M)_{00}(t,\sigma)|&\opnorm{(\partial_\sigma M)_{01}(t,\sigma)}_{Y_q\rightarrow \C}\\
\opnorm{(\partial_\sigma M)_{10}(t,\sigma)}_{\C\rightarrow Y_q} &\opnorm{(\partial_\sigma M)_{11}(t,\sigma)}_{Y_q}
\end{pmatrix}\leq
C(q,\beta)\begin{pmatrix}1 & \frac{1}{\sqrt{1+t}}\\ \frac{1}{\sqrt{1+t}} & \frac{1}{1+t}\end{pmatrix}
(t^{\frac{1}{2}}+t^{1-\beta})\rme^{-\widetilde{c}\sigma^2t}.
\eeqs
We now consider $(\partial_{\sigma}M)_{11}(t,\sigma)$.
For $\sigma\in[-\gamma_0,\gamma_0]$, we plug $\widehat{\mathscr{T}}_{\mathrm{dg}}(\sigma)$ from \eqref{e:T3},
$\widehat{\mathscr{T}}_{\mathrm{dg}}(\sigma)^{-1}$ from \eqref{e:T3inv}, and
$\widetilde{N}(\sigma,t,s)$ from the last equality in \eqref{e:Ntilde} into \eqref{e:Mde0}. We then obtain
\beqs
(\partial_{\sigma}M)_{11}(t,\sigma)=\int_0^t-\widehat{T}_{11}^{-1}\widehat{T}_{10}
\left(\widetilde{N}_{00}\widehat{T}_{01}+\widetilde{N}_{01}\widehat{T}_{11}\right)
+\widehat{T}_{11}^{-1}
\left(\widetilde{N}_{10}\widehat{T}_{01}+\widetilde{N}_{11}\widehat{T}_{11}\right)\rmd s \cdot \left(1+\rmO(\sigma)\right)
\eeqs
More precisely, for
$q\in[1,\infty]$ and $\alpha\in(0,1)$, there exists $C$ such that
\beqs
\bld
\|(\partial_\sigma M)_{11}(t,\sigma)\underline{\mathbf{W}}\|_{Y_q^\alpha}\leq &
C\int_{0}^t\left[\|\widehat{T}_{11}^{-1}\widehat{T}_{10}\|_{Y_q^\alpha}
\left(|\widetilde{N}_{00}||\widehat{T}_{01}\underline{\mathbf{W}}|+
|\widetilde{N}_{01}\widehat{T}_{11}\underline{\mathbf{W}}|\right)\right.\\
&\left.+\|\widehat{T}_{11}^{-1}\widetilde{N}_{10}\|_{Y_q^\alpha}|\widehat{T}_{01}\underline{\mathbf{W}}|+
\|\widehat{T}_{11}^{-1}\widetilde{N}_{11}\widehat{T}_{11}\underline{\mathbf{W}}\|_{Y_q^\alpha}\right]\rmd s\\
\stackrel{**}{\leq}&C\int_{0}^t\left[\|\widehat{T}_{10}\|_{Y_q^\alpha}
\left(|\widetilde{N}_{00}||\widehat{T}_{01}\underline{\mathbf{W}}|+
|\widetilde{N}_{01}\widehat{T}_{11}\underline{\mathbf{W}}|\right)\right.\\
&\left.+\|\widetilde{N}_{10}\|_{Y_q^\alpha}|\widehat{T}_{01}\underline{\mathbf{W}}|+
\|\widetilde{N}_{11}\widehat{T}_{11}\underline{\mathbf{W}}\|_{Y_q^\alpha}\right]\rmd s.
\eld
\eeqs
Here the inequality $(**)$ relies on the fact that
\beqs
\begin{matrix}
\widehat{T}_{11}(\sigma): & Y_q^\alpha & \longrightarrow & Y_q^\alpha \\
& \underline{\mathbf{v}} & \longmapsto & \underline{\mathbf{v}}-\langle\mathbf{e}(\sigma),\rme^{-\rmi\sigma x}\mathbf{u}_{\mathrm{ad}}\rangle^{-1}
\llangle\underline{\mathbf{v}},
\mathscr{F}_n(\rme^{-\rmi \sigma x}\mathbf{u}_{\mathrm{ad}})\rrangle
\mathscr{F}_n\mathbf{e}(\sigma)
\end{matrix}
\eeqs
is a uniformly bounded operator for $q\in[1,\infty]$ and $\alpha\in(0,1)$.
Using the explicit expressions of the entries of $\widehat{\mathscr{T}}_{\mathrm{dg}}(\sigma)$ in \eqref{e:T3con}
and the estimates on the entries of $\widetilde{\mathscr{N}}$ as shown above, we derive the following estimates,
\beqs
\bld
&\int_0^t\|\widehat{T}_{10}\|_{Y_q^\alpha}|\widetilde{N}_{00}||\widehat{T}_{01}\underline{\mathbf{W}}|\rmd s
\leq C(q,\alpha)|\sigma|^3t\rme^{-\frac{d}{2}\sigma^2t}\|\underline{\mathbf{W}}\|_{Y_q},\\
&\int_0^t\|\widehat{T}_{10}\|_{Y_q^\alpha}|\widetilde{N}_{01}\widehat{T}_{11}\underline{\mathbf{W}}|\rmd s
\leq C(q,\alpha,\beta)|\sigma|\left(\frac{t^{\frac{1}{2}}}{\sqrt{1+t}}+|\sigma|\frac{t^{1-\beta}}{(1+t)^{1-\beta}}\right)
\rme^{-\frac{d}{2}\sigma^2t}\|\underline{\mathbf{W}}\|_{Y_q},\\
&\int_0^t\|\widetilde{N}_{10}\|_{Y_q^\alpha}|\widehat{T}_{01}\underline{\mathbf{W}}|\rmd s\leq
C(q,\alpha) \frac{t^{\frac{1}{2}}}{\sqrt{1+t}}\rme^{-\frac{d}{2}\sigma^2 t}\|\underline{\mathbf{W}}\|_{Y_q},\\
&\int_0^t \|\widetilde{N}_{11}\widehat{T}_{11}\underline{\mathbf{W}}\|_{Y_q^\alpha}\rmd s \leq
C(q,\alpha,\beta)\left(t+t^{1-\beta}+\int_0^t(t-s)^{-\alpha}s^{-1/2}\rmd s\right)\rme^{-\frac{\gamma_1}{2}t}\|\underline{\mathbf{W}}\|_{Y_q}.
\eld
\eeqs
We now conclude that, for $\sigma\in[-\gamma_0,\gamma_0]$,
$q\in[1,\infty]$, $\beta\in(\frac{1}{2}(1-\frac{1}{q}),1)$
and $\alpha\in(0,1)$, there exists $C(q,\alpha,\beta)>0$ such that
\beqs
\opnorm{(\partial_\sigma M)_{11}(t,\sigma)}_{Y_q\rightarrow Y_q^\alpha}
\leq C(q,\alpha,\beta)[\frac{t^{\frac{1}{2}}+t^{1-\beta}}{1+t}\rme^{-\frac{d}{2}\sigma^2t}+
(t^{\frac{1}{2}-\alpha}+t^{1-\beta})\rme^{-\frac{\gamma_1}{2}t}].
\eeqs
For $\gamma_0\leq |\sigma |\leq \frac{1}{2}$, $q\in[1,\infty]$,
$\beta\in(\frac{1}{2}(1-\frac{1}{q}),1)$ and $\alpha\in(0,1)$, there exists $C(q,\alpha,\beta)>0$ such that
\beqs
\bld
\opnorm{(\partial_\sigma M)_{11}(t,\sigma)}_{Y_q\rightarrow Y_q^\alpha}= &
\opnorm{(\id-\widehat{\underline{E}}(\sigma)\widehat{F}(\sigma))
\rme^{\widehat{A}_{\mathrm{ch}}(\sigma)(t-s)}\mathscr{N}(\sigma)\rme^{\widehat{A}_{\mathrm{ch}}(\sigma)s}}_{Y_q\rightarrow Y_q^\alpha}\\
\leq &C(q,\alpha,\beta)\rme^{-\gamma_2t}\int_0^t(t-s)^{-\alpha}s^{-\frac{1}{2}}+s^{-\frac{1}{2}}
+1+|\sigma|s^{1-\beta}\rmd s\\
\leq &C(q,\alpha,\beta)(t^{\frac{1}{2}-\alpha}+t^{1-\beta})\rme^{-\frac{\gamma_2}{2}t}.
\eld
\eeqs
Altogether, for $\sigma \in[-\frac{1}{2},\frac{1}{2}]$, $q\in[1,\infty]$,
$\beta\in(\frac{1}{2}(1-\frac{1}{q}),1)$ and $\alpha\in(0,1)$, there exist $C(q,\alpha,\beta)>0$ and $\widetilde{\gamma}>0$ such that
\beqs
\bld
\opnorm{(\partial_\sigma M)_{11}(t,\sigma)}_{Y_q\rightarrow Y_q^\alpha}\leq &
C(q,\alpha,\beta)[\frac{t^{\frac{1}{2}}+t^{1-\beta}}{1+t}\rme^{-\frac{d}{2}\sigma^2t}+
(t^{\frac{1}{2}-\alpha}+t^{1-\beta})\rme^{-\widetilde{\gamma}t}],\text{ for all }t>0.
\eld
\eeqs
\epf

\section{Linear estimates in physical space}\label{s:5}

According to the outline at the beginning of Section \ref{s:4}, we are now ready to derive the linear estimates for
$\rme^{A_{\mathrm{nf}}t}$. To be more precise, we first show by Fubini's Theorem that
\beqs
\mathscr{M}(t)\begin{pmatrix} \underline{\theta}  \\ \underline{\mathbf{W}} \end{pmatrix}
=\check{M}(t)\ast\begin{pmatrix} \underline{\theta}  \\ \underline{\mathbf{W}} \end{pmatrix},
\eeqs
where $\check{M}(t)$ is the generalized ``inverse Fourier transform'' of $M(t,\sigma)$. We then employ an argument similar to, but more
intricate than, Young's inequality for the case of the scalar heat equation, exploiting the linear Fourier-Bloch estimates in
Proposition \ref{pro:31} and \ref{pro:32}, to obtain the general $L^p$--$L^q$ estimate on our linear normal form $\rme^{A_{\mathrm{nf}}t}$.

To this end, we first note that $A_{\mathrm{nf}}=\widetilde{A}_{\mathrm{nf}}|_{\ell^1\times X_{\mathrm{ch}}^{\perp}}$ and thus we have, by \eqref{e:lpf},
for any $(\underline{\theta},\underline{\mathbf{W}})\in \ell^1\times X_{\mathrm{ch}}^{\perp}$,
\beqs
\mathscr{M}(t)\begin{pmatrix}\underline{\theta}\\ \underline{\mathbf{W}}\end{pmatrix}=
\rme^{A_{\mathrm{nf}}t}\begin{pmatrix}\underline{\theta}\\ \underline{\mathbf{W}}\end{pmatrix}=
\rme^{\widetilde{A}_{\mathrm{nf}}t}\begin{pmatrix}\underline{\theta}\\ \underline{\mathbf{W}}\end{pmatrix}=
\mathscr{F}_{\mathrm{nf}}^{-1}\rme^{\widehat{A}_{\mathrm{nf}}t}\mathscr{F}_{\mathrm{nf}}
\begin{pmatrix}\underline{\theta}\\ \underline{\mathbf{W}}\end{pmatrix},
\text{ for all }t>0.
\eeqs
Recall the notation $\mathscr{M}(t)=\rme^{A_{\mathrm{nf}}t}$, the definition of $\mathscr{F}_{\mathrm{nf}}$ from \eqref{e:FT}, and
the definition of $\widetilde{A}_{\mathrm{nf}}$, $\widehat{A}_{\mathrm{nf}}$ from \eqref{e:315}.
In addition, by \eqref{e:foulin}, we have,
for any $(\theta(\sigma), \widehat{\underline{\mathbf{W}}}(\sigma))\in L^2(\T_1)\times L^2_{\perp}(\T_1,\ell^2)$,
\beqs
\left(\rme^{\widehat{A}_{\mathrm{nf}}t}\begin{pmatrix}\theta\\ \widehat{\underline{\mathbf{W}}}\end{pmatrix}\right)(\sigma)=
\rme^{\widehat{A}_{\mathrm{nf}}(\sigma)t}\begin{pmatrix}\theta(\sigma)\\ \widehat{\underline{\mathbf{W}}}(\sigma)\end{pmatrix}=
\mathscr{L}_{\mathrm{phd}}(\sigma)^{-1}\rme^{\widehat{A}_{\mathrm{ch}}(\sigma)t}\mathscr{L}_{\mathrm{phd}}(\sigma)
\begin{pmatrix}\theta(\sigma)\\ \widehat{\underline{\mathbf{W}}}(\sigma)\end{pmatrix},
\text{ for a.e. }\sigma\in[-\frac{1}{2},\frac{1}{2}].
\eeqs

To show that $\rme^{A_{\mathrm{nf}}t}$ is a generalized convolution, we first define
$M(t,\sigma)$'s ``generalized inverse Fourier transform''
$\check{M}(t):=\begin{pmatrix}\check{\underline{M}}_{00}&\check{\underline{M}}_{01}\\
\check{\underline{M}}_{10}&\check{\underline{M}}_{11}\\\end{pmatrix}$,
with expressions as follows
\beq
\label{e:41}
\bld
\check{\underline{M}}_{00}(t):=&\{\check{M}_{00}(t,j)\}_{j\in\Z}
:=\{\int^{\frac{1}{2}}_{-\frac{1}{2}}M_{00}(t, \sigma)\rme^{\rmi 2\pi \sigma j}\rmd\sigma\}_{j\in\Z},\\
\check{\underline{M}}_{01}(t,y):=&\{\check{M}_{01}(t,y,j)\}_{j\in\Z}:=\{\int^{\frac{1}{2}}_{-\frac{1}{2}}
\sum_{\ell\in\Z}(M_{01})_{\ell}(t,\sigma)\rme^{-\rmi(\sigma+\ell) y}\rme^{\rmi 2\pi j\sigma}\rmd\sigma\}_{j\in\Z},\\
\check{\underline{M}}_{10}(t,x):=&\{\check{M}_{10}(t,x,j)\}_{j\in\Z}
:=\{\frac{1}{2\pi}\int^{\frac{1}{2}}_{-\frac{1}{2}}\sum_{\ell\in\Z}(M_{10})_\ell(t,\sigma)\rme^{\rmi(\sigma+\ell)x}
\rme^{\rmi 2\pi j\sigma}\rmd\sigma\}_{j\in\Z},\\
\check{\underline{M}}_{11}(t,x,y):=&\{\check{M}_{11}(t,x,y,j)\}_{j\in\Z}
:=\{\frac{1}{2\pi}\int^{\frac{1}{2}}_{-\frac{1}{2}}\sum_{\ell,\eta\in\Z}(M_{11})_{\ell\eta}(t,\sigma)\rme^{\rmi (\sigma+\ell)x}
\rme^{-\rmi (\sigma+\eta)y}\rme^{\rmi 2\pi j\sigma}\rmd\sigma\}_{j\in\Z}.\\
\eld
\eeq
We then have the following lemma.
\bl
For any
$(\underline{\theta},\underline{\mathbf{W}})\in \ell^1\times X_{\mathrm{ch}}^{\perp}$ and all $t>0$,
\beq
\label{e:42}
\mathscr{M}(t) \begin{pmatrix} \underline{\theta}  \\ \underline{\mathbf{W}} \end{pmatrix}
=\check{M}(t)\ast\begin{pmatrix} \underline{\theta}  \\ \underline{\mathbf{W}} \end{pmatrix}
=\begin{pmatrix} \check{\underline{M}}_{00}\ast\underline{\theta}&
\check{\underline{M}}_{01}\ast\underline{\mathbf{W}}  \\
\check{\underline{M}}_{10}*\underline{\theta}& \check{\underline{M}}_{11}\ast\underline{\mathbf{W}}\end{pmatrix},
\eeq
where
\beqs
\bld
\check{\underline{M}}_{00}\ast\underline{\theta}=&\{\sum_{k\in\mathbb{Z}}\check{M}_{00}(t, j-k)\theta_k\}
_{j\in\mathbb{Z}}, \\
\check{\underline{M}}_{01}\ast\underline{\mathbf{W}}=&
\{\sum_{k\in\Z}\int_{-\pi}^\pi\check{M}_{01}(t,y,j-k)\mathbf{W}_k(y)\rmd y\}_{j\in\mathbb{Z}},\\
\check{\underline{M}}_{10}*\underline{\theta}=&\{\sum_{k\in\mathbb{Z}}\check{M}_{10}(t,x,j-k)\theta_k\}_{j\in\Z},\\
\check{\underline{M}}_{11}\ast\underline{\mathbf{W}}=&
\{\sum_{k\in\Z}\int_{-\pi}^\pi\check{M}_{11}(t,x,y,j-k)\mathbf{W}_k(y)\rmd y\}_{j\in\Z}.
\eld
\eeqs
\el
\bpf
The proof is a straightforward application of Fubini's theorem.
\epf

We are now ready to obtain the general $L^p-L^q$ linear estimates on $\mathscr{M}(t)$. We denote
\beqs X_q=(L^q(\Z,L^q(\T_{2\pi})))^n, \text{ for any }q\in[1,\infty], \eeqs
and prove the following proposition.
\bp[\textbf{general $L^p$--$L^q$ estimates}]
\label{pro:41}
For any $1\leq q \leq p \leq \infty$ and $(\underline{\theta},\underline{\mathbf{W}})\in \ell^1\times X_{\mathrm{ch}}^{\perp}$,
there exists a positive constant $C$ such that, for all $t>0$,
\beq
\label{e:le}
\begin{pmatrix}
\|\mathscr{M}_{00}(t)\underline{\theta}\|_{\ell^p}&\|\mathscr{M}_{01}(t)\underline{\mathbf{W}}\|_{\ell^p}\\
\|\mathscr{M}_{10}(t)\underline{\theta}\|_{X_p}&\|\mathscr{M}_{11}(t)\underline{\mathbf{W}}\|_{X_p}
\end{pmatrix}
\leq C
\begin{pmatrix}
(1+t)^{-\frac{1}{2}(\frac{1}{q}-\frac{1}{p})}\|\underline{\theta}\|_{\ell^q}
&(1+t)^{-\frac{1}{2}(\frac{1}{q}-\frac{1}{p})-\frac{1}{2}}\|\underline{\mathbf{W}}\|_{X_q}\\
(1+t)^{-\frac{1}{2}(\frac{1}{q}-\frac{1}{p})-\frac{1}{2}}\|\underline{\theta}\|_{\ell^q}
&t^{-\frac{1}{2}(\frac{1}{q}-\frac{1}{p})}(1+t)^{-1}\|\underline{\mathbf{W}}\|_{X_q}
\end{pmatrix}.
\eeq
\ep

\bpf
We illustrate the derivation of the estimates on $\mathscr{M}_{01}$ and sketch the estimates on $\mathscr{M}_{00}$ and
$\mathscr{M}_{10}$. Lastly, we show the estimates for $\mathscr{M}_{11}$.

We first notice that, for any $\underline{\mathbf{W}}\in X_{\mathrm{ch}}^{\perp}$ and $1\leq q, r\leq p\leq \infty$
satisfying $1+\frac{1}{p}=\frac{1}{q}+\frac{1}{r}$,
there exists a positive constant $C$ such that
\beq
\label{e:likeyoung}
\|\mathscr{M}_{01}(t)\underline{\mathbf{W}}\|_{\ell^p}
\leq C\|\check{\underline{M}}_{01}(t)\|_{X_\infty}^{\frac{1}{q}-\frac{1}{p}}\left(
\sum_{j}\sup_{|y|\leq\pi}|\check{M}_{01}(t,y,j)|\right)^{\frac{1}{r}}\|\underline{\mathbf{W}}\|_{X_q}.
\eeq
In fact, by H\"{o}lder's inequality, we have
\beqs
\bld
\|\mathscr{M}_{01}(t)\underline{\mathbf{W}}\|_{\ell^p}^p=&\|\check{\underline{M}}_{01}\ast\underline{\mathbf{W}}\|_{\ell^p}^p
=\sum_{j\in\mathbb{Z}}|\sum_{k\in\Z}\int_{-\pi}^\pi\check{M}_{01}(t,y,j-k)\mathbf{W}_k(y)\rmd y|^p\\
\leq &\sum_{j\in\mathbb{Z}}\left(\sum_{k\in\Z}\int_{-\pi}^{\pi}|\check{M}_{01}(t,y,j-k)|^{1-\frac{r}{p}}|\mathbf{W}_k(y)|^{1-\frac{q}{p}}
\left(|\check{M}_{01}(t,y,j-k)|^r|\mathbf{W}_k(y)|^q\right)^{\frac{1}{p}}\rmd y\right)^p\\
\leq &\|\check{\underline{M}}_{01}(t)\|_{X_r}^{p-r}\|\underline{\mathbf{W}}\|_{X_q}^{p-q}
\sum_{j,k\in\mathbb{Z}}\int_{-\pi}^\pi|\check{\underline{M}}_{01}(t,y,j-k)|^r|\mathbf{W}_k(y)|^q\rmd y\\
\leq&\|\check{\underline{M}}_{01}(t)\|_{X_r}^{p-r}\left(\sup_{|y|\leq\pi}\sum_{j\in\Z}|\check{M}_{01}(t,y,j)|^r\right)
\|\underline{\mathbf{W}}\|_{X_q}^{p}\\
\leq& \left[(2\pi)^{1-\frac{1}{q}}\|\check{\underline{M}}_{01}(t)\|_{X_\infty}^{\frac{1}{q}-\frac{1}{p}}\left(
\sup_{|y|\leq\pi}\sum_{j\in\Z}|\check{M}_{01}(t,y,j)|\right)^{\frac{1}{r}}\|\underline{\mathbf{W}}\|_{X_q}\right]^p.
\eld
\eeqs
Moreover, by (\ref{e:41}), we have
\beq\label{e:inftynorm}
\|\check{\underline{M}}_{01}(t)\|_{X_\infty}\leq \sup_{|y|\leq\pi}|\int^{\frac{1}{2}}_{-\frac{1}{2}}
|\sum_{\ell\in\Z}(M_{01})_\ell(t,\sigma)
\rme^{-\rmi(\sigma+\ell) y}|\rmd\sigma|
\leq \frac{C(\infty)}{\sqrt{1+t}}\int^{\frac{1}{2}}_{-\frac{1}{2}}\rme^{-c\sigma^2t}\rmd\sigma
\leq\frac{C}{1+t}.
\eeq
Here we use the fact that any bounded linear functional on $\ell_0^\infty$ can be viewed
as a bounded linear functional on $\ell^\infty$ with the same norm.
We now estimate the $X_1$ norm of $\{\check{M}_{01}(t,y,j)\}_{j\in\mathbb{Z}}$. By using Proposition \ref{pro:32},
there exists $C>0$, independent of the choice of $y\in[-\pi,\pi]$, such that
\beq\label{e:onenorm}
\bld
\sum_{j\neq 0}|\check{M}_{01}(t,y,j)|=&\sum_{j\neq 0}\left(1+\frac{(j-\frac{y}{2\pi})^2}{t}\right)^{-\frac{1}{2}}\left(1+
\frac{(j-\frac{y}{2\pi})^2}{t}\right)^{\frac{1}{2}}|\check{M}_{01}(t,y,j)|\\
\leq& C\left(\int_\R\frac{1}{1+\frac{x^2}{t}}\rmd x\right)^{\frac{1}{2}}\left[\sum_{j}(1+\frac{(j-
\frac{y}{2\pi})^2}{t})|\check{M}_{01}(t,y,j)|^2\right]^{\frac{1}{2}}\\
\leq&Ct^{\frac{1}{4}}\left(\int^{\frac{1}{2}}_{-\frac{1}{2}}\sum_{\alpha=0}^1t^{-\alpha}
|\sum_{\ell\in\Z}(\partial^\alpha_\sigma (M_{01})_\ell(t,\sigma))
\rme^{-\rmi(\sigma+\ell)y}|^2\rmd\sigma\right)^{\frac{1}{2}}\\
\stackrel{***}{\leq}&
Ct^{\frac{1}{4}}\left(\int_{-\frac{1}{2}}^{\frac{1}{2}}\frac{\rme^{-2c\sigma^2t}}{1+t}\rmd\sigma+\frac{1}{t}
\int_{-\frac{1}{2}}^{\frac{1}{2}}\frac{(t^{\frac{1}{2}}+t^{1-\frac{3}{4}})^2
\rme^{-2\widetilde{c}\sigma^2t}}{1+t}\rmd\sigma\right)^{\frac{1}{2}}\\
\leq &C\frac{t^{\frac{1}{4}}+1}{(1+t)^{\frac{3}{4}}}
\leq \frac{C}{\sqrt{1+t}},\text{ for all }t>0.
\eld
\eeq
Here in the inequality (***), we applied Proposition \ref{pro:32}
with $q=\infty$ and $\beta=\frac{3}{4}$ (actually, any fixed $\beta\in(\frac{1}{2}, \frac{3}{4}]$).
Combining \eqref{e:likeyoung}, \eqref{e:inftynorm}, and \eqref{e:onenorm}, we have that, for all $1\leq q\leq p\leq \infty$ and
$\underline{\mathbf{W}}\in X_{\mathrm{ch}}^{\perp}$, there exists a positive constant $C$ such that
\[
\|\sqrt{1+t}\mathscr{M}_{01}(t)\underline{\mathbf{W}}\|_{\ell^p}\leq \frac{C}{(1+t)^{\frac{1}{2}(\frac{1}{q}-\frac{1}{p})}}
\|\underline{\mathbf{W}}\|_{X_q},\text{ for all }t\geq 0.
\]
For $\mathscr{M}_{00}$, the steps are the same as above but easier. For $\mathscr{M}_{10}$, we point out two main differences
to the above calculation. First, instead of \eqref{e:likeyoung}, we use
\beqs
\|\mathscr{M}_{10}(t)\underline{\theta}\|_{X_p}
\leq C\|\check{\underline{M}}_{10}(t)\|_{X_\infty}^{\frac{1}{q}-\frac{1}{p}}\left(\int_{-\pi}^\pi
\left(\sum_{j}|\check{M}_{10}(t,x,j)|\right)^{\frac{p}{r}}\rmd x\right)^{\frac{1}{p}}\|\underline{\theta}\|_{\ell^q}.
\eeqs
Second, to estimate the $Y_1$ norm of $\{\check{M}_{10}(t,x,j)\}_{j\in\Z}$, we use Proposition \ref{pro:32}
with $q=1$ and $\beta=\frac{1}{2}$(actually, any fixed $\beta\in(0,\frac{3}{4}]$), instead of
$q=\infty$ and $\beta=\frac{3}{4}$.

The last step of the proof consists of deriving the estimates for $\mathscr{M}_{11}$. We first have
\beqs
\bld
\|\mathscr{M}_{11}(t)\underline{\mathbf{W}}\|_{X_p}
\leq & (2\pi)^{\frac{1}{r}}\left(\sup_{|x|,|y|\leq |\pi|}\sup_{j\in\Z}|\check{M}_{11}(t,x,y,j)|\right)^{\frac{1}{q}-\frac{1}{p}}
\left(\sup_{|x|,|y|\leq\pi}\sum_{j\in\Z}|\check{M}_{11}(t,x,y,j)|\right)^
\frac{1}{r}\|\underline{\mathbf{W}}\|_{X_q}.
\eld
\eeqs
On the one hand, we apply Proposition \ref{pro:31} with $q=\infty$ and $\alpha>\frac{1}{2}$ and have
\beqs
\bld
\sup_{|x|,|y|\leq |\pi|}\sup_{j\in\Z}|\check{M}_{11}(t,x,y,j)|\leq &
\sup_{|x|,|y|\leq |\pi|}
\int_{-\frac{1}{2}}^{\frac{1}{2}}|\sum_{\ell,\eta\in\Z}(M_{11})_{\ell\eta}(t,\sigma)
\rme^{\rmi (\sigma+\ell)x}\rme^{-\rmi (\sigma+\eta)y}|\rmd \sigma\\
\leq & C(\alpha)\int_{-\frac{1}{2}}^{\frac{1}{2}}\opnorm{M_{11}(t,\sigma)}_{Y_\infty\rightarrow Y_\infty^\alpha}\rmd \sigma\\
\leq &\frac{C(\alpha)}{t^{\alpha}(1+t)^{\frac{3}{2}-\alpha}}.
\eld
\eeqs
On the other hand, by applying Proposition \ref{pro:31} and \ref{pro:32} with $q=\infty$, $\alpha\in(\frac{1}{2},1)$ and $\beta=\frac{3}{4}$,
there exists $C(\alpha)$, independent of choices of $x, y\in[-\pi,\pi]$, such that
\beqs
\bld
\sum_{|j|>1}|\check{M}_{11}(t,x,y,j)|= &
\sum_{|j|>1}\left(1+\frac{(j+\frac{x-y}{2\pi})^2}{t}\right)^{-\frac{1}{2}}
\left(1+\frac{(j+\frac{x-y}{2\pi})^2}{t}\right)^{\frac{1}{2}}|\check{M}_{11}(t,x,y,j)|\\
\leq&Ct^{\frac{1}{4}}\left(\int^{\frac{1}{2}}_{-\frac{1}{2}}\sum_{\alpha=0}^1t^{-\alpha}
|\sum_{\ell,\eta\in\Z}\left(\partial^\alpha_\sigma (M_{11})_{\ell\eta}(t,\sigma)\right)
\rme^{\rmi(\sigma+\ell)x}\rme^{-\rmi(\sigma+\eta)y}|^2\rmd\sigma\right)^{\frac{1}{2}}\\
\leq &C(\alpha)t^{\frac{1}{4}}\left(\int_{-\frac{1}{2}}^{\frac{1}{2}}\opnorm{M_{11}(t,\sigma)}_
{Y_\infty\rightarrow Y_\infty^\alpha}^2\rmd \sigma
+\int_{-\frac{1}{2}}^{\frac{1}{2}}\opnorm{(\partial_\sigma M)_{11}(t,\sigma)}_
{Y_\infty\rightarrow Y_\infty^\alpha}^2\rmd \sigma\right)^{\frac{1}{2}}\\
\leq&C(\alpha)t^{\frac{1}{4}}\left(\frac{1}{t^{2\alpha}(1+t)^{\frac{5}{2}-2\alpha}}\right)^{\frac{1}{2}}\\
\leq &C(\alpha)\frac{1}{t^{\alpha-\frac{1}{4}}(1+t)^{\frac{5}{4}-\alpha}},\text{ for all }t>0.
\eld
\eeqs
Moreover, combining the above two estimates, we have that, for given $\alpha\in(\frac{1}{2},1)$, there exists
$C(\alpha)>0$ such that
\beqs
\bld
\sup_{|x|,|y|\leq\pi}\sum_{j\in\Z}|\check{M}_{11}(t,x,y,j)|\leq \frac{C(\alpha)}{t^{\alpha}(1+t)^{1-\alpha}}.
\eld
\eeqs
Therefore, for any $1\leq q\leq p\leq\infty$, $\alpha\in(\frac{1}{2},1)$ and $\underline{\mathbf{W}}\in X_{\mathrm{ch}}^{\perp}$,
there exists $C(\alpha)>0$ such that
\beqs
\|\mathscr{M}_{11}(t)\underline{\mathbf{W}}\|_{X_p}\leq \frac{C(\alpha)}{(1+t)^{\frac{1}{2}(\frac{1}{q}-\frac{1}{p})}}
\frac{1}{t^\alpha(1+t)^{1-\alpha}}\|\underline{\mathbf{W}}\|_{X_q}.
\eeqs
Moreover, we can improve the above estimate for $t$ close to zero. Note that for the Laplacian operator, we
have the general $L^p$-$L^q$ estimate for all $t>0$. As a perturbation of the Laplacian operator, $\mathscr{M}_{11}$
has the same estimate for sufficiently small $t$, which can be seen by using the variation of constant formula as follows.
\beqs
\|\mathscr{M}_{11}(t)\underline{\mathbf{W}}\|_{X_p}=\|(\id-\mathbf{E}*F)\rme^{A_{\mathrm{ch}} t}\underline{\mathbf{W}}\|_{X_p}
\leq C \|\rme^{A_{\mathrm{ch}} t}\underline{\mathbf{W}}\|_{X_p}=C \|\rme^{A t}\mathbf{W}\|_{L^p},
\eeqs
where $\underline{\mathbf{W}}=\{\mathbf{W}_j(x)\}_{j\in\Z}$ and $\mathbf{W}(2\pi j+x)=\mathbf{W}_j(x)$ for all $j\in\Z$ and $x\in[-\pi,\pi]$. We now
let $\mathbf{V}(t,x)=\rme^{A t}\mathbf{W}(x)$ and have
\beqs
\mathbf{V}(t)=\rme^{D\partial_{xx} t}\mathbf{W}+\int_{0}^t \rme^{D\partial_{xx} (t-s)}\mathbf{f}^\prime(\mathbf{u}_\star)\mathbf{V}(s)\rmd s.
\eeqs
from which we derive
\beqs
\sup_{0<t\leq T}t^{\frac{1}{2}(\frac{1}{q}-\frac{1}{p})}\|\mathbf{V}(t)\|_{L^p}\leq \|\mathbf{W}\|_{L^q}+
CT^{1-\frac{1}{2}(\frac{1}{q}-\frac{1}{p})}\sup_{0<t\leq T}t^{\frac{1}{2}(\frac{1}{q}-\frac{1}{p})}\|\mathbf{V}(t)\|_{L^p}.
\eeqs
Taking $T$ sufficiently small such that $CT^{1-\frac{1}{2}(\frac{1}{q}-\frac{1}{p})}\leq \frac{1}{2}$, we obtain
\beqs
\|\rme^{A t}\mathbf{W}\|_{L^p}\leq \frac{C}{t^{\frac{1}{2}(\frac{1}{q}-\frac{1}{p})}}\|\mathbf{W}\|_{L^q},
\eeqs
which implies that
\beqs
\|\mathscr{M}_{11}(t)\underline{\mathbf{W}}\|_{X_p}
\leq \frac{C}{t^{\frac{1}{2}(\frac{1}{q}-\frac{1}{p})}}\|\underline{\mathbf{W}}\|_{X_q}, \text{ for all }0<t\leq T.
\eeqs
Therefore, for any $1\leq q\leq p\leq\infty$, there exists $C>0$ such that
\beqs
\|\mathscr{M}_{11}(t)\underline{\mathbf{W}}\|_{X_p}\leq \frac{C}{t^{\frac{1}{2}(\frac{1}{q}-\frac{1}{p})}}
\frac{1}{(1+t)}\|\underline{\mathbf{W}}\|_{X_q}.
\eeqs
\epf
\br
\label{r:led}
By \eqref{e:41},\eqref{e:42} and a similar argument as in Proposition \ref{pro:41}, it is not hard to conclude that,
for any $j\in\Z^+$, $1\leq q \leq p \leq \infty$ and $(\underline{\theta},\underline{\mathbf{W}})\in \ell^1\times X_{\mathrm{ch}}^{\perp}$,
there exists a positive constant $C$ such that, for all $t>0$,
\beq
\label{e:led}
\bld
&\begin{pmatrix}
\|\delta_+^j\mathscr{M}_{00}(t)\underline{\theta}\|_{\ell^p}&\|\delta_+^j\mathscr{M}_{01}(t)\underline{\mathbf{W}}\|_{\ell^p}\\
\|\delta_+^j\mathscr{M}_{10}(t)\underline{\theta}\|_{X_p}&\|\delta_+^j\mathscr{M}_{11}(t)\underline{\mathbf{W}}\|_{X_p}
\end{pmatrix}\\
&\quad\quad\quad\quad\quad\quad\quad\quad\leq  C
\begin{pmatrix}
(1+t)^{-\frac{1}{2}(\frac{1}{q}-\frac{1}{p}+j)}\|\underline{\theta}\|_{\ell^q}
&(1+t)^{-\frac{1}{2}(\frac{1}{q}-\frac{1}{p}+j+1)}\|\underline{\mathbf{W}}\|_{X_q}\\
(1+t)^{-\frac{1}{2}(\frac{1}{q}-\frac{1}{p}+j+1)}\|\underline{\theta}\|_{\ell^q}
&t^{-\frac{1}{2}(\frac{1}{q}-\frac{1}{p})}(1+t)^{-(1+\frac{j}{2})}\|\underline{\mathbf{W}}\|_{X_q}
\end{pmatrix}.
\eld
\eeq
\er

\section{Maximal regularity and nonlinear stability}\label{s:6}

In this section, we prove the main theorem--Theorem \ref{t:1}.
To achieve this, we first introduce a Banach space that our argument will be based on.
We then collect several maximal regularity results since the normal form system is quasilinear.
Based on our normal form and the general $L^p-L^q$ linear estimates, we can apply a fixed point argument to
the variation of constant formula, thus obtaining the nonlinear stability result.

We choose $r\in (4,+\infty)$ and define
\beqs
Z=\{(\underline{\theta},\underline{\mathbf{W}})\in C((0,+\infty), \ell^1\times (X_{\mathrm{ch}}\cap\mathscr{T}_{\mathrm{ch}}^{-1}(H^2)))
\mid \|(\underline{\theta},\underline{\mathbf{W}})\|_Z< \infty \},
\eeqs
where
\beqs
\bld
\|(\underline{\theta},\underline{\mathbf{W}})\|_Z=&\sup_{t> 0}\|\underline{\theta}(t)\|_{\ell^1}+
\sup_{t> 0}(1+t)^{\frac{1}{2}}\|\underline{\theta}(t)\|_{\ell^\infty}+
\sup_{t>0}(1+t)^{\frac{5}{4}}\|\delta^2\underline{\theta}\|_{\ell^2}\\
&+\sup_{t>0}(1+t)^{\frac{1}{2}}\|\underline{\mathbf{W}}\|_{X_1}+\sup_{t>0}(1+t)\|\underline{\mathbf{W}}\|_{X_\infty}+
+\sup_{t>0}(1+t)^{\frac{5}{4}}\|\delta_+\underline{\mathbf{W}}\|_{X_2}\\
&+\left(\int_0^{\infty}(1+t)^r\|\delta_+\partial_{xx}\underline{\mathbf{W}}(t)\|^r_{X_2}\right)^{1/r}.
\eld
\eeqs
Here we have $\delta^2:=\delta_-\delta_+$, where $\delta_{\pm}$ is defined in \eqref{e:deltaGamma}.
\bl[\textbf{maximal regularity}]
\label{l:mr}
For any given $T > 0$ and $r \in (1, +\infty)$, there exists a positive constant $C$ such that the following holds.
If $(\underline{\eta},\underline{\mathbf{v}})\in L^r((0, T), \ell^2\times X_2)$ and if $(\underline{\theta},\underline{\mathbf{w}})$ satisfies
\beqs
\begin{pmatrix}\underline{\theta}(t) \\ \underline{\mathbf{w}}(t)\end{pmatrix}
=\int_{0}^{t}\mathscr{M}(t-s)\begin{pmatrix}\underline{\eta}(s)\\ \underline{\mathbf{v}}(s)\end{pmatrix}\rmd s,
\quad t\in [0, T],
\eeqs
then
\beqs
\int_0^T\|\partial_{xx}\underline{\mathbf{w}}(t)\|_{X_2}^r\rmd t \leq C
\int_0^T\left(\|\underline{\eta}(t)\|_{\ell^2}+\|\underline{\mathbf{v}}(t)\|_{X_2}\right)^r\rmd t.
\eeqs
\el
\bpf
The result just follows from the standard maximal regularity results on the Laplacian operator and the robustness of maximal regularity
with respect to lower order perturbations. To see that, we first recall $\mathscr{M}(t)=\rme^{A_{\mathrm{nf}}t}$, where $A_{\mathrm{nf}}$ is defined in \eqref{e:A2}.
By \cite{mielke_1987}, the maximal regularity result holds when we just replace $A_{\mathrm{nf}}$ by $A_0$, which is defined as
\beqs
A_0=\begin{pmatrix}0&0\\0&D\partial_{xx}\end{pmatrix}.
\eeqs
Viewing $A_{\mathrm{nf}}$ as a perturbation of $A_0$, we have
\beqs
\begin{pmatrix}\underline{\theta}(t) \\ \underline{\mathbf{w}}(t)\end{pmatrix}
=\int_{0}^{t}\mathscr{M}(t-s)\begin{pmatrix}\underline{\eta}(s)\\ \underline{\mathbf{v}}(s)\end{pmatrix}\rmd s
=\int_0^t\rme^{A_0(t-s)}\left((A_{\mathrm{nf}}-A_0)\begin{pmatrix}\underline{\theta}(s) \\ \underline{\mathbf{w}}(s)\end{pmatrix}+
\begin{pmatrix}\underline{\eta}(s)\\ \underline{\mathbf{v}}(s)\end{pmatrix}\right)\rmd s.
\eeqs
Then by the maximal regularity property of $A_0$, we obtain
\beqs
\int_0^T\|\partial_{xx}\underline{\mathbf{w}}(t)\|_{X_2}^r\rmd t \leq C
\int_0^T\left(\|(A_{\mathrm{nf}}-A_0)\begin{pmatrix}\underline{\theta}(s)\\ \underline{\mathbf{w}}(s)\end{pmatrix}\|_{\ell^2\times X_2}+
\|\begin{pmatrix}\underline{\eta}(s)\\ \underline{\mathbf{v}}(s)\end{pmatrix}\|_{\ell^2\times X_2}\right)^r\rmd s	.
\eeqs
We observe that, for any $\epsilon>0$, there exists $K(\epsilon)>0$ such that
\beqs
\|(A_{\mathrm{nf}}-A_0)\begin{pmatrix}\underline{\theta} \\ \underline{\mathbf{w}}\end{pmatrix}\|_{\ell^2\times X_2}\leq
\epsilon \|A_0\begin{pmatrix}\underline{\theta} \\ \underline{\mathbf{w}}\end{pmatrix}\|_{\ell^2\times X_2}+
K(\epsilon)\|\begin{pmatrix}\underline{\theta} \\ \underline{\mathbf{w}}\end{pmatrix}\|_{\ell^2\times X_2}.
\eeqs
In addition, it is straightforward to see that
\beqs
\int_0^T \|\begin{pmatrix}\underline{\theta}(t) \\ \underline{\mathbf{w}}(t)\end{pmatrix}\|_{\ell^2\times X_2}^r\rmd t\leq C
\int_0^T \|\begin{pmatrix}\underline{\eta}(t)\\ \underline{\mathbf{v}}(t)\end{pmatrix}\|_{\ell^2\times X_2}^r\rmd t.
\eeqs
The conclusion follows by combing the above three inequalities and taking $\epsilon$ sufficiently small.
\epf

We also prove a corollary which will be useful in the proof of nonlinear stability.

\begin{corollary}
\label{c:mr}
  For given $\alpha\in\R$ and $r\in(1,\infty)$, there exists a positive constant $C$ such that, if
\beqs
\begin{pmatrix}\underline{\theta}(t) \\ \underline{\mathbf{w}}(t)\end{pmatrix}
=\int_{t-1}^{t}\mathscr{M}(t-s)\begin{pmatrix}\underline{\eta}(s)\\ \underline{\mathbf{v}}(s)\end{pmatrix}\rmd s,
\quad t\geq 1,
\eeqs
then
\beqs
\int_1^\infty(1+t)^\alpha\|\partial_{xx}\underline{\mathbf{w}}(t)\|_{X_2}^r\rmd t \leq C
\int_0^\infty(1+t)^\alpha\left(\|\underline{\eta}(t)\|_{\ell^2}+\|\underline{\mathbf{v}}(t)\|_{X_2}\right)^r\rmd t.
\eeqs
\end{corollary}
\bpf
We first note that, for $t\in[n,n+1)$, $n\in\N\backslash\{0\}$,
\beqs
\bld
\begin{pmatrix}\underline{\theta}(t) \\ \underline{\mathbf{w}}(t)\end{pmatrix}&
=\left(\int_{n-1}^{t}-\int_{n-1}^{t-1}\right)\mathscr{M}(t-s)\begin{pmatrix}\underline{\eta}(s)\\ \underline{\mathbf{v}}(s)\end{pmatrix}\rmd s\\
&=\left(\int_{0}^{t-n+1}-\int_0^{t-n}\right)\mathscr{M}(t-n+1-s)
\begin{pmatrix}\underline{\eta}(n-1+s)\\ \underline{\mathbf{v}}(n-1+s)\end{pmatrix}\rmd s.
\eld
\eeqs
Applying Lemma \ref{l:mr} to the above expression, we obtain
\beqs
\int_n^{n+1}\|\partial_{xx}\underline{\mathbf{w}}(t)\|_{X_2}^r\rmd t \leq C
\int_{n-1}^{n+1}\left(\|\underline{\eta}(t)\|_{\ell^2}+\|\underline{\mathbf{v}}(t)\|_{X_2}\right)^r\rmd t.
\eeqs
The conclusion follows from multiplying both sides with $n^\alpha\sim (1+t)^\alpha$ and summing over $n\in\N\backslash\{0\}$.
\epf

\bl
\label{l:fpl}
If $(\underline{\theta}_0,\underline{\mathbf{W}}_0)\in\ell^1\times (X_{\mathrm{ch}}^{\perp}\cap \mathscr{T}_{\mathrm{ch}}^{-1}(H^2))$,
the solution of the linear system
\beqs
\begin{pmatrix} \underline{\theta}(t)\\ \underline{\mathbf{W}}(t) \end{pmatrix}
=\mathscr{M}(t)\begin{pmatrix} \underline{\theta}_0 \\ \underline{\mathbf{W}}_0 \end{pmatrix}
\eeqs
belongs to $Z$ and there exists a positive constant $C_1>0$ such that
\beq
\label{e:smallterm}
\|\begin{pmatrix} \underline{\theta}(t) \\ \underline{\mathbf{W}}(t) \end{pmatrix}\|_Z\leq C_1
\|\begin{pmatrix} \underline{\theta}_0 \\ \underline{\mathbf{W}}_0
\end{pmatrix}\|_{\ell^1\times(X_{\mathrm{ch}}\cap \mathscr{T}_{\mathrm{ch}}^{-1}(H^2))}.
\eeq
\el
\bpf
By Proposition \ref{pro:41}, it is straightforward to see that
\beqs
\bld
\|\mathscr{M}_{00}(t)\underline{\theta}_0\|_{\ell^1}&\leq C\|\underline{\theta}_0\|_{\ell^1},&
\|\mathscr{M}_{01}(t)\underline{\mathbf{W}}_0\|_{\ell^1}&\leq \frac{C}{(1+t)^{1/2}}\|\underline{\mathbf{W}}_0\|_{X_1},\\
\|\mathscr{M}_{00}(t)\underline{\theta}_0\|_{\ell^\infty}&\leq \frac{C}{(1+t)^{1/2}}\|\underline{\theta}_0\|_{\ell^1},&
\|\mathscr{M}_{01}(t)\underline{\mathbf{W}}_0\|_{\ell^\infty}&\leq \frac{C}{1+t}\|\underline{\mathbf{W}}_0\|_{X_1},\\
\|\delta^2\mathscr{M}_{00}(t)\underline{\theta}_0\|_{\ell^2}&\leq \frac{C}{(1+t)^{5/4}}\|\underline{\theta}_0\|_{\ell^1},&
\|\delta^2\mathscr{M}_{01}(t)\underline{\mathbf{W}}_0\|_{\ell^2}&\leq \frac{C}{(1+t)^{7/4}}\|\underline{\mathbf{W}}_0\|_{X_1},\\
\|\mathscr{M}_{10}(t)\underline{\theta}_0\|_{X_1}&\leq \frac{C}{(1+t)^{1/2}}\|\underline{\theta}_0\|_{\ell^1},&
\|\mathscr{M}_{11}(t)\underline{\mathbf{W}}_0\|_{X_1}&\leq \frac{C}{1+t}\|\underline{\mathbf{W}}_0\|_{X_1},\\
\|\mathscr{M}_{10}(t)\underline{\theta}_0\|_{X_\infty}&\leq \frac{C}{1+t}\|\underline{\theta}_0\|_{\ell^1},&
\|\mathscr{M}_{11}(t)\underline{\mathbf{W}}_0\|_{X_\infty}&\leq \frac{C}{1+t}\|\underline{\mathbf{W}}_0\|_{X_\infty}.\\
\eld
\eeqs
Moreover, we have
\beqs
\bld
\|\delta_+\partial_{xx}\underline{\mathbf{W}}(t)\|_{X_2}&\leq C
\|\delta_+A_0\mathscr{M}(t)
\begin{pmatrix}\underline{\theta}_0 \\ \underline{\mathbf{W}}_0\end{pmatrix}\|_{\ell^2\times X_2}\\
&\leq C(\|\delta_+A_{\mathrm{nf}}\mathscr{M}(t)\begin{pmatrix}\underline{\theta}_0 \\ \underline{\mathbf{W}}_0\end{pmatrix}\|_{\ell^2\times X_2}
+\|\delta_+(A_{\mathrm{nf}}-A_0)\mathscr{M}(t)\begin{pmatrix}\underline{\theta}_0 \\ \underline{\mathbf{W}}_0\end{pmatrix}\|_{\ell^2\times X_2}).
\eld
\eeqs
We need to show that the two terms on the right hand side of the above inequality decay sufficiently fast.
On the one hand, we claim that
\beqs
\|\delta_+A_{\mathrm{nf}}\mathscr{M}(t)\begin{pmatrix}\underline{\theta}_0 \\ \underline{\mathbf{W}}_0\end{pmatrix}\|_{\ell^2\times X_2}\leq
\frac{C}{(1+t)^{\frac{3}{2}}}\|\begin{pmatrix}\underline{\theta}_0\\ \underline{\mathbf{W}}
\end{pmatrix}\|_{\ell^2\times \mathscr{T}_{\mathrm{ch}}^{-1}(H^2)},
\text{ for all } t\geq 0.
\eeqs
Actually, for $t\in [0,1]$, the above inequality is true since $\delta_+$ is bounded and
\beqs
\|A_{\mathrm{nf}}\mathscr{M}(t)\begin{pmatrix}\underline{\theta}_0 \\ \underline{\mathbf{W}}_0\end{pmatrix}\|_{\ell^2\times X_2}\leq C
\|A_{\mathrm{nf}}\begin{pmatrix}\underline{\theta}_0 \\ \underline{\mathbf{W}}_0\end{pmatrix}\|_{\ell^2\times X_2}\leq C
\|\begin{pmatrix}\underline{\theta}_0\\ \underline{\mathbf{W}}\end{pmatrix}\|_{\ell^2\times \mathscr{T}_{\mathrm{ch}}^{-1}(H^2)}.
\eeqs
For $t\in[1,\infty]$, we first point out that, to show
$\opnorm{A_{\mathrm{nf}}\mathscr{M}(t)=A_{\mathrm{nf}}\rme^{A_{\mathrm{nf}}t}}_{\ell^2\times X_2}$
decays with rate $t^{-1}$ as $t$ goes to $\infty$,
we only have to show that the supremum norm of its Fourier-Bloch counterpart
$\widehat{A}_{\mathrm{nf}}(\sigma)M(t,\sigma)$ decays with rate $t^{-1}$ as $t$ goes to $\infty$,
just as in the scalar heat equation case.	
This is true by applying the steps in Lemma \ref{l:pen0} and Lemma \ref{l:pea0}
to $\widehat{A}_{\mathrm{nf}}(\sigma)M(t,\sigma)$. Second, it is straightforward to see that
the discrete derivative operator $\delta_+$ gives an extra $t^{-1/2}$ decay, which concludes our justification.
On the other hand, we have the explicit expression, using that $\delta_+$ and $A_{\mathrm{nf}}$, $A_0$ commute,
\beqs
\delta_+(A_{\mathrm{nf}}-A_0)\mathscr{M}(t)\begin{pmatrix}\underline{\theta}_0 \\ \underline{\mathbf{W}}_0\end{pmatrix}=
\begin{pmatrix}0&\delta_+\Gamma\\A_{\mathrm{ch}}\mathbf{E}*&\mathbf{f}^\prime(\mathbf{u}_\star)-\mathbf{E}*\delta_+\Gamma\end{pmatrix}
\begin{pmatrix}\delta_+\underline{\theta}(t)\\ \delta_+\underline{\mathbf{W}}(t)\end{pmatrix}.
\eeqs
We apply Proposition \ref{pro:41} again and obtain
\beqs
\bld
\|A_{\mathrm{ch}}\mathbf{E}*(\delta_+\underline{\theta}(t))\|_{\ell^2}&\leq C\|\delta^2\underline{\theta}(t)\|_{\ell^2}
\leq \frac{C}{(1+t)^{\frac{5}{4}}}\left(\|\underline{\theta}_0\|_{\ell^1}+\|\underline{\mathbf{W}}_0\|_{X_2}\right),\\
\|\mathbf{f}^\prime(\mathbf{u}_\star)(\delta_+\underline{\mathbf{W}}(t))\|_{X_2}&\leq C\|\delta_+\underline{\mathbf{W}}(t)\|_{X_2}
\leq \frac{C}{(1+t)^{\frac{5}{4}}}(\|\underline{\theta}_0\|_{\ell^1}+\|\underline{\mathbf{W}}_0\|_{X_2}).
\eld
\eeqs
In addition, recalling that $\Gamma$ is defined in \eqref{e:deltaGamma}, we conclude that,
for any $\epsilon>0$, there exists $K(\epsilon)>0$ such that
\beqs
\bld
\|\delta_+\Gamma(\delta_+\underline{\mathbf{W}}(t))\|_{X_2}&\leq C\|\partial_x\delta_+\underline{\mathbf{W}}(t)\|_{X_2}
\leq \epsilon\|\delta_+\partial_{xx}\underline{\mathbf{W}}(t)\|_{X_2}+K(\epsilon)\|\delta_+\underline{\mathbf{W}}(t)\|_{X_2}.
\eld
\eeqs
Therefore, by choosing $\epsilon$ sufficiently small, we conclude that
\beqs
\bld
\|\delta_+\partial_{xx}\underline{\mathbf{W}}(t)\|_{X_2}&\leq
\frac{C}{(1+t)^{\frac{5}{4}}}\left(\|\underline{\theta}_0\|_{\ell^1}+\|\underline{\mathbf{W}}_0\|_{ \mathscr{T}_{\mathrm{ch}}^{-1}(H^2)}\right),
\eld
\eeqs
which shows that
\beqs
\left(\int_0^{\infty}(1+t)^r\|\delta_+\partial_{xx}\underline{\mathbf{W}}(t)\|_{X_2}^r\rmd t\right)^{1/r}\leq
C\left(\|\underline{\theta}_0\|_{\ell^1}+\|\underline{\mathbf{W}}_0\|_{\mathscr{T}_{\mathrm{ch}}^{-1}(H^2)}\right).
\eeqs
This proves the lemma.
\epf

\bl
\label{l:fpn}
For $\|(\underline{\theta}(t),\underline{\mathbf{W}}(t))\|_Z<\varepsilon$,
where $\varepsilon$ is sufficiently small ($0<\varepsilon\leq\varepsilon_0$), there exists a positive constant $C_2\geq 1$ such that
\beq
\label{e:IE}
\|\int_0^t\mathscr{M}(t-s)\begin{pmatrix} \mathbf{N}^\theta(\underline{\theta}(s),\underline{\mathbf{W}}(s))
\\ \mathbf{N}^{\mathbf{w}}(\underline{\theta}(s),\underline{\mathbf{W}}(s)) \end{pmatrix}
\rmd s\|_Z\leq C_2\|(\underline{\theta}(t), \underline{\mathbf{W}}(t))\|_Z^2 .
\eeq
Moreover, for $(\underline{\theta}_1,\underline{\mathbf{W}}_1)$, $(\underline{\theta}_2, \underline{\mathbf{W}}_2)$
with their norms in $Z$ smaller than $\varepsilon$, we have
\beq
\label{e:IEL}
\bld
&\|\int_0^t\mathscr{M}(t-s)\begin{pmatrix}
\mathbf{N}^\theta(\underline{\theta}_1(s),\underline{\mathbf{W}}_1(s))-\mathbf{N}^\theta(\underline{\theta}_2(s),\underline{\mathbf{W}}_2(s))
\\ \mathbf{N}^{\mathbf{w}}(\underline{\theta}_1(s),\underline{\mathbf{W}}_1(s))
-\mathbf{N}^{\mathbf{w}}(\underline{\theta}_2(s),\underline{\mathbf{W}}_2(s)) \end{pmatrix}
\rmd s\|_Z\\
\leq & C_2\left(\sum_{j=1}^2\|(\underline{\theta}_j(t), \underline{\mathbf{W}}_j(t))\|_Z\right)
\|(\underline{\theta}_1(t)-\underline{\theta}_2(t), \underline{\mathbf{W}}_1(t)-\underline{\mathbf{W}}_2(t))\|_Z .
\eld
\eeq
\el
\bpf
We start with proving the estimate \eqref{e:IE}.
The proof is fairly straightforward. The strategy is to use estimates for the linear part $\mathscr{M}(t)$ in Proposition \ref{pro:41},
the estimates for the nonlinear terms in Lemma \ref{l:61e} from the appendix, and the maximal regularity estimates in Lemma \ref{l:mr},
Corollary \ref{c:mr}. For simplicity, we denote
\beqs
\mathbf{N}^\theta(s)=\mathbf{N}^\theta(\underline{\theta}(s),\underline{\mathbf{W}}(s)),\quad
\mathbf{N}^{\mathbf{w}}(s)=\mathbf{N}^{\mathbf{w}}(\underline{\theta}(s),\underline{\mathbf{W}}(s)).
\eeqs

By Lemma \ref{l:61e}, we have that
\beq
\label{e:NN}
\bld
\|\mathbf{N}^\theta(s)\|_{\ell^1}&\leq \frac{C}{(1+s)^{\frac{3}{2}}}\|(\underline{\theta}(t), \underline{\mathbf{W}}(t))\|_Z^2+
\frac{C}{(1+s)^{\frac{5}{4}}}\|(\underline{\theta}(t), \underline{\mathbf{W}}(t))\|_Z(1+s)\|\delta_+\partial_{xx}\underline{\mathbf{W}}(s)\|_{X_2},\\
&\\
\|\mathbf{N}^\theta(s)\|_{\ell^2}&\leq \frac{C}{(1+s)^{\frac{3}{2}}}\|(\underline{\theta}(t), \underline{\mathbf{W}}(t))\|_Z^2+
\frac{C}{(1+s)^{\frac{3}{2}}}\|(\underline{\theta}(t), \underline{\mathbf{W}}(t))\|_Z(1+s)\|\delta_+\partial_{xx}\underline{\mathbf{W}}(s)\|_{X_2},\\
&\\
\|\mathbf{N}^{\mathbf{w}}(s)\|_{X_1}&\leq \frac{C}{1+s}\|(\underline{\theta}(t), \underline{\mathbf{W}}(t))\|_Z^2+
\frac{C}{(1+s)^{\frac{5}{4}}}\|(\underline{\theta}(t), \underline{\mathbf{W}}(t))\|_Z(1+s)\|\delta_+\partial_{xx}\underline{\mathbf{W}}(s)\|_{X_2},\\
&\\
\|\mathbf{N}^{\mathbf{w}}(s)\|_{X_2}&\leq \frac{C}{(1+s)^{\frac{5}{4}}}\|(\underline{\theta}(t), \underline{\mathbf{W}}(t))\|_Z^2+
\frac{C}{(1+s)^{\frac{3}{2}}}\|(\underline{\theta}(t), \underline{\mathbf{W}}(t))\|_Z(1+s)\|\delta_+\partial_{xx}\underline{\mathbf{W}}(s)\|_{X_2}.
\eld
\eeq
We also exploit the linear estimates from Proposition \ref{pro:41} and obtain the following estimates.
\beq
\label{e:IE1}
\bld
N_1&=\|\int_0^t\mathscr{M}_{00}(t-s)\mathbf{N}^\theta(s)\rmd s\|_{\ell^1}\leq  C\int_0^t\|\mathbf{N}^\theta(s)\|_{\ell^1}\rmd s,\\
&\\
N_2&=\|\int_0^t\mathscr{M}_{01}(t-s)\mathbf{N}^{\mathbf{w}}(s)\rmd s\|_{\ell^1}\leq
C\int_0^t\frac{\|\mathbf{N}^{\mathbf{w}}(s)\|_{X_1}}{(1+t-s)^{\frac{1}{2}}}\rmd s,\\
&\\
N_3&=(1+t)^{\frac{1}{2}}\|\int_0^t\mathscr{M}_{00}(t-s)\mathbf{N}^\theta(s)\rmd s\|_{\ell^\infty}\leq
C(1+t)^{\frac{1}{2}}\int_0^t\frac{\|\mathbf{N}^\theta(s)\|_{\ell^1}}{(1+t-s)^{\frac{1}{2}}}\rmd s,\\
&\\
N_4&=(1+t)^{\frac{1}{2}}\|\int_0^t\mathscr{M}_{01}(t-s)\mathbf{N}^{\mathbf{w}}(s)\rmd s\|_{\ell^\infty}\leq
C(1+t)^{\frac{1}{2}}\int_0^t\frac{\|\mathbf{N}^{\mathbf{w}}(s)\|_{X_1}}{1+t-s}\rmd s,\\
&\\
N_5&=(1+t)^{\frac{5}{4}}\|\delta^2\int_0^t\mathscr{M}_{00}(t-s)\mathbf{N}^\theta(s)\rmd s\|_{\ell^2}\leq
C(1+t)^{\frac{5}{4}}\int_0^t\frac{\|\mathbf{N}^\theta(s)\|_{\ell^1}}{(1+t-s)^{\frac{5}{4}}}\rmd s,\\
&\\
N_6&=(1+t)^{\frac{5}{4}}\|\delta^2\int_0^t\mathscr{M}_{01}(t-s)\mathbf{N}^{\mathbf{w}}(s)\rmd s\|_{\ell^2}\leq
C(1+t)^{\frac{5}{4}}\int_0^t\frac{\|\mathbf{N}^{\mathbf{w}}(s)\|_{X_2}}{(1+t-s)^{\frac{3}{2}}}\rmd s,\\
&\\
N_7&=(1+t)^{\frac{1}{2}}\|\int_0^t\mathscr{M}_{10}(t-s)\mathbf{N}^\theta(s)\rmd s\|_{X_1} \leq
C(1+t)^{\frac{1}{2}}\int_0^t\frac{\|\mathbf{N}^\theta(s)\|_{\ell^1}}{(1+t-s)^{\frac{1}{2}}}\rmd s,\\
&\\
N_8&=(1+t)^{\frac{1}{2}}\|\int_0^t\mathscr{M}_{11}(t-s)\mathbf{N}^{\mathbf{w}}(s)\rmd s\|_{X_1} \leq
C(1+t)^{\frac{1}{2}}\int_0^t\frac{\|\mathbf{N}^{\mathbf{w}}(s)\|_{X_1}}{1+t-s}\rmd s,\\
&\\
N_9&=(1+t)^{\frac{5}{4}}\|\delta_+\int_0^t\mathscr{M}_{10}(t-s)\mathbf{N}^\theta(s)\rmd s\|_{X_2} \leq
C(1+t)^{\frac{5}{4}}\int_0^t\frac{\|\mathbf{N}^\theta(s)\|_{\ell^1}}{(1+t-s)^{\frac{5}{4}}}\rmd s,\\
&\\
N_{10}&=(1+t)^{\frac{5}{4}}\|\delta_+\int_0^t\mathscr{M}_{11}(t-s)\mathbf{N}^{\mathbf{w}}(s)\rmd s\|_{X_2}\leq
C(1+t)^{\frac{5}{4}}\int_0^t\frac{\|\mathbf{N}^{\mathbf{w}}(s)\|_{X_2}}{(1+t-s)^{\frac{3}{2}}}\rmd s,\\
&\\
N_{11}&=(1+t)\|\int_0^t\mathscr{M}_{10}(t-s)\mathbf{N}^\theta(s)\rmd s\|_{X_\infty} \leq
C(1+t)\int_0^t\frac{\|\mathbf{N}^\theta(s)\|_{\ell^1}}{1+t-s}\rmd s,\\
&\\
N_{12}&=(1+t)\|\int_0^t\mathscr{M}_{11}(t-s)\mathbf{N}^{\mathbf{w}}(s)\rmd s\|_{X_\infty} \leq
C(1+t)\int_0^t\frac{\|\mathbf{N}^{\mathbf{w}}(s)\|_{X_1}}{(1+t-s)(t-s)^{\frac{1}{2}}}\rmd s.
\eld
\eeq
At this point, we substitute \eqref{e:NN} into \eqref{e:IE1}, estimate the resulting integrals, and find
\beq
\label{e:IE2}
N_j\leq C\|(\underline{\theta}(t), \underline{\mathbf{W}}(t))\|_Z^2, \text{ for all }1\leq j\leq 12.
\eeq
The calculations establishing the estimates for $N_1$, ..., $N_{11}$ are based on the following elementary integral estimates.
\beqs
\int_0^t\frac{1}{(1+t-s)^{\alpha}}\frac{1}{(1+s)^{\beta}}\rmd s\leq
\frac{C}{(1+t)^\alpha}\int_0^{\frac{t}{2}}\frac{1}{(1+s)^{\beta}}\rmd s+
\frac{C}{(1+t)^\beta}\int_0^{\frac{t}{2}}\frac{1}{(1+s)^{\alpha}}\rmd s.
\eeqs
For the estimate on $N_{12}$, we just need to show that the following integral expression
\beqs
h(t)=(1+t)\left(\int_0^t\frac{1}{(1+t-s)(t-s)^{\frac{1}{2}}}\frac{1}{1+s}\rmd s+
\left(\int_0^t\left(\frac{1}{(1+t-s)(t-s)^{\frac{1}{2}}}\frac{1}{(1+s)^{\frac{5}{4}}}\right)^{\frac{r}{r-1}}\rmd s\right)^{1-\frac{1}{r}}\right)
\eeqs
has a uniform upper bound for $t\in(0,\infty)$.
First, for all $t\in(0, 1]$, there exists $C>0$ such that,
\beqs
h(t)\leq 2\left(\int_0^1(t-s)^{-\frac{1}{2}}\rmd s+\left(\int_0^1(t-s)^{-\frac{r}{2(r-1)}}\rmd s\right)^{1-\frac{1}{r}}\right) \leq C
\eeqs
Second, for $t\in[1,\infty)$, we have
\beqs
\bld
(1+t)\int_0^t\frac{1}{(1+t-s)(t-s)^{\frac{1}{2}}}\frac{1}{1+s}\rmd s
&\leq \frac{C}{(1+t)^{\frac{1}{2}}}\int_0^{\frac{t}{2}}\frac{1}{1+s}\rmd s+C\int_{\frac{t}{2}}^t\frac{1}{(1+t-s)(t-s)^{\frac{1}{2}}}\rmd s\\
&\leq C\left(1+\int_0^\infty\frac{1}{(1+s)s^{\frac{1}{2}}}\rmd s\right)\leq C.
\eld
\eeqs
Similar arguments show that the second part of $h(t)$ is also uniformly bounded on $[1,\infty)$.

The estimates on $N_1$, ..., $N_{12}$ bound the $Z$-norm
of the left-hand side of \eqref{e:IE}, except for the maximal regularity component.
Thus it remains to show that
\beq
\label{e:NM}
\left(\int_0^\infty(1+t)^r\|\delta_+\partial_{xx}\mathscr{W}(t)\|_{X_2}^r\rmd t\right)^{\frac{1}{r}}
\leq C\|(\underline{\theta}(t), \underline{\mathbf{W}}(t))\|_Z^2,
\eeq
where
\beqs
\mathscr{W}(t)=\int_0^t\mathscr{M}_{10}(t-s)\mathbf{N}^\theta(s)+\mathscr{M}_{11}(t-s)\mathbf{N}^{\mathbf{w}}(s)\rmd s.
\eeqs
For $t\in[0,1]$, by maximal regularity in Lemma \ref{l:mr}, we have
\beq
\label{e:NM1}
\bld
\int_0^1(1+t)^r\|\delta_+\partial_{xx}\mathscr{W}(t)\|_{X_2}^r\rmd t\leq & C\int_0^1\|\partial_{xx}\mathscr{W}(t)\|_{X_2}^r\rmd t\\
\leq & C\int_0^1\left(\|\mathbf{N}^\theta(t)\|_{\ell^2}+\|\mathbf{N}^{\mathbf{w}}(t)\|_{X_2}\right)^r\rmd t\\
\leq &C\|(\underline{\theta}(t), \underline{\mathbf{W}}(t))\|_Z^{2r}.
\eld
\eeq
For $t\in[1,\infty)$, we split $\mathscr{W}$ into two parts, that is,
\beqs
\mathscr{W}=\left(\int_0^{t-1}+\int_{t-1}^{t}\right)\mathscr{M}_{10}(t-s)\mathbf{N}^\theta(s)+\mathscr{M}_{11}(t-s)\mathbf{N}^{\mathbf{w}}(s)\rmd s
=\mathscr{W}_1+\mathscr{W}_2.	
\eeqs
By Corollary \ref{c:mr}, we have
\beqs
\bld
\int_1^\infty(1+t)^r\|\delta_+\partial_{xx}\mathscr{W}_2(t)\|_{X_2}^r\rmd t\leq &
C\int_1^\infty(1+t)^r\|\partial_{xx}\mathscr{W}_2(t)\|_{X_2}^r\rmd t\\
\leq & C\int_0^\infty(1+t)^r\left(\|\mathbf{N}^\theta(t)\|_{\ell^2}+\|\mathbf{N}^{\mathbf{w}}(t)\|_{X_2}\right)^r\rmd t\\
\leq & C\|(\underline{\theta}(t), \underline{\mathbf{W}}(t))\|_Z^{2r}.
\eld
\eeqs
By similar arguments as in Lemma \ref{l:fpl} and with the condition that $t-s>1$, we can show that
\beqs
\bld
&\|\delta_+\partial_{xx}\left(\mathscr{M}_{10}(t-s)\mathbf{N}^\theta(s)+\mathscr{M}_{11}(t-s)\mathbf{N}^{\mathbf{w}}(s)\right)\|_{X_2}\leq
\frac{C}{(1+t-s)^{\frac{5}{4}}}\left(\|\mathbf{N}^\theta(s)\|_{\ell^1}+\|\mathbf{N}^{\mathbf{w}}(s)\|_{X_2}\right)\\ \leq &
\frac{C}{(1+t-s)^{\frac{5}{4}}(1+s)^{\frac{5}{4}}}\|(\underline{\theta}(t), \underline{\mathbf{W}}(t))\|_Z
\bigg(\|(\underline{\theta}(t), \underline{\mathbf{W}}(t))\|_Z
+(1+s)\|\delta_+\partial_{xx}\underline{\mathbf{W}}(s)\|_{X_2}\bigg).
\eld
\eeqs
As a result, we obtain
\beqs
\bld
&(1+t)\|\delta_+\partial_{xx}\mathscr{W}_1(t)\|_{X_2}\leq (1+t)\int_0^{t-1}
\|\delta_+\partial_{xx}\left(\mathscr{M}_{10}(t-s)\mathbf{N}^\theta(s)+\mathscr{M}_{11}(t-s)\mathbf{N}^{\mathbf{w}}(s)\right)\|_{X_2}\rmd s\\
\leq &(1+t)\int_0^{t-1}\frac{C}{(1+t-s)^{\frac{5}{4}}(1+s)^{\frac{5}{4}}}\|(\underline{\theta}(t), \underline{\mathbf{W}}(t))\|_Z
\left(\|(\underline{\theta}(t), \underline{\mathbf{W}}(t))\|_Z
+(1+s)\|\delta_+\partial_{xx}\underline{\mathbf{W}}(s)\|_{X_2}\right)\rmd s\\
\leq &\frac{C}{(1+t)^{\frac{1}{4}}}\|(\underline{\theta}(t), \underline{\mathbf{W}}(t))\|_Y^{2r},
\eld
\eeqs
which immediately implies that
\beqs
\int_1^\infty(1+t)^r\|\delta_+\partial_{xx}\mathscr{W}_1(t)\|_{X_2}^r
\rmd t\leq C\|\left(\underline{\theta}(t), \underline{\mathbf{W}}(t)\right)\|_Z^{2r}.
\eeqs
Together with \eqref{e:NM1}, this establishes \eqref{e:NM} and concludes the proof.

In a completely analogous fashion, one establishes the Lipshitz estimates.
\epf
We now prove our main theorem.
\bpf[Proof of Theorem \ref{t:1}]
The proof is a fixed-point-theorem argument. We first recall the variation of constant formula,
\beqs
\begin{pmatrix}\underline{\theta}(t)\\ \underline{\mathbf{W}}(t)\end{pmatrix}=
\mathscr{t}\begin{pmatrix}\underline{\theta}_0\\ \underline{\mathbf{W}}_0\end{pmatrix}+
\int_0^t\mathscr{M}(t-s)\begin{pmatrix}\mathbf{N}^\theta(\underline{\theta}(s),\underline{\mathbf{W}}(s))
\\ \mathbf{N}^{\mathbf{w}}(\underline{\theta}(s),\underline{\mathbf{W}}(s))\end{pmatrix}\rmd s.
\eeqs
Let $\mathscr{P}$ be the right-hand side of the formula, that is
\beqs
\mathscr{P}\begin{pmatrix}\underline{\theta}(t)\\ \underline{\mathbf{W}}(t)\end{pmatrix}=
\mathscr{t}\begin{pmatrix}\underline{\theta}_0\\ \underline{\mathbf{W}}_0\end{pmatrix}+
\int_0^t\mathscr{M}(t-s)\begin{pmatrix}
\mathbf{N}^\theta(\underline{\theta}(s),\underline{\mathbf{W}}(s))
\\ \mathbf{N}^{\mathbf{w}}(\underline{\theta}(s),\underline{\mathbf{W}}(s))\end{pmatrix}\rmd s.
\eeqs
Assume that now the initial value is sufficiently small, that is, for some small $\epsilon>0$,
\beqs
\|(\underline{\theta}_0, \underline{\mathbf{W}}_0)\|_{\ell^1\times(Z\cap \mathscr{T}^{-1}(H^2))}\leq \epsilon.
\eeqs
If $(\underline{\theta}(t),\underline{\mathbf{W}}(t))\in Z$ with norm smaller than $\varepsilon$, we know that
\beqs
\mathscr{P}\begin{pmatrix}\underline{\theta}(t)\\ \underline{\mathbf{W}}(t)\end{pmatrix}\in Z.
\eeqs
By Lemma \ref{l:fpl} and \ref{l:fpn}, we have that
\beq
\label{e:fp1}
\|\mathscr{P}\begin{pmatrix}\underline{\theta}(t)\\ \underline{\mathbf{W}}(t)\end{pmatrix}\|_Z
\leq C_1\epsilon + C_2\|\begin{pmatrix}\underline{\theta}(t)\\ \underline{\mathbf{W}}(t)\end{pmatrix}\|_Z^2
\eeq
Moreover, we have
\beq
\label{e:fp2}
\|\mathscr{P}\begin{pmatrix}\underline{\theta}_1(t)\\ \underline{\mathbf{W}}_1(t)\end{pmatrix}-
\mathscr{P}\begin{pmatrix}\underline{\theta}_2(t)\\ \underline{\mathbf{W}}_2(t)\end{pmatrix}
\|_Z\leq C_2\left(\sum_{j=1}^2\|(\underline{\theta}_j(t), \underline{\mathbf{W}}_j(t))\|_Z\right)
\|(\underline{\theta}_1(t)-\underline{\theta}_2(t), \underline{\mathbf{W}}_1(t)-\underline{\mathbf{W}}_2(t))\|_Z .
\eeq
We denote $B=\{\underline{\mathbf{V}}\in Z\mid \|\underline{\mathbf{V}}\|_Z\leq R\}$,
where $R=\min(2C_1\epsilon , \varepsilon )$. We now take $\epsilon> 0$ small enough so that $2C_2 R < 1$ and readily
conclude, based on \eqref{e:fp1} and \eqref{e:fp2}, that $\mathscr{P}(B) \subset B$ and
that $\mathscr{P}$ is a strict contraction in $B$. By Banach's fixed point theorem, there is a
unique fixed point of $\mathscr{P}$ in $B$, denoted as $(\underline{\theta}(t),\underline{\mathbf{W}}(t))$.
Then $(\underline{\theta}(t),\underline{\mathbf{W}}(t))$ is a global solution of \eqref{e:A2},
and if we return to the original variables, we obtain a global solution of \eqref{e:21} which satisfies the decay
estimate in Theorem \ref{t:1}.
This concludes the proof.
\epf

\section{Appendix}

\subsection{Estimates on nonlinear terms}\label{ss:61}
In this section, we derive the estimates on the nonlinear terms
$\mathbf{N}^\theta$ and $\mathbf{N}^\mathbf{w}$ in our normal form (\ref{e:34}).
\bl
\label{l:61}
For $\|\underline{\mathbf{W}}\|_{X_{\mathrm{ch}}}, \|\underline{\theta}\|_{\ell^1}<\varepsilon$,
where $\varepsilon$ is sufficiently small($0<\varepsilon\leq\varepsilon_0$), there exists a nondecreasing
function $C(\varepsilon)>0$ such that, for all $1\leq p\leq \infty$, the nonlinear terms in system (\ref{e:34}) have the following estimates.
\beq
\label{e:EN}
\bld
|\mathbf{N}^\theta_j|\leq &C(\varepsilon)\bigg[\sum_{k=-1}^{0}|(\delta_+\underline{\theta})_{j+k}|^2+\left(\sum_{k=-1}^{1}|\theta_{j+k}|^3\right)
\left(\sum_{k=-1}^{0}|(\delta_+\underline{\theta})_{j+k}|\right)\\
&+\left(\sum_{k=-1}^{1}|\theta_{j+k}|\right)
\bigg(|(\delta_+\underline{\mathbf{W}})_j(-\pi)|+|(\delta_+\partial_x\underline{\mathbf{W}})_j(-\pi)| \bigg)\\
&+\|\mathbf{W}_j\|_{L^p}|(\delta_+\underline{\mathbf{W}})_j(-\pi)|+\|\mathbf{W}_j^2\|_{L^p}\bigg],\\
\|\mathbf{N}^\mathbf{w}_j\|_{L^p}\leq &
C(\varepsilon)\bigg[\left(\sum_{k=-1}^{0}|(\delta_+\underline{\theta})_{j+k}|\right)\left(\sum_{k=-1}^1|\theta_{j+k}|\right)+
|\theta_j|\|\mathbf{W}_j\|_{L^p}\\
&+\left(\sum_{k=-1}^1|\theta_{j+k}|\right)
\left(\sum_{k=-1}^1|(\delta_+\underline{\mathbf{W}})_{j+k}(-\pi)|+|(\delta_+\partial_x\underline{\mathbf{W}})_{j+k}(-\pi)|\right)\\
&+\|\mathbf{W}_j\|_{L^p}|(\delta_+\underline{\mathbf{W}})_{j}(-\pi)|
+\|\mathbf{W}^2_j\|_{L^p}+
|\mathbf{N}^\theta_j|+|\mathbf{N}^\theta_{j+1}|+|\mathbf{N}^\theta_{j-1}|\bigg].\\
\eld
\eeq
\el

\bpf
We point out that throughout the proof, we repeatedly exploit the fact that
the $L^2$ scalar product of an even function and an odd function are zero.
We also recall that $\mathbf{u}_\star$ is even and $\mathbf{u}_{\mathrm{ad}}$ is odd.
By equations (\ref{e:310}),(\ref{e:PD}) and (\ref{e:34}), we obtain
\beqs
\bld
\mathbf{N}^\theta_j=&\Rmnum{1}_j+\left(\Rmnum{2}_j+\Rmnum{3}_j+\Rmnum{4}_j+\Rmnum{5}_j\right)\mathscr{S}_j \text{ and }\\
\mathbf{N}^\mathbf{w}_j=&\left(\id-\frac{\partial\mathbf{G}_j}{\partial \mathbf{W}_j}\right)^{-1}\left(\Rmnum{6}_j+\Rmnum{7}_j
+\Rmnum{8}_j+\Rmnum{9}_j+\frac{\partial\mathbf{G}_j}{\partial \mathbf{W}_j}\Rmnum{10}_j\right), \text{ where}\\
\mathscr{S}_j=&\left(-1+\langle\mathbf{W}_j(x)+\mathbf{H}_j(x),\mathbf{u}_{\mathrm{ad}}^\prime(x-\theta_j)\rangle\right)^{-1};\\
\mathbf{G}_j=&\mathbf{G}(\theta_j, \mathbf{W}_j)=\langle\mathbf{W}_j(x),\mathbf{u}_{\mathrm{ad}}(x-\theta_j)-\mathbf{u}_{\mathrm{ad}}(x)\rangle\mathbf{\psi}(x-\theta_j);\\
\Rmnum{1}_j=&(-\mathscr{S}_j-1)(\delta_+\Gamma W)_j;\\
\Rmnum{2}_j=&-\left(\mathbf{W}_{j+1}(-\pi)-\mathbf{W}_j(-\pi),D\left(\mathbf{u}_{\mathrm{ad}}^\prime(\pi-\theta_j)-\mathbf{u}_{\mathrm{ad}}^\prime(\pi)\right)\right);\\
\Rmnum{3}_j=&\left(\partial_x\mathbf{W}_{j+1}(-\pi)-\partial_x\mathbf{W}_{j}(-\pi),D\mathbf{u}_{\mathrm{ad}}(\pi-\theta_j)\right);\\
\Rmnum{4}_j=&(\partial_x\mathbf{H}_{j}(\pi)-\partial_x\mathbf{H}_{j}(-\pi),D\mathbf{u}_{\mathrm{ad}}(\pi-\theta_j))-
          (\mathbf{H}_{j}(\pi)-\mathbf{H}_{j}(-\pi),D\mathbf{u}_{\mathrm{ad}}^\prime(\pi-\theta_j));\\
\Rmnum{5}_j=&\langle\tilde{g}(\theta_j,\mathbf{W}_j+\mathbf{H}_j),\mathbf{u}_{\mathrm{ad}}(x-\theta_j)\rangle;\\
\Rmnum{6}_j=&A(\mathbf{H}_j-(\underline{\mathbf{E}}*\underline{\theta})_j);\\
\Rmnum{7}_j=&-\left(\left(\dot{\mathbf{H}}_j-\mathbf{u}^\prime_\star(x-\theta_j)\dot{\theta}_j\right)
-(\underline{\mathbf{E}}*\underline{\dot{\theta}})_j+\langle\dot{\mathbf{W}}_j(x),\mathbf{u}_{\mathrm{ad}}(x-\theta_j)-
\mathbf{u}_{\mathrm{ad}}(x)\rangle\mathbf{\psi}(x-\theta_j)\right);\\
\Rmnum{8}_j=&(\underline{\mathbf{E}}*(\delta_+\Gamma\underline{\mathbf{W}}-\underline{\dot{\theta}}))_j;\\
\Rmnum{9}_j=&\tilde{g}(\theta_j,\mathbf{W}_j+\mathbf{H}_j)+\left[\mathbf{f}^\prime(\mathbf{u}_\star(x-\theta_j))-\mathbf{f}^\prime
(\mathbf{u}_\star(x))\right](\mathbf{W}_j+\mathbf{H}_j);\\
\Rmnum{10}_j=&A(\underline{\mathbf{E}}*\underline{\theta})_j+A\mathbf{W}_j-(\underline{\mathbf{E}}*\delta_+\Gamma\underline{\mathbf{W}})_j.\\
\eld
\eeqs\\
We recall here that $\underline{\mathbf{E}}$ is defined in \eqref{e:E} and point out that
the term in $\Rmnum{7}_j$ involving $\dot{\mathbf{W}}_j$ in fact cancels with a contribution from $\dot{\mathbf{H}}_j$.
We now prove the estimate of $\mathbf{N}^\theta_j$.\\
\textbf{Estimate on $\Rmnum{1}_j$:  $|\Rmnum{1}_j|\leq C(\varepsilon)\bigg(|(\delta_+\underline{\theta})_j|+
|(\delta_-\underline{\theta})_j|+\|\mathbf{W}_j\|_{L^p}\bigg)|(\delta_+\underline{\mathbf{W}})_j(-\pi)|$}.\\
We first recall that $\mathbf{H}_j$ is defined in \eqref{e:PD} and \eqref{e:HD}.
We claim that the number $c_j$, appearing in the definition of $\mathbf{H}_j^2$ as in \eqref{e:H2} and \eqref{e:cj}, can be estimated as
\beqs
|c_j|\leq C(\varepsilon)\left[|(\delta^2\underline{\theta})_j|+\bigg(|(\delta_+\underline{\theta})_{j}|+
|(\delta_-\underline{\theta})_{j}|\bigg)\sum_{k=-1}^1\theta_{j+k}+|\theta_j|\|\mathbf{W}_j\|_{L^p}\right],
\eeqs
where we use notation $\delta^2=\delta_+\delta_-$. In fact, we have
\beqs
\bld
|\langle\phi(x)(\mathbf{u}_\star(x+\theta_j-\theta_{j+1})-\mathbf{u}_\star(x+\theta_j-\theta_{j-1})),\mathbf{u}_{\mathrm{ad}}(x)\rangle|
&\leq C(|(\delta_+\underline{\theta})_{j}|^2+|(\delta_-\underline{\theta})_{j}|^2);\\
|\langle(\phi(x+\theta_j)-\phi(x))(\mathbf{u}_\star(x+\theta_j-\theta_{j+1})-\mathbf{u}_\star(x+\theta_j-\theta_{j-1})),\mathbf{u}_{\mathrm{ad}}(x)\rangle|
&\leq C|\theta_j|(|(\delta_+\underline{\theta})_{j}|+|(\delta_-\underline{\theta})_{j}|);\\
|\langle(\mathbf{u}_\star(x+\theta_j-\theta_{j+1})+\mathbf{u}_\star(x+\theta_j-\theta_{j-1})),\mathbf{u}_{\mathrm{ad}}(x)\rangle|
&\leq C|(\delta^2\underline{\theta})_j|.
\eld
\eeqs
We also have $|\mathbf{H}_j(x)|\leq C(\varepsilon)\bigg(|(\delta_+\underline{\theta})_j|
+|(\delta_-\underline{\theta})_j|+|\theta_j|\|\mathbf{W}_j\|_{L^p}\bigg)$,
from which we obtain the estimate.

\textbf{Estimate on $\Rmnum{2}_j$:  $|\Rmnum{2}_j|\leq C|\theta_j|^2|(\delta_+\underline{\mathbf{W}})_j(-\pi)|$}.\\
This is straightforward.

\textbf{Estimate on $\Rmnum{3}_j$:  $|\Rmnum{3}_j|\leq C|\theta_j||(\delta_+\partial_x\underline{\mathbf{W}})_j(-\pi)|$}.\\
This is straightforward.

\textbf{Estimate on $\Rmnum{4}_j$:  $|\Rmnum{4}_j|\leq C\left[|(\delta_+\underline{\theta})_j|^2+
|(\delta_-\underline{\theta})_{j}|^2+|(\delta_+\underline{\theta})_{j}+(\delta_-\underline{\theta})_{j}|
\bigg(|\theta_{j+1}|^3+|\theta_{j}|^3+|\theta_{j-1}|^3\bigg)\right]$}.\\
We first simplify $\Rmnum{4}_j$ and obtain
\beqs
\bld
\Rmnum{4}_j=
&\frac{1}{2}(\mathbf{u}_{\star}^\prime(\pi-\theta_{j+1})-\mathbf{u}_{\star}^\prime(\pi-\theta_{j-1}),D\mathbf{u}_{\mathrm{ad}}(\pi-\theta_j))-
\frac{1}{2}(\mathbf{u}_\star(\pi-\theta_{j+1})-\mathbf{u}_\star(\pi-\theta_{j-1}),D\mathbf{u}_{\mathrm{ad}}^\prime(\pi-\theta_j)).\\
\eld
\eeqs
Then, it is not hard to see that


\beqs
\bld
 &\left|\frac{1}{2}\bigg(\mathbf{u}_{\star}^\prime(\pi-\theta_{j+1})-\mathbf{u}_{\star}^\prime(\pi-\theta_{j-1}),D\mathbf{u}_{\mathrm{ad}}(\pi-\theta_j)\bigg)
 - \frac{1}{2}\bigg(\mathbf{u}_{\star,\theta\theta}(\pi)(\theta_{j-1}-\theta_{j+1}),-D\mathbf{u}_{\mathrm{ad}}^\prime(\pi)\theta_j\bigg)\right|\\
\leq & C\bigg(|\theta_j||\theta_{j+1}^3-\theta_{j-1}^3|+|\theta_j|^3|\theta_{j+1}-\theta_{j-1}|\bigg),\\
 &\left|\frac{1}{2}\bigg(\mathbf{u}_\star(\pi-\theta_{j+1})-\mathbf{u}_\star(\pi-\theta_{j-1}),D\mathbf{u}_{\mathrm{ad}}^\prime(\pi-\theta_j)\bigg)
-\frac{1}{2}\bigg(\frac{1}{2}\mathbf{u}_{\star,\theta\theta}(\pi)(\theta_{j+1}^2-\theta_{j-1}^2),D\mathbf{u}_{\mathrm{ad}}^\prime(\pi)\bigg)\right|\\
\leq & C\bigg(|\theta_{j+1}^4-\theta_{j-1}^4|+|\theta_j|^2|\theta_{j+1}^2-\theta_{j-1}^2|\bigg),\\
&\frac{1}{2}\bigg(\mathbf{u}_{\star,\theta\theta}(\pi)(\theta_{j-1}-\theta_{j+1}),-D\mathbf{u}_{\mathrm{ad}}^\prime(\pi)\theta_j\bigg)
-\frac{1}{2}\bigg(\frac{1}{2}\mathbf{u}_{\star,\theta\theta}(\pi)(\theta_{j+1}^2-\theta_{j-1}^2),D\mathbf{u}_{\mathrm{ad}}^\prime(\pi)\bigg)\\
=&\frac{1}{4}\bigg(\mathbf{u}_{\star,\theta\theta}(\pi),D\mathbf{u}_{\mathrm{ad}}^\prime(\pi)\bigg)\bigg[(\delta_-\underline{\theta})^2_j-
(\delta_+\underline{\theta})^2_j\bigg],\\
\eld
\eeqs
which establishes the estimate on $\Rmnum{4}_j$ as claimed.

\textbf{Estimate on $\Rmnum{5}_j$:  $|\Rmnum{5}_j|\leq C(\varepsilon)\bigg(
|(\delta_+\underline{\theta})_j|^2+|(\delta_-\underline{\theta})_j|^2+\|\mathbf{W}_j^2\|_{L^p}\bigg)$}.\\
Noting that $|\Rmnum{5}_j|\leq C(\varepsilon)\|(\mathbf{W}_j+\mathbf{H}_j)^2\|_{L^p}$ and
applying the estimate of $\mathbf{H}_j$ into the inequality lead to the above estimate.

\textbf{Estimate on $\mathscr{S}_j$:  $|\mathscr{S}_j|\leq C(\varepsilon)$}.\\
This is straightforward.

Combining our estimates of $\Rmnum{1}_j-\Rmnum{5}_j$ and $\mathscr{S}_j$, we obtain the first inequality in \eqref{e:EN}.

Now, we have to show that the estimate of $\mathbf{N}^\mathbf{w}_j$ in \eqref{e:EN} is true.\\
\textbf{Estimate on $\Rmnum{6}_j$:  }
\beqs
|\Rmnum{6}_j|\leq C(\varepsilon)\left[\bigg(|(\delta_+\underline{\theta})_j|+
|(\delta_-\underline{\theta})_j|\bigg)\sum_{k=-1}^1|\theta_{j+k}|+|\theta_j|\|\mathbf{W}_j\|_{L^p}\right].
\eeqs
First, for $f$ $2\pi$-periodic and smooth, we have
\beqs
\bld
|f(x-\theta_1)-f(x-\theta_2)-f^\prime(x)(\theta_2-\theta_1)|\leq C(|\theta_2-\theta_1|^2+|\theta_2||\theta_2-\theta_1|).
\eld
\eeqs
If in addition, $f$ is odd, we have
\beqs
\bld
|f(\theta_1)-f(\theta_2)-f^\prime(0)(\theta_1-\theta_2)|\leq C|\theta_2^3-\theta_1^3|.
\eld
\eeqs
The latter implies that
\beqs
\bld
|c_j-\frac{1}{4}(\delta^2\underline{\theta})_j|\leq & C(\varepsilon)\bigg(|(\delta_+\underline{\theta})_{j}|^2+
|(\delta_-\underline{\theta})_{j}|^2+|\theta_j|\|\mathbf{W}_j\|_{L^p}\bigg).\\
\eld
\eeqs
Moreover, by the former inequality, we have
\beqs
\bld
|\Rmnum{6}_j|\leq & C\bigg(|(\delta_+\underline{\theta})_j|^2+|(\delta_-\underline{\theta})_j|^2+
|\theta_j||(\delta_+\underline{\theta})_j|+|\theta_j||(\delta_-\underline{\theta})_j|\bigg)
+\left|c_jA\psi(x-\theta_j)-\frac{1}{4}(\delta^2\underline{\theta})_jA\psi(x)\right|\\
\leq&C(\varepsilon)\left[\bigg(|(\delta_+\underline{\theta})_j|+
|(\delta_-\underline{\theta})_j|\bigg)\sum_{k=-1}^1|\theta_{j+k}|+|\theta_j|\|\mathbf{W}_j\|_{L^p}\right].\\
\eld
\eeqs
\textbf{Estimate on $\Rmnum{7}_j$: }
\beqs
|\Rmnum{7}_j|\leq C\left[\left(\sum_{k=-1}^1|\theta_{j+k}|\right)\left(\sum_{k=-1}^1|\dot{\theta}_{j+k}|\right)
+|\dot{\theta}_j|\|\mathbf{W}_j\|_{L^p}\right]
\eeqs
Noting that $(\underline{\mathbf{E}}*\underline{\theta})_j$ is the linear part of
$\mathbf{H}_j+\mathbf{u}_\star(x-\theta_j)-\mathbf{u}_\star(x)$ and there is no term invovling $\dot{\mathbf{W}}_j$ in $\Rmnum{7}_j$, we have
\beqs
\bld
|\Rmnum{7}_j|\leq & C\bigg(|\theta_{j+1}||\dot{\theta}_{j+1}|+|\theta_{j-1}||\dot{\theta}_{j-1}|
+|\theta_{j}||\dot{\theta}_{j}|\bigg)+\left|c_j\psi^\prime(x-\theta_j)\dot{\theta}_j\right|
+\left|\frac{1}{4}(\delta^2\dot{\underline{\theta}})_j\psi(x)-\tilde{\dot{c}}_j\psi(x-\theta_j)\right|,\\
\eld
\eeqs
where $\tilde{\dot{c}}_j=\dot{c}_j+\langle\dot{\mathbf{W}}_j(x),\mathbf{u}_{\mathrm{ad}}(x-\theta_j)-\mathbf{u}_{\mathrm{ad}}(x)\rangle$.\\
First, we note that
\beqs
\bld
\left|c_j\psi^\prime(x-\theta_j)\dot{\theta}_j\right|\leq & C|\dot{\theta}_j|\left[|(\delta^2\underline{\theta})_j|
+\bigg(|(\delta_+\underline{\theta})_{j}|+
|(\delta_-\underline{\theta})_{j}|\bigg)\sum_{k=-1}^1\theta_{j+k}+|\theta_j|\|\mathbf{W}_j\|_{L^p}\right].\\
\eld
\eeqs
Moreover, we claim that
\beqs
\bld
|\tilde{\dot{c}}_j|
\leq &C\left[|(\delta_+\dot{\underline{\theta}})_j|+|(\delta_-\dot{\underline{\theta}})_j|+
\left(\sum_{k=-1}^1|\theta_{j+k}|\right)\left(\sum_{k=-1}^1|\dot{\theta}_{j+k}|\right)
+|\dot{\theta}_j|\|\mathbf{W}_j\|_{L^p}\right],\\
|\tilde{\dot{c}}_j-\frac{1}{4}(\delta^2\dot{\underline{\theta}})_j|
\leq &
C\left[\left(\sum_{k=-1}^1|\theta_{j+k}|\right)\left(\sum_{k=-1}^1|\dot{\theta}_{j+k}|\right)
+|\dot{\theta}_j|\|\mathbf{W}_j\|_{L^p}\right].\\
\eld
\eeqs
In fact, we have
\beqs
\bld
&\left|\left\langle\phi(x)(\dot{\mathbf{u}}_\star(x+\theta_j-\theta_{j+1})-
\dot{\mathbf{u}}_\star(x+\theta_j-\theta_{j-1})),\mathbf{u}_{\mathrm{ad}}(x)\right\rangle\right|
\leq C\bigg(|(\delta_+\underline{\theta})_{j}||(\delta_+\dot{\underline{\theta}})_{j}|+
|(\delta_-\underline{\theta})_{j}||(\delta_-\dot{\underline{\theta}})_{j}|\bigg),\\
&|\langle\phi^\prime(x+\theta_j)\dot{\theta}_j(\mathbf{u}_\star(x+
\theta_j-\theta_{j+1})-\mathbf{u}_\star(x+\theta_j-\theta_{j-1})),
\mathbf{u}_{\mathrm{ad}}(x)\rangle|
\leq C|\dot{\theta}_j|\bigg(|(\delta_+\underline{\theta})_{j}|+|(\delta_-\underline{\theta})_{j}|\bigg),\\
&\left|\left\langle(\phi(x+\theta_j)-\phi(x))(\dot{\mathbf{u}}_\star(x+\theta_j-\theta_{j+1})-\dot{\mathbf{u}}_\star(x+\theta_j-\theta_{j-1})),
\mathbf{u}_{\mathrm{ad}}(x)\right\rangle\right|
\leq C|\theta_j|\bigg(|(\delta_+\dot{\underline{\theta}})_{j}|+|(\delta_-\dot{\underline{\theta}})_{j}|\bigg),\\
&\left|\left\langle(\dot{\mathbf{u}}_\star(x+\theta_j-\theta_{j+1})+
\dot{\mathbf{u}}_\star(x+\theta_j-\theta_{j-1})),\mathbf{u}_{\mathrm{ad}}(x)\right\rangle\right|
\leq C\bigg(|(\delta_+\dot{\underline{\theta}})_j|+|(\delta_-\dot{\underline{\theta}})_j|\bigg),\\
&|\left\langle(\dot{\mathbf{u}}_\star(x+\theta_j-\theta_{j+1})+\dot{\mathbf{u}}_\star(x+
\theta_j-\theta_{j-1})),\mathbf{u}_{\mathrm{ad}}(x)\right\rangle
+\delta^2\underline{\dot{\theta}}_j|
\leq C\bigg(|(\delta_+\underline{\theta})_j||(\delta_+\dot{\underline{\theta}})_j|+
|(\delta_-\underline{\theta})_j||(\delta_-\dot{\underline{\theta}})_j|\bigg),
\eld
\eeqs
which establishes the claim and thus the estimate on $\Rmnum{7}_j$.

\textbf{Estimate on $\Rmnum{8}_j$: }
\beqs
|\Rmnum{8}_j|\leq C\bigg(|\mathbf{N}^\theta_j|+|\mathbf{N}^\theta_{j+1}|+|\mathbf{N}^\theta_{j-1}|\bigg)
\eeqs
The calculation is straightforward using the expressions for $\mathbf{K}_j$ and $\dot{\theta}_j$.

\textbf{Estimate on $\Rmnum{9}_j$:}
\beqs
|\Rmnum{9}_j|\leq C(\varepsilon)\left[\left(\sum_{k=0}^1|(\delta_-\underline{\theta})_{j+k}|\right)\left(\sum_{k=-1}^1|\theta_{j+k}\right)+
+|\theta_j||\mathbf{W}_j|+|\theta_j|^2\|\mathbf{W}_j\|_{L^p}+|\mathbf{W}_j|^2
\right]
\eeqs
The calculation is straightforward using the estimate on $\mathbf{H}_j$.

\textbf{Estimate on $\frac{\partial\mathbf{G}_j}{\partial\mathbf{W}_j}\Rmnum{10}_j$:}
\beqs
|\frac{\partial\mathbf{G}_j}{\partial\mathbf{W}_j}\Rmnum{10}_j|
\leq C(\varepsilon)\left(|\theta_j|\sum_{k=-1}^1\bigg(|(\delta_+\underline{\theta})_{j+k}|+|(\delta_+\underline{\mathbf{W}})_{j+k}(-\pi)|\bigg)+
 |\langle A\mathbf{W}_j(x),\mathbf{u}_{\mathrm{ad}}(x-\theta_j)-\mathbf{u}_{\mathrm{ad}}(x)\rangle|\right).
\eeqs
Integrating by parts, we have
\beqs
\bld
\langle A\mathbf{W}_j(x),\mathbf{u}_{\mathrm{ad}}(x-\theta_j)-\mathbf{u}_{\mathrm{ad}}(x)\rangle=&
\langle\mathbf{W}_j(x),A^*\left(\mathbf{u}_{\mathrm{ad}}(x-\theta_j)-\mathbf{u}_{\mathrm{ad}}(x)\right)\rangle+\\
&(\partial_x\mathbf{W}_{j+1}(-\pi)-\partial_x\mathbf{W}_{j}(-\pi),D(\mathbf{u}_{\mathrm{ad}}(\pi-\theta_j)-\mathbf{u}_{\mathrm{ad}}(\pi)))-\\
&(\mathbf{W}_{j+1}(-\pi)-\mathbf{W}_{j}(-\pi),D(\mathbf{u}^\prime_{\mathrm{ad}}(\pi-\theta_j)-\mathbf{u}^\prime_{\mathrm{ad}}(\pi))).
\eld
\eeqs
Thereofore, we have
\beqs
|\frac{\partial\mathbf{G}_j}{\partial\mathbf{W}_j}\Rmnum{10}_j|
\leq C(\varepsilon)|\theta_j|\left(\sum_{k=-1}^1\bigg(|(\delta_+\underline{\theta})_{j+k}|+|(\delta_+\underline{\mathbf{W}})_{j+k}(-\pi)|\bigg)+
 \|\mathbf{W}_j\|_{L^p}+|(\delta_+\partial_x\underline{\mathbf{W}})_{j}(-\pi)|\right).
\eeqs
\textbf{Estimate on $(\id-\frac{\partial\mathbf{G}_j}{\partial\mathbf{W}_j})^{-1}$:}
For any $\underline{\theta}\in\ell^\infty$ and $p\in[1, \infty]$,
there exists a constant $C>0$ such that
\beqs
\opnorm{(\id-\frac{\partial\mathbf{G}_j}{\partial\mathbf{W}_j})^{-1}}_{L^p}\leq C.
\eeqs

Combining estimates on $\Rmnum{6}_j$ to $\Rmnum{9}_j$, $\frac{\partial \mathbf{G}_j}{\partial \mathbf{W}_j}$
and $(\id-\frac{\partial\mathbf{G}_j}{\partial\mathbf{W}_j})^{-1}$, we obtain the second inequality in \eqref{e:EN}.
\epf

Moreover, we have the following lemma.
\bl
\label{l:61e} There exist $C>0$ and $\eta>0$ such that, for all $(\underline{\theta},\underline{\mathbf{W}})\in Y$ with its $Y$-norm smaller than
$\eta$, we have
\beqs
\bld
\|\mathbf{N}^\theta(s)\|_{\ell^1}&\leq \frac{C}{(1+s)^{\frac{3}{2}}}\|(\underline{\theta}(t), \underline{\mathbf{W}}(t))\|_Y^2+
\frac{C}{(1+s)^{\frac{5}{4}}}\|(\underline{\theta}(t), \underline{\mathbf{W}}(t))\|_Y(1+s)\|\delta_+\partial_{xx}\underline{\mathbf{W}}(s)\|_{X_2},\\
&\\
\|\mathbf{N}^\theta(s)\|_{\ell^2}&\leq \frac{C}{(1+s)^{\frac{3}{2}}}\|(\underline{\theta}(t), \underline{\mathbf{W}}(t))\|_Y^2+
\frac{C}{(1+s)^{\frac{3}{2}}}\|(\underline{\theta}(t), \underline{\mathbf{W}}(t))\|_Y(1+s)\|\delta_+\partial_{xx}\underline{\mathbf{W}}(s)\|_{X_2},\\
&\\
\|\mathbf{N}^{\mathbf{w}}(s)\|_{X_1}&\leq \frac{C}{1+s}\|(\underline{\theta}(t), \underline{\mathbf{W}}(t))\|_Y^2+
\frac{C}{(1+s)^{\frac{5}{4}}}\|(\underline{\theta}(t), \underline{\mathbf{W}}(t))\|_Y(1+s)\|\delta_+\partial_{xx}\underline{\mathbf{W}}(s)\|_{X_2},\\
&\\
\|\mathbf{N}^{\mathbf{w}}(s)\|_{X_2}&\leq \frac{C}{(1+s)^{\frac{5}{4}}}\|(\underline{\theta}(t), \underline{\mathbf{W}}(t))\|_Y^2+
\frac{C}{(1+s)^{\frac{3}{2}}}\|(\underline{\theta}(t), \underline{\mathbf{W}}(t))\|_Y(1+s)\|\delta_+\partial_{xx}\underline{\mathbf{W}}(s)\|_{X_2}.
\eld
\eeqs
\el
\bpf
The estimates are obtained through a direct calculation from the estimates in Lemma \ref{l:61}. 
We sketch the computation for $\|\mathbf{N}^\theta(s)\|_{\ell^1}$, and
the others follow similarly.

First, for terms only involving $\underline{\theta}$, we notice that
\beqs
\label{e:NLE11}
\bld
&\sum_{j\in\Z}|(\delta_+\underline{\theta})_j|^2=-\sum_{j\in\Z}\theta_j(\delta^2\underline{\theta})_j\leq
\|\underline{\theta}\|_{\ell^2}\|\delta^2\underline{\theta}\|_{\ell^2}\leq
\frac{1}{(1+s)^{\frac{3}{2}}}\|(\underline{\theta}(t), \underline{\mathbf{W}}(t))\|_Y^2,\\
&\sum_{j\in\Z}|\theta_j|^3|(\delta_+\underline{\theta})_j|\leq \|\underline{\theta}\|_{\ell^\infty}^2
\|\underline{\theta}\|_{\ell^2}\|\delta_+\underline{\theta}\|_{\ell^2}\leq \|\underline{\theta}\|_{\ell^\infty}^2
\|\underline{\theta}\|_{\ell^2}^{\frac{3}{2}}\|\delta^2\underline{\theta}\|_{\ell^2}^{\frac{1}{2}}\leq \frac{1}{(1+s)^2}
\|(\underline{\theta}(t), \underline{\mathbf{W}}(t))\|_Y^4.\\
\eld
\eeqs
Second, for terms involving $\underline{\mathbf{W}}$, we observe that
\beqs
\label{e:NLE12}
\bld
\sum_{j\in\Z}|\theta_j||(\delta_+\partial_x\underline{\mathbf{W}})_j(-\pi)|&\leq
\|\underline{\theta}\|_{\ell^2}\left(\sum_{j\in\Z}\left(\int_{-\pi}^\pi\left(\partial_{xx}\mathbf{W}_{j+1}(x)
-\partial_{xx}\mathbf{W}_j(x)\right)\rmd x\right)^2\right)^{\frac{1}{2}}\\
&\leq \sqrt{2\pi}\|\underline{\theta}\|_{\ell^2}\|\delta_+\partial_{xx}\underline{\mathbf{W}}\|_{X_2}\\
&\leq \frac{\sqrt{2\pi}}{(1+s)^{\frac{5}{4}}}\|(\underline{\theta}(t), \underline{\mathbf{W}}(t))\|_Y
(1+s)\|\delta_+\partial_{xx}\underline{\mathbf{W}}(s)\|_{X_2}.\\
\eld
\eeqs
Similarly, for $\sum_{j\in\Z}|\theta_j||(\delta_+\underline{\mathbf{W}})_j(-\pi)|$, we have
\beqs
\sum_{j\in\Z}|\theta_j||(\delta_+\underline{\mathbf{W}})_j(-\pi)|\leq
\sqrt{2\pi}\|\underline{\theta}\|_{\ell^2}\|\delta_+\partial_{x}\underline{\mathbf{W}}\|_{X_2}.
\eeqs
Using the ``homogeneous matching boundary conditions " \eqref{e:33a}, we have
\beqs
\bld
\|\delta_+\partial_{x}\underline{\mathbf{W}}\|_{X_2}&=\left(-\sum_{j\in\Z}\int_{-\pi}^\pi\left(\delta_+\underline{\mathbf{W}}\right)_j(x)
\left(\delta_+\partial_{xx}\underline{\mathbf{W}}\right)_j(x)\rmd x\right)^{\frac{1}{2}}\\
&\leq \|\delta_+\underline{\mathbf{W}}\|_{X_2}^{\frac{1}{2}}\|\delta_+\partial_{xx}\underline{\mathbf{W}}\|_{X_2}^{\frac{1}{2}}\\
&\leq \|\delta_+\underline{\mathbf{W}}\|_{X_2}+\|\delta_+\partial_{xx}\underline{\mathbf{W}}\|_{X_2}.
\eld
\eeqs
We plug the latter estimate into the former one and obtain that
\beqs
\label{e:NLE13}
\bld
\sum_{j\in\Z}|\theta_j||(\delta_+\underline{\mathbf{W}})_j(-\pi)|&\leq
\sqrt{2\pi}\|\underline{\theta}\|_{\ell^2}\bigg( \|\delta_+\underline{\mathbf{W}}\|_{X_2}+\|\delta_+\partial_{xx}\underline{\mathbf{W}}\|_{X_2}\bigg)\\
&\leq \frac{\sqrt{2\pi}}{(1+s)^{\frac{3}{2}}}\|(\underline{\theta}(t), \underline{\mathbf{W}}(t))\|_Y^2+
\frac{\sqrt{2\pi}}{(1+s)^{\frac{5}{4}}}\|(\underline{\theta}(t), \underline{\mathbf{W}}(t))\|_Y
(1+s)\|\delta_+\partial_{xx}\underline{\mathbf{W}}(s)\|_{X_2}.
\eld
\eeqs
For $\sum_{j\in\Z}\|\mathbf{W}_j\|_{L^p}|(\delta_+\underline{\mathbf{W}})_j(-\pi)|$, we take $p=2$ and follow steps as above, obtaining
the following estimate.
\beqs
\label{e:NLE14}
\sum_{j\in\Z}\|\mathbf{W}_j\|_{L^2}|(\delta_+\underline{\mathbf{W}})_j(-\pi)|\leq
\frac{\sqrt{2\pi}}{(1+s)^{2}}\|(\underline{\theta}(t), \underline{\mathbf{W}}(t))\|_Y^2+
\frac{\sqrt{2\pi}}{(1+s)^{\frac{7}{4}}}\|(\underline{\theta}(t), \underline{\mathbf{W}}(t))\|_Y
(1+s)\|\delta_+\partial_{xx}\underline{\mathbf{W}}(s)\|_{X_2}.
\eeqs
For $\sum_{j\in\Z}\|\mathbf{W}_j^2\|_{L^p}|$, we take $p=1$ and obtain that
\beqs
\sum_{j\in\Z}\|\mathbf{W}_j^2\|_{L^1}\leq \|\underline{\mathbf{W}}\|_{X_2}^2\leq \frac{1}{(1+s)^{\frac{3}{2}}}
\|(\underline{\theta}(t), \underline{\mathbf{W}}(t))\|_Y^2.
\eeqs
Combining the above estimate, we establish the first inequality in the lemma.
\epf
\subsection{Bloch wave decomposition}\label{ss:62}
In this section, we present the Bloch wave decomposition of the linear operator $\widetilde{A}$. We first
recall that $\widetilde{A}$, as in \eqref{e:opL2}, is defined as
\beqs
\begin{matrix}
 \widetilde{A}:& (H^2(\R))^n&\longrightarrow&(L^2(\R))^n\\
 &\mathbf{v}&\longmapsto&D\partial_{xx}\mathbf{v}-\mathbf{f}^\prime(\mathbf{u}_\star)\mathbf{v}.
\end{matrix}
\eeqs
We introduce the direct integral \cite[\Rmnum{13}.16.]{reedsimon}
\begin{equation}
\label{e:DI}
 \begin{matrix}
  \mathscr{B}: &L^2(\mathbb{T}_1, (L^2(\mathbb{T}_{2\pi}))^n)&\longrightarrow &(L^2(\mathbb{R}))^n\\
  &\mathbf{U}(\sigma,x)&\longmapsto& \int_{\sigma\in\mathbb{T}_1}\rme^{\rmi\sigma\cdot x}\mathbf{U}(\sigma,\cdot)\rmd\sigma
 \end{matrix}.
\end{equation}
The direct interal is an isometric isomorphism with inverse
\begin{equation*}
\begin{matrix}
 \mathscr{B}^{-1}: &(L^2(\mathbb{R}))^n &\longrightarrow&L^2(\mathbb{T}_1, (L^2(\mathbb{T}_{2\pi}))^n)\\
 &\mathbf{u}(x)&\longmapsto&\frac{1}{2\pi}\sum_{\ell\in \mathbb{Z}^m}\rme^{\rmi\ell\cdot x}\widehat{\mathbf{u}}(\sigma+\ell).
\end{matrix}
\end{equation*}

The following result from \cite{scarpellini_1994, mielke_1997} characterizes the Bloch wave decomposition of $\widetilde{A}$.
\begin{theorem}[\bf Bloch wave decomposition]
\label{t:4}The linear operator $\widetilde{A}$ is diagonal in Bloch wave space. To be precise,
\begin{equation}
\mathscr{B}^{-1}\widetilde{A}\mathscr{B}=\widehat{A}=\int_{-\frac{1}{2}}^{\frac{1}{2}}B(\sigma)\rmd\sigma,
\end{equation}
where by $\widehat{A}=\int_{-\frac{1}{2}}^{\frac{1}{2}}B(\sigma)\rmd\sigma$,
we mean that, given any $\mathbf{u}\in L^2(\mathbb{T}_1, (L^2(\mathbb{T}_{2\pi}))^n)$,
\beqs
(\widehat{A}\mathbf{u})(\sigma)=B(\sigma)\mathbf{u}(\sigma),\text{ a.e. }\sigma\in[-\frac{1}{2},\frac{1}{2}].
\eeqs
Moreover, we have the following spectral mapping property.
\begin{equation}
 \spec(\widetilde{A})=\spec(\widehat{A})=\bigcup_{\sigma\in[-\frac{1}{2}, \frac{1}{2}]}\spec(B(\sigma)).
\end{equation}
\end{theorem}
\subsection{ Spectral properties of \texorpdfstring{$\{\widehat{A}_{\mathrm{ch}}(\sigma)\}
            _{\sigma\in[-\frac{1}{2},\frac{1}{2}]}$}{A(sigma)}  }
\label{ss:64}
We recall that $\widehat{A}_{\mathrm{ch}}(\sigma)$ is defined in \eqref{e:Asigma} as
$\widehat{A}_{\mathrm{ch}}(\sigma)=\mathscr{F}_nB(\sigma)\mathscr{F}_n^{-1}$ and $Y_q$ in \eqref{e:space} for $1\leq q\leq \infty$.
We are concerned with their spectral properties as unbounded operators in $Y_q$,
which is useful for the derivation of the estimates for $M(t,\sigma)$ as defined in \eqref{e:M}.

We first show the well-definedness of $\widehat{A}_{\mathrm{ch}}(\sigma)$ in $Y_q$ in the following lemma.
\bl
\label{l:64}
For any given $\sigma\in [-\frac{1}{2},\frac{1}{2}]$, $\widehat{A}_{\mathrm{ch}}(\sigma)$ is an unbounded
closed operator in $Y_2$, that is,
\begin{equation}
\begin{matrix}
\widehat{A}_{\mathrm{ch}}(\sigma):&\mathscr{D}_2(\widehat{A}_{\mathrm{ch}}(\sigma))\subset Y_2&\longrightarrow &Y_2\\
&\underline{\mathbf{w}}&\longmapsto &\{-(\sigma+\ell)^2D\mathbf{w}_\ell+
\sum_{k\in\Z}\mathbf{h}_{\ell-k}\mathbf{w}_k\}_{\ell\in\Z},
\end{matrix}
\end{equation}
where $\mathscr{D}_2(\widehat{A}_{\mathrm{ch}}(\sigma))=\{\mathbf{w}\in Y_2\mid\{(1+m^2)\mathbf{w}_m\}_{m\in\mathbb{Z}}\in Y_2\}$
and $\underline{\mathbf{h}}=\{\mathbf{h}_\ell\}_{\ell\in\Z}=
\frac{1}{2\pi}\int_{-\pi}^{\pi}\mathbf{f}^\prime(\mathbf{u}_\star(x))\rme^{-\rmi kx} \rmd x$.
Moreover, $\widehat{A}_{\mathrm{ch}}(\sigma)$ can naturally be considered as an unbounded closed operator in $Y_q$,
with
$\mathscr{D}_q(\widehat{A}_{\mathrm{ch}}(\sigma))=\{\mathbf{w}\in Y_q\mid\{(1+m^2)\mathbf{w}_m\}_{m\in\mathbb{Z}}\in Y_q\}$,
for all $1\leq q \leq \infty$.
\el
\bpf
The expression for $\widehat{A}_{\mathrm{ch}}(\sigma)$ in $Y_2$ follows from a direct calculation.
The extension to $Y_q$ follows from the fact that the set
$\{\underline{\mathbf{w}}\in Y_\infty\mid \underline{\mathbf{w}} \text{ has finitely many nonzero entries}\}$ is dense
in $Y_q$ and $\mathscr{D}_q(\widehat{A}_{\mathrm{ch}}(\sigma))$, for all $q\in[1,\infty]$.
\epf

We then have the following proposition.
\bp
\label{p:601}
For any fixed $\sigma\in[-\frac{1}{2},\frac{1}{2}]$ and $p\in[1,\infty]$,
$\widehat{A}_{\mathrm{ch}}(\sigma)$ defined in $Y_q$ is
sectorial and has compact resolvent. In fact,
there exist $C>0$, $\omega\in(\pi/2,\pi)$ and $\lambda_0\in\R$, independent of $\sigma$ and $q$, such that
the sector $S(\lambda_0,\omega)=\{\lambda\in\C\mid0\leq |\arg(\lambda-\lambda_0)|
\leq\omega, \lambda\neq \lambda_0\}\subseteq \rho(\widehat{A}_{\mathrm{ch}}(\sigma))$ and
\beq
\opnorm{(\widehat{A}_{\mathrm{ch}}(\sigma)-\lambda)^{-1}}_{Y_q}\leq C|\lambda-\lambda_0|^{-1}, \text{ for all }
\lambda\in S(\lambda_0,\omega),\sigma\in[-\frac{1}{2},\frac{1}{2}]\text{ and }q\in[1,\infty].
\eeq
Moreover, for any fixed $\sigma\in[-\frac{1}{2},\frac{1}{2}]$, the spectrum of
$\widehat{A}_{\mathrm{ch}}(\sigma)$ is independent of the choice of its underlying space $Y_q$
and thus denoted as $\spec(\widehat{A}_{\mathrm{ch}}(\sigma))$, for any
$q\in[1,\infty]$, with $\spec(\widehat{A}_{\mathrm{ch}}(\sigma))=\spec(B(\sigma))$
consisting only of isolated eigenvalues with finite multiplicity.
\ep
\bpf
We view $\widehat{A}_{\mathrm{ch}}(\sigma)$ as a perturbation of the Laplacian in the discrete Fourier space, that is,
\beqs
\widehat{A}_{\mathrm{ch}}(\sigma)=L(\sigma)+H,
\eeqs
where $L(\sigma)\underline{\mathbf{w}}=\{-(\sigma+\ell)^2D\mathbf{w}_\ell\}_{\ell\in\Z}$ and
$H\underline{\mathbf{w}}=\{\sum_{k\in\Z}\mathbf{h}_{\ell-k}\mathbf{w}_k\}_{\ell\in\Z}$.
It is straightforward to verify that the proposition holds for the Laplacian $L(\sigma)$.
We only have to show that the perturbation $H$ is good enough to preserve these properties.
Noting that $H\in \mathscr{L}((\ell^p)^n)$ for any $p\in[1,\infty]$ with its norm uniformly bounded,
we have, for $\lambda\in \rho(L(\sigma))$, $|\lambda|$ sufficiently large,
\begin{equation}
\label{app:01}
(\widehat{A}_{\mathrm{ch}}(\sigma)-\lambda)^{-1}=(L(\sigma)+H-\lambda)^{-1}=(L(\sigma)-\lambda)^{-1}(\id+H(L(\sigma)-\lambda)^{-1})^{-1}.
\end{equation}
All assertions in the proposition easily follows from this expression \eqref{app:01}, except for the fact that
the spectrum of $\widehat{A}_{\mathrm{ch}}(\sigma)$ is independent of $q$.

To prove this property, we denote the spectrum of $\widehat{A}_{\mathrm{ch}}(\sigma)$
defined on $Y_q$ as $\spec(\widehat{A}_{\mathrm{ch}}(\sigma),q)$, which consists of eigenvalues with finite multiplicity, accumulating at
infinity, only. Given any eigenfunction $\underline{\mathbf{v}}=\{\mathbf{v}_j\}_{j\in\Z}$,
$\underline{\mathbf{v}}$ belongs to $\bigcap_{q\in[1,\infty]}Y_q$
since $\underline{\mathbf{v}}$ are smooth, that is, $\mathbf{v}_j$ decays algebraically with any rate.
This establishes $\spec(\widehat{A}(\sigma),q)=\spec(\widehat{A}(\sigma),p)$, for any $p, q\in[1,\infty]$.
\epf

\subsection{Perturbation results}
\label{ss:65}
We apply perturbation theory to the Bloch wave operator $B(\sigma)$ for $\sigma$ near $0$
and obtain more detailed spectral information, including the Taylor expansion of $d$ in Hypotheses \ref{h:22}.

To this end, we define
\beqs
\begin{matrix}
F: &[-\frac{1}{2},\frac{1}{2}]\times\mathbb{C}\times H^2_{\perp}&\longrightarrow &L^2\\
&(\sigma, \lambda, \mathbf{w})&\longmapsto &(B(\sigma)-\lambda )(\mathbf{w}+\mathbf{u}_\star^\prime),
\end{matrix}
\eeqs
where $H^2_{\perp}=\{\mathbf{w}\in (H^2(\T_{2\pi}))^n\mid \langle\mathbf{w}, \mathbf{u}_\star^\prime\rangle=0\}.$
A standard implicit-function-theorem argument shows that there are a small neighborhood of $\sigma$ at the origin and a
smooth function $(\lambda(\sigma), \mathbf{w}(\sigma))$ with $(\lambda(\sigma), \mathbf{w}(\sigma))=0$ on this neighborhood
such that $F(\sigma, \lambda(\sigma),\mathbf{w}(\sigma))=0$. We denote
$\mathbf{e}(\sigma)=\mathbf{u}_\star^\prime+\mathbf{w}(\sigma)$. Similarly, replacing $B(\sigma)$ with
its adjoint $B^*(\sigma)$, we obtain a smooth continuation of $\mathbf{u}_{\mathrm{ad}}$,
denoted as $\mathbf{e}^*(\sigma)$. Without loss of generality, we can assume that
 $\langle\mathbf{e}(\sigma),\mathbf{e}^*(\sigma)\rangle=1$. Moreover, we have the following proposition.
\bp
\label{p:604}
There exist positive numbers $\gamma_0$ and $\gamma_1$
such that for any $|\sigma|\leq \gamma_0$ in $\R$,
$B(\sigma)$ has only one simple eigenvalue within the strip $|\re\lambda|\leq \gamma_1$ in $\C$,
which is exactly the continuation $\lambda(\sigma)$ of the eigenvalue $\lambda(0)=0$.
Moreover, $\lambda(\sigma)$ has the Taylor expansion,
\beqs
\lambda(\sigma)=-d\sigma^2+\rmO(|\sigma|^3),
\eeqs
where $-\gamma_1/4\leq-2d\sigma^2<Re\lambda(\sigma)<-\frac{d}{2}\sigma^2$,
for all $\sigma\in[-\gamma_0,\gamma_0]$ and
\beqs
d=-\langle 2\rmi \frac{\partial^2\mathbf{e}(0,x)}{\partial x\partial\sigma}
-\mathbf{u}_\star^\prime(x),D\mathbf{u}_{\mathrm{ad}}(x)\rangle.
\eeqs
\ep
\bpf
We first derive the explicit expression of $d$. To do that,
taking first and second derivative with respect to $\sigma$
of $F(\sigma,\lambda(\sigma),\mathbf{w}(\sigma))=0$,
taking the inner product of the derivatives with $\mathbf{u}_{\mathrm{ad}}$ and letting $\sigma=0$,
we have
\beqs
\bld
&\lambda^\prime(0)=\langle B(0)\partial_\sigma\mathbf{e}(0,x)+
2\rmi D\mathbf{u}_\star^{\prime\prime}(x),\mathbf{u}_{\mathrm{ad}}(x)\rangle,\\
&\lambda^{\prime\prime}(0)=\langle B(0)\partial_\sigma^2\mathbf{e}(0,x)+(4\rmi D\partial_x-2\lambda^\prime(0))\partial_\sigma\mathbf{e}(0,x)-
2D\mathbf{u}_\star^\prime(x),\mathbf{u}_{\mathrm{ad}}(x)\rangle.
\eld
\eeqs
Noting that $\span\{\mathbf{u}_{\mathrm{ad}}\}\perp Rg(B(0))$ and
the inner product of an even function and an odd function is always 0, we have
\beqs
\bld
\lambda^\prime(0)=0, \quad \lambda^{\prime\prime}(0)
=2\langle 2\rmi \frac{\partial^2\mathbf{e}(0,x)}{\partial x\partial\sigma}
-\mathbf{u}_\star^\prime(x),D\mathbf{u}_{\mathrm{ad}}(x)\rangle.
\eld
\eeqs
It remains to prove the uniqueness of the eigenvalue of $B(\sigma)$ in a vertical strip
centered at the origin for sufficiently small $\sigma$.
First, there is no eigenvalue within the strip far away from the origin due to the fact that, by Proposition \ref{p:601},
$\spec(B(\sigma))$ is in the same sector for every $\sigma\in [-\frac{1}{2}, \frac{1}{2}]$.
Secondly, the uniqueness within a small neighborhood of the origin follows from the above perturbation results.
For the region inbetween, compactness and the local robustness of resolvent guarantee the absence of eigenvalues within this area.
\epf
\br
\begin{enumerate}
 \item We stress that we may choose $\gamma_0$ as small as desired.
 \item The uniqueness implies that, for $|\sigma|$ sufficiently small, $\lambda(\sigma)$ is a real number since its complex conjugate is also
 an eigenvalue.
\end{enumerate}
\er

\subsection{Properties of analytic semigroups \texorpdfstring{$\{\rme^{\widehat{A}_{\mathrm{ch}}(\sigma) t}\}
_{\sigma\in[-\frac{1}{2},\frac{1}{2}]}$}{of A(sigma)'s}}\label{ss:66}
In this section, we will derive various estimates on $\rme^{\widehat{A}_{\mathrm{ch}}(\sigma) t}$.
We first note that by \cite[1.4]{henry} the interpolation space  $\mathscr{D}_q(\widehat{A}_{\mathrm{ch}}(\sigma)^\alpha)$
is independent of $\sigma$,
\beqs
\mathscr{D}_q(\widehat{A}_{\mathrm{ch}}(\sigma)^\alpha)=
\{\underline{\mathbf{w}}\in Y_q\mid \{(1+m^2)^\alpha\mathbf{w}_m\}_{m\in\mathbb{Z}}\in Y_q\}=:Y^\alpha_q,
\quad \|\underline{\mathbf{w}}\|_{Y^\alpha_q}=\|\{(1+m^2)^\alpha \mathbf{w}_m\}_{m\in\mathbb{Z}}\|_{Y_q}.
\eeqs
We then recall the definitions of
$Y_{q,\mathrm{c}}(\sigma)$, $Y_{q,\mathrm{s}}(\sigma)$, $\widehat{A}_{\mathrm{c}}(\sigma)$ and $\widehat{A}_{\mathrm{s}}(\sigma)$
from \eqref{e:cs}.
We now have the following proposition.
\bp
\label{p:605}
For every $q\in[1, +\infty]$ and $\alpha>0$,
there exist positive constants $\epsilon\in(0,1)$, $\gamma_2$, $C(q)$, $C(\alpha)$ and $C(\alpha, q)$ such that
\begin{equation*}
\bld
\opnorm{\rme^{\widehat{A}_{\mathrm{c}}(\sigma)t}}_{Y_{q,\mathrm{c}}(\sigma)}
&\leq \rme^{-\frac{d}{2}\sigma^2 t}, \text{ for all }|\sigma|\leq \gamma_0,t\geq 0,\\
\opnorm{\rme^{\widehat{A}_{\mathrm{c}}(\sigma)t}}_{Y_{q,\mathrm{c}(\sigma)}\rightarrow Y_q^\alpha}&\leq C(\alpha)\rme^{-\frac{d}{2}\sigma^2 t},
\text{ for all }|\sigma|\leq \gamma_0,t\geq 0,\\
\opnorm{\rme^{\widehat{A}_{\mathrm{s}}(\sigma)t}}_{Y_{q,\mathrm{s}}(\sigma)} &\leq
C(q)e^{-\frac{\gamma_1}{2} t},\text{ for all }|\sigma|\leq \gamma_0,t\geq 0,\\
\opnorm{\rme^{\widehat{A}_{\mathrm{s}}(\sigma)t}}_{Y_{q,\mathrm{s}}(\sigma)\rightarrow Y_q^\alpha}&\leq C(\alpha,q)t^{-\alpha}\rme^{-\gamma_1 t/2},
\text{ for all }|\sigma|\leq \gamma_0,t>0,\\
\opnorm{\rme^{\widehat{A}_{\mathrm{ch}}(\sigma)t}}_{Y_q}&\leq C(q)\rme^{-\epsilon d\sigma^2 t},\text{ for all }|\sigma|\leq \gamma_0,t\geq 0,\\
\opnorm{\rme^{\widehat{A}_{\mathrm{ch}}(\sigma)t}}_{Y_q}&\leq C(q)\rme^{-\gamma_2 t},\text{ for all }\gamma_0\leq|\sigma|\leq\frac{1}{2}, t\geq 0,\\
\opnorm{\rme^{\widehat{A}_{\mathrm{ch}}(\sigma)t}}_{Y_q\rightarrow Y_q^\alpha}&\leq C(\alpha,p)t^{-\alpha}\rme^{-\gamma_2 t},
\text{ for all }\gamma_0\leq|\sigma|\leq\frac{1}{2}, t> 0.
\eld
\end{equation*}
\ep
\bpf
We first derive estimates for the case $|\sigma|\leq \gamma_0$.
For $\widehat{A}_{\mathrm{c}}(\sigma)$, we have $\rme^{\widehat{A}_{\mathrm{c}}(\sigma)t}=\rme^{\lambda(\sigma)t}$.
The first two inequalities follow directly from the fact that $\re\lambda(\sigma)<-\frac{d}{2}\sigma^2$ and $\mathbf{e}(\sigma)$ is smooth,
by Proposition \ref{p:604}, for $|\sigma|\leq \gamma_0$.

For $\widehat{A}_{\mathrm{s}}(\sigma)$, by Proposition \ref{p:601} and \ref{p:604}, for any $\sigma\in(-\gamma_0,\gamma_0)$ and $q\in[1,\infty]$,
\beqs
\spec(\widehat{A}_{\mathrm{s}}(\sigma),q)\subset \C\backslash S(-\frac{\gamma_1}{2},\tilde{\omega}), \text{ where }\tilde{\omega}\in(\frac{\pi}{2},\pi).
\eeqs
Moreover, for every $q\in[1, +\infty]$,
there exists a positive constant $C(q)$ such that
$$\opnorm{(\widehat{A}_{\mathrm{s}}(\sigma)-\lambda)^{-1}}_{Y_{q,\mathrm{s}}(\sigma)}\leq C(q)|\lambda+\frac{\gamma_1}{2}|^{-1},
\text{ for all }|\sigma|\leq\gamma_0\text{ and }\lambda
\in  S(-\frac{\gamma_1}{2},\tilde{\omega}).$$
Thus, by \cite[Thm.1.3.4, 1.4.3]{henry}, we immediately obtain the two inequalities for $\widehat{A}_{\mathrm{s}}(\sigma)$. The
first inequality on $\widehat{A}_{\mathrm{ch}}(\sigma)$ follows directly by combining the first inequality for
$\widehat{A}_{\mathrm{c}}(\sigma)$ and the first inequality for $\widehat{A}_{\mathrm{s}}(\sigma)$.

We now derive the estimates for the case $\gamma_0<|\sigma|\leq \frac{1}{2}$.
By a similar analysis as in Proposition \ref{p:604}, there exists a positive constant $\gamma_2$ such that
$$\re(\spec \widehat{A}_{\mathrm{ch}}(\sigma))<-2\gamma_2, \text{ for all }\gamma_0<|\sigma|\leq \frac{1}{2}.$$
It is then not hard to conclude that
\beqs
\spec(\widehat{A}_{\mathrm{ch}}(\sigma))\subset \C\backslash S(-\gamma_2,\tilde{\omega}_1),\text{ where }
\tilde{\omega}_1\in(\frac{\pi}{2},\pi).
\eeqs
Moreover, for every $q\in[1, +\infty]$,
there exists a positive constant $C(q)$ such that
$$\opnorm{(\widehat{A}_{\mathrm{ch}}(\sigma)-\lambda)^{-1}}_{Y_q}\leq C(q)|\lambda+\gamma_2|^{-1},
\text{ for all }\gamma_0<|\sigma|\leq \frac{1}{2}\text{ and }\lambda
\in  S(-\gamma_2,\tilde{\omega}_1).$$
Therefore, again by \cite[Thm.1.3.4, 1.4.3]{henry}, we immediately obtain the last two inequalities for $\widehat{A}_{\mathrm{ch}}(\sigma)$,
which concludes the proof.
\epf

\bibliographystyle{siam}
\bibliography{myref}

\begin{thebibliography}{10}

\bibitem{bricmontkupiainen_1992}
{\sc J.~Bricmont and A.~Kupiainen}, {\em Renormalization group and the
  {G}inzburg-{L}andau equation}, Comm. Math. Phys., 150 (1992), pp.~193--208.

\bibitem{bricmontkupiainen_1994}
\leavevmode\vrule height 2pt depth -1.6pt width 23pt, {\em Stability of moving
  fronts in the {G}inzburg-{L}andau equation}, Comm. Math. Phys., 159 (1994),
  pp.~287--318.

\bibitem{denglevine_2000}
{\sc K.~Deng and H.~A. Levine}, {\em The role of critical exponents in blow-up
  theorems: the sequel}, J. Math. Anal. Appl., 243 (2000), pp.~85--126.

\bibitem{gallayscheel_2011}
{\sc T.~Gallay and A.~Scheel}, {\em Diffusive stability of oscillations in
  reaction-diffusion systems}, Trans. Amer. Math. Soc., 363 (2011),
  pp.~2571--2598.

\bibitem{henry}
{\sc D.~Henry}, {\em Geometric theory of semilinear parabolic equations},
  vol.~840 of Lecture Notes in Mathematics, Springer-Verlag, Berlin, 1981.

\bibitem{hervel_1993}
{\sc M.~A. Herrero and J.~J.~L. Vel{\'a}zquez}, {\em Some results on blow up
  for semilinear parabolic problems}, in Degenerate diffusions ({M}inneapolis,
  {MN}, 1991), vol.~47 of IMA Vol. Math. Appl., Springer, New York, 1993,
  pp.~105--125.

\bibitem{johnson_2009}
{\sc M.~A. Johnson}, {\em Nonlinear stability of periodic traveling wave
  solutions of the generalized {K}orteweg-de {V}ries equation}, SIAM J. Math.
  Anal., 41 (2009), pp.~1921--1947.

\bibitem{johnsonzumbrun_2010}
{\sc M.~A. Johnson and K.~Zumbrun}, {\em Nonlinear stability of periodic
  traveling wave solutions of systems of viscous conservation laws in the
  generic case}, J. Differential Equations, 249 (2010), pp.~1213--1240.

\bibitem{johnsonzumbrun_2011sj}
\leavevmode\vrule height 2pt depth -1.6pt width 23pt, {\em Nonlinear stability
  of periodic traveling-wave solutions of viscous conservation laws in
  dimensions one and two}, SIAM J. Appl. Dyn. Syst., 10 (2011), pp.~189--211.

\bibitem{johnsonzumbrun_2011}
\leavevmode\vrule height 2pt depth -1.6pt width 23pt, {\em Nonlinear stability
  of spatially-periodic traveling-wave solutions of systems of
  reaction-diffusion equations}, Ann. Inst. H. Poincar\'e Anal. Non Lin\'eaire,
  28 (2011), pp.~471--483.

\bibitem{johnsonzumbrunpascal_2011}
{\sc M.~A. Johnson, K.~Zumbrun, and P.~Noble}, {\em Nonlinear stability of
  viscous roll waves}, SIAM J. Math. Anal., 43 (2011), pp.~577--611.

\bibitem{lunardi}
{\sc A.~Lunardi}, {\em Analytic semigroups and optimal regularity in parabolic
  problems}, Progress in Nonlinear Differential Equations and their
  Applications, 16, Birkh\"auser Verlag, Basel, 1995.

\bibitem{mielke_1987}
{\sc A.~Mielke}, {\em \"{U}ber maximale {$L^p$}-{R}egularit\"at f\"ur
  {D}ifferentialgleichungen in {B}anach- und {H}ilbert-{R}\"aumen}, Math. Ann.,
  277 (1987), pp.~121--133.

\bibitem{mielke_1997}
\leavevmode\vrule height 2pt depth -1.6pt width 23pt, {\em Instability and
  stability of rolls in the {S}wift-{H}ohenberg equation}, Comm. Math. Phys.,
  189 (1997), pp.~829--853.

\bibitem{murray1}
{\sc J.~D. Murray}, {\em Mathematical biology. {I}}, vol.~17 of
  Interdisciplinary Applied Mathematics, Springer-Verlag, New York, third~ed.,
  2002.
\newblock An introduction.

\bibitem{murray2}
\leavevmode\vrule height 2pt depth -1.6pt width 23pt, {\em Mathematical
  biology. {II}}, vol.~18 of Interdisciplinary Applied Mathematics,
  Springer-Verlag, New York, third~ed., 2003.
\newblock Spatial models and biomedical applications.

\bibitem{reedsimon}
{\sc M.~Reed and B.~Simon}, {\em Methods of modern mathematical physics. {IV}.
  {A}nalysis of operators}, Academic Press [Harcourt Brace Jovanovich
  Publishers], New York, 1978.

\bibitem{sslsp_2001}
{\sc B.~Sandstede and A.~Scheel}, {\em On the stability of periodic travelling
  waves with large spatial period}, J. Differential Equations, 172 (2001),
  pp.~134--188.

\bibitem{sssu_2012}
{\sc B.~Sandstede, A.~Scheel, G.~Schneider, and H.~Uecker}, {\em Diffusive
  mixing of periodic wave trains in reaction-diffusion systems}, J.
  Differential Equations, 252 (2012), pp.~3541--3574.

\bibitem{scarpellini_1994}
{\sc B.~Scarpellini}, {\em $l^2$-perturbations of periodic equilibria of
  reaction diffusion systems}, NoDEA Nonlinear Differential Equations Appl., 1
  (1994), pp.~281--311.

\bibitem{schneider_1996}
{\sc G.~Schneider}, {\em Diffusive stability of spatial periodic solutions of
  the {S}wift-{H}ohenberg equation}, Comm. Math. Phys., 178 (1996),
  pp.~679--702.

\bibitem{turing_1952}
{\sc A.~Turing}, {\em The chemical basis of morphogenesis}, Philosophical
  Transactions of the Royal Society of London. Series B, Biological Sciences,
  237 (1952), pp.~37--72.

\bibitem{uecker_1999}
{\sc H.~Uecker}, {\em Diffusive stability of rolls in the two-dimensional real
  and complex {S}wift-{H}ohenberg equation}, Comm. Partial Differential
  Equations, 24 (1999), pp.~2109--2146.

\bibitem{zemi_2009}
{\sc S.~Zelik and A.~Mielke}, {\em Multi-pulse evolution and space-time chaos
  in dissipative systems}, Mem. Amer. Math. Soc., 198 (2009), pp.~vi+97.

\end{thebibliography}

\end{document}